\documentclass[reqno]{amsart}%
\usepackage{palatino, mathpazo}
\usepackage{amsfonts}
\usepackage{amsmath}
\usepackage{amssymb,latexsym,xcolor}
\usepackage{graphicx}
\usepackage{harpoon}
\usepackage[mathscr]{eucal}
\usepackage{amssymb}%
\usepackage[
linktocpage=true,colorlinks,citecolor=magenta,linkcolor=blue,urlcolor=magenta]{hyperref}
\hypersetup{
 colorlinks=true,
 linkcolor=blue,
 citecolor=blue,
 urlcolor=blue,
 pdftitle={Arbitrary Dirichlet Eigenvalue Multiplicities on Smooth Star-Shaped Euclidean Domains},
 pdfsubject={Dirichlet eigenvalue, multiplicity, star-shaped domain, spectral prescription, capacity},
 pdfkeywords={Dirichlet eigenvalues, multiplicity, star-shaped domains, spectral prescription, capacity}
}

\providecommand{\U}[1]{\protect \rule{.1in}{.1in}}

\usepackage{amsmath,amssymb,mathrsfs}

\usepackage{amsmath,amssymb,amsthm,mathtools}
\usepackage{esint}
\usepackage{microtype}
\usepackage{enumitem}
\usepackage{xurl}
\newtheorem{Theorem}{Theorem}

\numberwithin{equation}{section}
\newtheorem{theorem}{Theorem}[section]
\newtheorem{proposition}[theorem]{Proposition}
\newtheorem{lemma}[theorem]{Lemma}
\newtheorem{corollary}[theorem]{Corollary}
\theoremstyle{definition}

\usepackage{amsmath,amssymb}
\usepackage{graphicx}
\usepackage{tikz}
\usetikzlibrary{arrows.meta,calc}

\definecolor{sec31ink}{RGB}{38,59,76}
\definecolor{sec31room}{RGB}{232,241,248}
\definecolor{sec31tube}{RGB}{247,224,185}

\theoremstyle{remark}
\newtheorem{remark}[theorem]{Remark}
\newcommand{\R}{\mathbb R}
\newcommand{\Sph}{\mathbb S}

\newcommand{\dd}{\,dz}
\newcommand{\ds}{\,d\sigma}
\newcommand{\avint}{\mathop{\fint}}

\newcommand{\Sym}{\operatorname{Sym}}

\newcommand{\dist}{\operatorname{dist}}

\newcommand{\eps}{\varepsilon}
\newcommand{\norm}[1]{\left\|#1\right\|}
\allowdisplaybreaks[1]
\setlist[enumerate]{itemsep=3pt,topsep=5pt}
\title[Dirichlet multiplicities on star-shaped domains]{Arbitrary Dirichlet Eigenvalue Multiplicities on Smooth Star-Shaped Euclidean Domains}
\author{Zhijie Chen}
\address{Department of Mathematical Sciences, Yau Mathematical Sciences Center,
Tsinghua University, Beijing, 100084, China }
\email{zjchen2016@tsinghua.edu.cn}

\author{
Chen-Gang Long} 
\address{Yau Mathematical Sciences Center,
Tsinghua University, Beijing, 100084, China}
\email{lcg21@tsinghua.org.cn}

\author{Shihong Zhang}
\address{Yau Mathematical Sciences Center,
Tsinghua University, Beijing, 100084, China}
\email{shihong-zhang@tsinghua.edu.cn}
\date{}

\begin{document}
\maketitle
\begin{abstract}
In every fixed dimension \(N\ge3\), for any finite
list of positive integers \(m_1,\ldots,m_\ell\),  we construct smooth strictly
star-shaped Euclidean domains whose first
\(\ell+1\) distinct Dirichlet eigenvalues can have multiplicities
\(1,m_1,\ldots,m_\ell\) respectively. The condition "star-shaped" is almost optimal in the sense that this statement can not hold for bounded convex domains.
In fact, we show that the $k$-th Dirichlet eigenvalue
of any bounded convex domain
has multiplicity at most $(N-1)(k-1)+1$, and this bound $N$  for the second Dirichlet eigenvalue is sharp.
\end{abstract}

\section{Introduction and main results}\label{sec:intro}

For a bounded connected open set $\Omega\subset\R^N$, let
\[
0<\lambda_1(\Omega)<\lambda_2(\Omega)
\le\lambda_3(\Omega)\le\cdots
\]
be the eigenvalues of the Dirichlet Laplacian, counted with
multiplicity. For each $k\ge1$, write
\[
m_k(\Omega)
=\dim\ker\bigl(-\Delta_\Omega^D-\lambda_k(\Omega)\bigr).
\]

A natural question in spectral geometry is whether a finite initial
spectrum, including its multiplicities, can be prescribed by varying
the domain or the metric. For the Neumann Laplacian, Colin de
Verdi\`ere \cite{CdV} answered this question affirmatively for Euclidean
domains of dimension at least three.

\begin{Theorem}\label{thm:cdv-neumann}
\cite[Theorem~1.4]{CdV}
Let $N\ge3$ and $m\ge1$ be integers. For every finite sequence
\[
0=a_1<a_2\le\cdots\le a_{m+1},
\]
there exists a bounded connected open set $\Omega\subset\R^N$
with piecewise $C^1$ boundary such that
\[
\lambda_j^{\mathrm N}(\Omega)=a_j,
\qquad 1\le j\le m+1,
\]
where $\lambda_j^{\mathrm N}(\Omega)$ are the eigenvalues of
the Neumann Laplacian, arranged in nondecreasing order and
counted with multiplicity.
\end{Theorem}

Finite spectral prescription is also possible for the Dirichlet
problem when the Riemannian metric is allowed to vary. Recently, 
He and Wang \cite{HW2024} proved that, on any fixed compact smooth
manifold of dimension $N\ge3$ with nonempty smooth boundary,
a finite initial Dirichlet spectrum and the volume can be
prescribed simultaneously by a suitable choice of metrics.

For Euclidean domains, however, universal eigenvalue inequalities obstruct such an arbitrary prescription of the Dirichlet spectrum. There is a substantial literature on this topic, and we mention only a few representative results.
Indeed, let $\mathbb B$ be the unit ball in $\R^N$. The Faber--Krahn inequality
(see, e.g., \cite{BrascoDePhilippisVelichkov2015}) and the
Ashbaugh--Benguria inequality \cite{AshbaughBenguria1992} give
\[
\lambda_1(\Omega)|\Omega|^{2/N}
\ge \lambda_1(\mathbb{B})|\mathbb{B}|^{2/N},
\qquad
\frac{\lambda_2(\Omega)}{\lambda_1(\Omega)}
\le \frac{\lambda_2(\mathbb{B})}{\lambda_1(\mathbb{B})}.
\]
The first inequality prevents arbitrary prescription of
$\lambda_1(\Omega)$ at fixed volume, while the second rules out arbitrary
prescription of the first two eigenvalues even without a volume
constraint. 
For a broader discussion of extremal problems for
Dirichlet eigenvalues on Euclidean domains, see
F. Lin \cite{FanghuaLinSurvey}.
   More generally, for lower bounds on the eigenvalues, the Li--Yau inequality \cite{LiYau1983}
yields, for every $k\ge1$,
\[
\begin{aligned}
\sum_{i=1}^k\lambda_i(\Omega)
\ge \frac{4\pi^2 N}{(N+2)\bigl(|\mathbb{B}||\Omega|\bigr)^{2/N} }\,
k^{1+2/N}.
\end{aligned}
\]
On the other hand, the upper bound of Cheng and Yang \cite{ChengYang2007} gives
\begin{align*}
    \lambda_{k+1}(\Omega)
    \le \left(1+\frac4N\right)\frac1k
    \sum_{i=1}^k\lambda_i(\Omega).
\end{align*}
Recently, Funano \cite{Funano} obtained further universal eigenvalue inequalities for convex domains.
These estimates indicate that  the Dirichlet eigenvalues for bounded domains in $\R^N$ cannot be prescribed arbitrary like  He-Wang 's result \cite{HW2024}. 

Therefore, the remaining natural question is whether the multiplicity $m_k(\Omega)$ can be prescribed arbitrarily. It is interesting that the answer to this problem is negative for the dimension $N=2$, because Nadirashvili \cite{N} already proved in 1988 that $m_k(\Omega)\leq 2k-1$ for any $k$, provided $\Omega\subset \mathbb R^2$ is a bounded domain. See also \cite{HMN}. To the best of our knowledge, the general case $N\geq 3$ seems still open so far.

  In this paper, we give an affirmative answer to this problem for $N\geq 3$ by constructing star-shaped domains. More precisely, we consider the following type of strict star-shapedness:
\begin{equation}\label{eq:strictstar}
 \Omega=\{r\theta:\theta\in\Sph^{N-1},\ 0\le r<R(\theta)\},
 \qquad R\in C^\infty(\Sph^{N-1}),\qquad R>0.
\end{equation}
The outward unit normal then satisfies
\begin{equation}\label{eq:normal}
 x\cdot\nu(x)=
 \frac{R(\theta)^2}
 {\sqrt{R(\theta)^2+|\nabla_{\Sph^{N-1}}R(\theta)|^2}}>0,
 \qquad x=R(\theta)\theta.
\end{equation}

Let $\mathcal D_{N}$ denote the class of domains of the form
\eqref{eq:strictstar} whose closures are
$C^\infty$-diffeomorphic to the closed unit ball. The main result of this paper
prescribes an arbitrary finite list of multiplicities within this class.

\begin{theorem}
\label{thm:finite}
Let \(N\ge3\), and let
\(m_1,\ldots,m_\ell\) be positive integers, where \(\ell\ge1\).
Set
\begin{equation*}
 s_0=1,\qquad s_j=1+\sum_{i=1}^j m_i
 \quad 1\le j\le\ell.
\end{equation*}
Then there exist a domain \(\Omega\in\mathcal D_{N}\) and real numbers
\(\mu_1,\ldots,\mu_\ell\) such that
\begin{align} \label{eq:finiteblocks}
\begin{cases}
  \displaystyle  0<\lambda_1(\Omega)
 <\mu_1<\cdots<\mu_\ell<\lambda_{s_\ell+1}(\Omega),
\\
 \displaystyle\lambda_k(\Omega)=\mu_j \;\;\text{ for }\;\; s_{j-1}+1\le k\le s_j,\quad\forall\, 1\le j\le\ell.
\end{cases}
\end{align}
In particular,
\begin{equation}\label{eq:finitemultiplicities}
 \dim\ker(-\Delta_\Omega^D-\mu_j)=m_j,
 \quad(1\le j\le\ell).
\end{equation}
\end{theorem}

Theorem~\ref{thm:finite} prescribes multiplicities, not the numerical
values of the Dirichlet eigenvalues. In 2016, 
Hassannezhad, Kokarev and Polterovich
\cite[Open Problem~3, p.~817]{HKP}
asked whether Dirichlet eigenvalue multiplicities on Euclidean
domains admit an upper bound depending only on the dimension and the
eigenvalue index. The consequence of Theorem~\ref{thm:finite} gives a negative
answer to \cite[Open Problem~3]{HKP} at every fixed index \(k\ge2\), already in the class
\eqref{eq:strictstar}. 

\begin{corollary}
\label{cor:fixedindex}\label{cor:volumetopology}
For every $N\ge3$ and integer $k\ge2$,
\[
\sup_{\Omega\in\mathcal D_{N}}m_k(\Omega)=+\infty.
\]
\end{corollary}
The construction underlying Theorem~\ref{thm:cdv-neumann} of Colin de Verdi\`ere \cite{CdV}
approximates weighted graph Laplacians by rooms connected
through narrow passages, whose capacities encode the edge
weights. Further developments of this approach in spectral prescription
can be found in the works of Lohkamp \cite{Lohkamp1996}
and Jammes \cite{Jammes2011,Jammes2014}.   Our proof of Theorem~\ref{thm:finite} adapts this
graph approximation method to the Dirichlet problem.  We encode the graph weights in
densities on a fixed configuration of rooms and tubes, pass to a
Schr\"odinger operator on a ball, and use radial cuts with suitably
prescribed capacities to obtain smooth strictly star-shaped domains
with the desired Dirichlet eigenvalue multiplicities.

Next, we consider the more special case:  the multiplicity of Dirichlet eigenvalues on convex domains. For planar convex domains, Lin \cite{Lin} proved in 1987 that $m_2(\Omega)\leq 2$. For convex domains in higher dimensions, Funano \cite[Corollary 1.3]{Funano} recently proved that
\begin{equation}\label{eq:Funano}
m_k(\Omega)\le C_Nk
\qquad\text{provided that}\qquad
\lambda_k(\Omega)>\alpha_N\lambda_1(\Omega),
\end{equation}
where
\[
\alpha_N=\frac{16N^3}{\pi^2}j_{\frac N2-1,1}^{\,2},
\qquad
C_N=\left(12N^3j_{\frac N2-1,1}^{\,2}\right)^{N/2},
\]
and $j_{\nu,1}$ denotes the first positive zero of the
Bessel function $J_\nu$.
In a complementary direction, Frank and Larson
\cite{FrankLarsonConvex} established semiclassical
spectral inequalities for Dirichlet and Neumann
Laplacians on convex Euclidean domains.
  Our second theorem gives a new upper bound for $m_k(\Omega)$ in arbitrary dimensions, thereby generalizing Lin's theorem \cite{Lin} and  improving Funano's \cite[Corollary 1.3]{Funano}.

\begin{theorem}\label{thm:convex-all-multiplicities}
Let $\Omega\subset\mathbb{R}^N$, $N\ge2$, be a nonempty
bounded convex domain, and let $D=\operatorname{diam}(\Omega)$.
Then, for every integer $k\ge2$,
\begin{equation}\label{eq:convex-all-multiplicity-bound}
m_k(\Omega)\le (N-1)(k-1)+1,
\end{equation}
and
\begin{equation}\label{eq:convex-quantitative-separation}
\lambda_{N(k-1)+2}(\Omega)-\lambda_k(\Omega)
\ge \frac{2\pi^2}{D^2}.
\end{equation}
The multiplicity bound \eqref{eq:convex-all-multiplicity-bound} for $k=2$ is $m_2(\Omega)\leq N$, which is
is sharp by taking $\Omega=\mathbb B$.
\end{theorem}

\begin{remark} Theorem \ref{thm:convex-all-multiplicities} indicates that the condition "star-shaped" in Theorem~\ref{thm:finite} is almost optimal.
Clearly the estimate~\eqref{eq:convex-all-multiplicity-bound} improves
Funano's result \eqref{eq:Funano} in two respects: It removes the
restriction $\lambda_k(\Omega)>\alpha_N\lambda_1(\Omega)$ and gives the smaller upper bound
\[
m_k(\Omega)\le (N-1)(k-1)+1\le (N-1)k
\qquad\text{for every }k\ge2.
\]
In particular, the coefficient in the linear bound is reduced
from $C_N$ to $N-1$.
\end{remark}

The organization of this paper is as follows. In Section~\ref{sec:overview}, we begin by reviewing Colin de Verdi\`ere's construction and outlining our approach. In Section~\ref{sec:graph}, we construct a weighted graph model on a fixed configuration of rooms and tubes, and derive the corresponding energy estimate. In Section~\ref{sec:potential}, we establish the spectral reduction for the weighted Dirichlet problem on the unit ball, and then transform it into a family of Schr\"odinger operators with smooth potentials. In Section~\ref{sec:radialconstruction}, we use radial cuts and boundary smoothing to construct smooth strictly star-shaped domains. Their Dirichlet Laplacians approximate these Schr\"odinger operators after a common constant shift. In Section~\ref{sec:finish}, we apply Brouwer's fixed point theorem to obtain the prescribed multiplicities and prove Theorem~\ref{thm:finite}. Finally, in Section~\ref{sec:convexproof}, we prove the multiplicity bound for the $k$-th Dirichlet eigenvalue on convex domains, thereby proving Theorem~\ref{thm:convex-all-multiplicities}. The auxiliary spectral perturbation result is collected in Appendix~\ref{sec:spectral}.

\textbf{Notation}: Let $N\geq 3$ and $\mathbb{B}_r:=\left\{x\in \R^{N}:|x|<r\right\}$ be the ball in $\R^{N}$, and let  $\mathbb{B}=\mathbb{B}_1$. Let $m\geq 1$ and
$
\mathbf{1}=(1,\ldots,1)^T\in\mathbb{R}^{m+1}.
$
Its orthogonal complement is given by
\begin{align*}
    H:=\mathbf{1}^{\perp}
    =\left\{x\in\mathbb{R}^{m+1}:\sum_{i=1}^{m+1}x_i=0\right\}.
\end{align*}
Then $\dim H=m$. We write $\Sym(H)$ for the space of self-adjoint operators on the finite-dimensional real inner product space $H$. Fix an orthonormal basis $h_1,\ldots,h_m$ of $H$ as
\begin{equation}\label{eq:hk-def}
h_k:=\frac{e_1+\cdots+e_k-ke_{k+1}}{\sqrt{k(k+1)}},
\qquad 1\le k\le m,
\end{equation}
where $e_1,\ldots,e_{m+1}$ are the standard Euclidean basis of
$\R^{m+1}$. 
Throughout, we identify $\Sym(H)$ with the space of real symmetric
$m\times m$ matrices through
\[
p\longleftrightarrow
\bigl(\langle ph_l,h_k\rangle\bigr)_{k,l=1}^m,
\]
and equip this space with the operator norm
$$\|p\|:=\max_{x\in H\setminus\{0\}}\frac{|p\, x|}{|x|}.$$
Throughout, $x,y\in\mathbb{R}^{m+1}$ denote fixed vectors, while $z\in\mathbb{R}^{N}$ denotes the variable of integration.

\section{Colin de Verdi\`ere's approach and ideas of our proof}
\label{sec:overview}

In this section, we will explain  the Neumann approximation method of Colin de Verdi\`ere
\cite{CdV} (see also \cite{CdV1986,CdV1988,ColboisCdV1988}) and illustrate our main modifications of his approach. Roughly speaking, his main idea is to transform the eigenvalue problem of $-\Delta_{\Omega}^{N}$ into the finite dimension case, namely, the eigenvalues of matrices. To achieve this, he needed to start on the finite graph to construct desired small eigenvalues with high order perturbation. 

\subsection{Colin de Verdi\`ere's construction}
\begin{figure}[htbp]
\centering
\resizebox{0.75\linewidth}{!}{%
\begin{tikzpicture}[
    x=1.12cm,y=1.12cm,
    font=\small,
    >=Stealth,
    line cap=round,line join=round,
    wall/.style={draw=black!80,line width=0.9pt}
]
\definecolor{roomfill}{RGB}{231,239,249}
\definecolor{tubefill}{RGB}{251,225,180}
\definecolor{profileblue}{RGB}{42,91,153}

\path[fill=roomfill]
(-2,1.3)
.. controls (-3,1.6) and (-4.4,1.6) .. (-5,0.8)
.. controls (-5.3,0) and (-5,-1.2) .. (-4.2,-1.4)
.. controls (-3.4,-1.6) and (-2.7,-1.5) .. (-2,-1.3)
-- cycle;

\path[fill=roomfill]
(2,1.3)
.. controls (3,1.6) and (4.4,1.6) .. (5,0.8)
.. controls (5.3,0) and (5,-1.2) .. (4.2,-1.4)
.. controls (3.4,-1.6) and (2.7,-1.5) .. (2,-1.3)
-- cycle;

\path[fill=tubefill]
(-2,-0.16) rectangle (2,0.16);

\draw[wall]
(-2,0.16)--(-2,1.3)
.. controls (-3,1.6) and (-4.4,1.6) .. (-5,0.8)
.. controls (-5.3,0) and (-5,-1.2) .. (-4.2,-1.4)
.. controls (-3.4,-1.6) and (-2.7,-1.5) .. (-2,-1.3)
--(-2,-0.16);

\draw[wall]
(2,0.16)--(2,1.3)
.. controls (3,1.6) and (4.4,1.6) .. (5,0.8)
.. controls (5.3,0) and (5,-1.2) .. (4.2,-1.4)
.. controls (3.4,-1.6) and (2.7,-1.5) .. (2,-1.3)
--(2,-0.16);

\draw[wall] (-2,0.16)--(2,0.16);
\draw[wall] (-2,-0.16)--(2,-0.16);

\node[align=center] at (-4,0)
{$\widehat R_i$\\[3pt]$x_i/\sqrt v$};

\node[align=center] at (4,0)
{$\widehat R_j$\\[3pt]$x_j/\sqrt v$};

\node[align=center] at (0,0.72)
{$T_{e,\tau}$\\[-1pt]\footnotesize radius $\eta_{e,\tau}$};

\draw[<->,black!65]
(-2,-0.58)--(2,-0.58)
node[midway,below] {$L_e$};

\node at (0,1.82)
{Two rooms connected by a thin tube};

\begin{scope}[yshift=-3.6cm]
\node at (0,1.35)
{Trial profile along the connection};

\draw[->,black!50]
(-4.8,-0.28)--(5,-0.28)
node[right] {$s$};

\draw[black!35,densely dashed]
(-2,-0.28)--(-2,0.9);
\draw[black!35,densely dashed]
(2,-0.28)--(2,0);

\draw[profileblue,line width=1.3pt]
(-4.5,0.9)--(-2,0.9)--(2,0)--(4.5,0);

\node[text=profileblue,above] at (-3.7,0.9)
{$x_i/\sqrt v$};
\node[text=profileblue,above] at (3.7,0)
{$x_j/\sqrt v$};

\node[below] at (-2,-0.28) {$0$};
\node[below] at (2,-0.28) {$L_e$};
\node[text=profileblue,above] at (0,0.72)
{affine in $s$};
\end{scope}

\end{tikzpicture}%
}
\caption{\small{Two rooms joined by a thin tube and the corresponding
trial profile.}}
\label{fig:cdv-connector}
\end{figure}
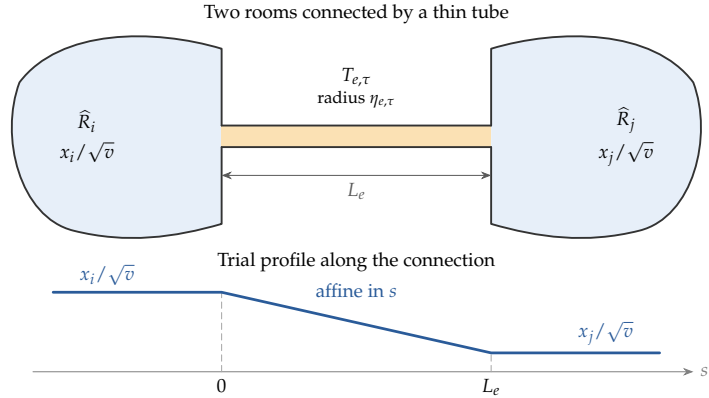

 As in the introduction, $x,y\in\mathbb R^{m+1}$
denote room values and $z\in\mathbb R^N$ is the spatial variable.
We write $\lambda_j^{\mathrm N}$ for Neumann eigenvalues and reserve
$\lambda_j$ for Dirichlet eigenvalues. In the Neumann case, the
rooms are denoted by $\widehat R_i$, and the  connections are indexed by
$e=\{i,j\}$, see Figure~\ref{fig:cdv-connector}.

\textbf{Two rooms.}
Two  connected rooms provide the simplest example:
their interaction produces a small eigenvalue. Let us explain this in detail.

Let $\widehat{R}_1,\widehat{R}_2$ be disjoint connected rooms of common volume $v>0$.
On $\Omega_0=\widehat{R}_1\cup \widehat{R}_2$, given $(x_1,x_2)\in\mathbb R^2$, the function $u_x: \Omega_0\to\mathbb R$ of the form
\[
u_x|_{\widehat{R}_i}=\frac{x_i}{\sqrt v},
\qquad i=1,2,
\]
has zero Dirichlet energy. The zero eigenspace is two-dimensional, so
\[
\lambda_1^{\mathrm N}(\Omega_0)
=
\lambda_2^{\mathrm N}(\Omega_0)
=
0.
\]

Denote $e=\{1,2\}$ and given $c_e=c_{12}>0$. For $\tau>0$ small, connect the rooms by a cylindrical tube
$T_{e,\tau}$ of fixed length $L_e$. In coordinates along the
tube, write
\[
T_{e,\tau}=(0,L_e)\times\Sigma_{e,\tau},
\qquad
\Sigma_{e,\tau}
=\{z'\in\mathbb R^{N-1}:|z'|<\eta_{e,\tau}\},
\]
and $\Omega_{\tau}=\Omega_0\cup T_{e,\tau}$. Choose its radius $\eta_{e,\tau}=O(\tau^{\frac{1}{N-1}})$ so that
\[
\frac{|\Sigma_{e,\tau}|}{L_e}
=\tau v c_{e}.
\]
Taking
\[
u=\frac1{\sqrt v}
\quad\text{in }\widehat{R}_1,
\qquad
u=-\frac1{\sqrt v}
\quad\text{in }\widehat{R}_2,
\]
and
\[
u(s)=\frac1{\sqrt v}
\left(1-\frac{2s}{L_e}\right),
\qquad 0\le s\le L_e,
\]
in the tube, we obtain
\[
\int_{\Omega_\tau}u\,dz=0,
\qquad
\int_{\Omega_\tau}u^2\,dz\ge2,
\qquad
\int_{\Omega_\tau}|\nabla u|^2\,dz
=\frac{4|\Sigma_{e,\tau}|}{vL_e}=
4\tau c_{12}.
\]
The variational characterization therefore gives
\[
0<
\lambda_2^{\mathrm N}(\Omega_\tau)
\le
2\tau c_{12}.
\]

To obtain the matching lower bound, let $u_\tau$ be a real
eigenfunction corresponding to $\lambda_2^{\mathrm N}(\Omega_\tau)$,
normalized by
\[
\int_{\Omega_\tau}u_\tau^2\,dz=1,
\qquad
\int_{\Omega_\tau}u_\tau\,dz=0.
\]
The preceding upper bound implies
\[
\int_{\Omega_\tau}|\nabla u_\tau|^2\,dz
=\lambda_2^{\mathrm N}(\Omega_\tau)
=O(\tau).
\]
Define the room averages by
\[
x_{i,\tau}
:=\frac1{\sqrt v}\int_{\widehat R_i}u_\tau\,dz,
\qquad
a_{i,\tau}:=\frac{x_{i,\tau}}{\sqrt v},
\qquad i=1,2.
\]
Then \(|x_{i,\tau}|\leq1\), and the Poincar\'e inequality on
each fixed room gives
\[
\int_{\widehat R_i}|u_\tau-a_{i,\tau}|^2\,dz
\leq C\int_{\widehat R_i}|\nabla u_\tau|^2\,dz
=O(\tau).
\]

Let \(\Gamma_{i,\tau}\) be the interface between
\(\widehat R_i\) and the tube.
Since \(N\geq3\), the Sobolev trace inequality with
\(q=2(N-1)/(N-2)\) gives
\[
\|u_\tau-a_{i,\tau}\|_{L^q(\partial\widehat R_i)}^2
\leq C\int_{\widehat R_i}|\nabla u_\tau|^2\,dz=O(\tau).
\]
Since \(|\Gamma_{i,\tau}|\sim\eta_{e,\tau}^{N-1}\), H\"older's
inequality implies
\begin{align}\label{Sec 2 equ-s}
\int_{\Gamma_{i,\tau}}|u_\tau-a_{i,\tau}|^2\,dz'
\leq C\eta_{e,\tau}
\int_{\widehat R_i}|\nabla u_\tau|^2\,dz.
\end{align}
Since \(|\Sigma_{e,\tau}|/L_e=\tau vc_{12}\), we have
\(\eta_{e,\tau}\sim\tau^{1/(N-1)}\) and \(|\Sigma_{e,\tau}|=O(\tau)\).
In particular, \(a_{1,\tau}=O(1)\) gives
\[
\int_{\Gamma_{1,\tau}}u_\tau^2\,dz'
\leq 2|\Sigma_{e,\tau}|a_{1,\tau}^2
+2\int_{\Gamma_{1,\tau}}|u_\tau-a_{1,\tau}|^2\,dz'
=O(\tau).
\]
Integration along the tube and Cauchy--Schwarz inequality therefore yield
\[
\int_{T_{e,\tau}}u_\tau^2\,dz
\leq
2L_e\int_{\Gamma_{1,\tau}}u_\tau^2\,dz'
+2L_e^2\int_{T_{e,\tau}}|\partial_su_\tau|^2\,dz
=O(\tau).
\]
Using the definition of the room averages, we obtain
\[
1=\int_{\Omega_\tau}u_\tau^2\,dz=\sum_{i=1}^2
\left(x_{i,\tau}^2+
\int_{\widehat R_i}|u_\tau-a_{i,\tau}|^2\,dz\right)
+\int_{T_{e,\tau}}u_\tau^2\,dz,
\]
and hence
\[
x_{1,\tau}^2+x_{2,\tau}^2=1+O(\tau).
\]
Moreover, \(|T_{e,\tau}|=O(\tau)\), so the zero-mean condition gives
\[
\sqrt v\,|x_{1,\tau}+x_{2,\tau}|
=\left|\int_{T_{e,\tau}}u_\tau\,dz\right|
\leq |T_{e,\tau}|^{1/2}
\|u_\tau\|_{L^2(T_{e,\tau})}
=O(\tau).
\]
Consequently,
\[
(x_{1,\tau}-x_{2,\tau})^2
=2(x_{1,\tau}^2+x_{2,\tau}^2)
-(x_{1,\tau}+x_{2,\tau})^2
=2+O(\tau).
\]

To estimate the tube energy, define the cross-sectional average
\[
w_\tau(s):=\frac1{|\Sigma_{e,\tau}|}
\int_{\Sigma_{e,\tau}}u_\tau(s,z')\,dz',
\qquad 0<s<L_e.
\]
Since $w_\tau(0)-a_{1,\tau}=\frac1{|\Sigma_{e,\tau}|}
\int_{\Gamma_{1,\tau}}(u_\tau-a_{1,\tau})\,dz'$ and $w_\tau(L_e)-a_{2,\tau}=\frac1{|\Sigma_{e,\tau}|}
\int_{\Gamma_{2,\tau}}(u_\tau-a_{2,\tau})\,dz'$,
the interface estimates \eqref{Sec 2 equ-s} imply
\[
|w_\tau(0)-a_{1,\tau}|^2
+|w_\tau(L_e)-a_{2,\tau}|^2\leq C\eta_{e,\tau}^{2-N}\sum_{i=1}^2\int_{\widehat R_i}|\nabla u_\tau|^2\,dz=o(1).
\]
Since \(x_{i,\tau}\) are bounded, it follows that
\[
|w_\tau(L_e)-w_\tau(0)|^2
=\frac{(x_{1,\tau}-x_{2,\tau})^2}{v}+o(1).
\]
Applying the Cauchy--Schwarz inequality first on each cross-section and then
along the tube, we obtain
\begin{align}\label{Sec 2 energy}
\frac{\lambda_2^{\mathrm N}(\Omega_\tau)}{\tau}
&\geq \frac{|\Sigma_{e,\tau}|}{\tau}
\int_0^{L_e}|w_\tau'(s)|^2\,ds
\geq \frac{|\Sigma_{e,\tau}|}{\tau L_e}
|w_\tau(L_e)-w_\tau(0)|^2 \nonumber\\
&=c_{12}(x_{1,\tau}-x_{2,\tau})^2+o(1)
=2c_{12}+o(1).
\end{align}
Together with the upper bound, this proves
\[
\lambda_2^{\mathrm N}(\Omega_\tau)
=2\tau c_{12}+o(\tau),
\qquad \tau\downarrow0.
\]

\textbf{Multiple rooms.}
The same construction can be applied to $m+1\geq 2$ fixed, disjoint,
bounded connected Lipschitz rooms $\widehat R_1,\dots,\widehat R_{m+1}$,
each of volume $v$. Here we only sketch the construction without details. See \cite{CdV} for rigorous proofs.

For every pair $e=\{i,j\}$ with $1\leq i<j\leq m+1$, join the two rooms
by a cylindrical tube
\[
T_{e,\tau}=(0,L_e)\times\Sigma_{e,\tau},
\qquad
\frac{|\Sigma_{e,\tau}|}{L_e}=\tau v c_e,
\qquad c_e=c_{ij}>0.
\]
The lengths $L_e$ are fixed. The rooms and tubes are arranged so
that the tubes meet flat portions of the room boundaries and are
pairwise disjoint outside the rooms for sufficiently small $\tau$.
Let $\Omega_\tau$ be the resulting connected domain.
For $x=(x_1,\dots,x_{m+1})\in\mathbb R^{m+1}$, define $J_\tau x:\Omega_\tau\to\mathbb R$ by
\[
\begin{aligned}
(J_\tau x)|_{\widehat R_i}
=\frac{x_i}{\sqrt v},\quad
(J_\tau x)(s,z')
=\frac{(1-s/L_e)x_i+(s/L_e)x_j}{\sqrt v},
\quad (s,z')\in T_{e,\tau}.
\end{aligned}
\]
Thus the two-room interpolation is applied separately to each tube,
with a common constant value on each room. Figure~\ref{fig:cdv-connector} illustrates one such connection.
The trial function is constant on each room and interpolates
linearly along the axial coordinate $s$ of the cylindrical tube
$T_{e,\tau}$.

To describe the combined interaction between the rooms, recall that $e_1,\ldots,e_{m+1}$ are the standard basis of $\R^{m+1}$.
For $1\le i<j\le m+1$, define a $(m+1)\times(m+1)$ matrix
\[
E_{ij}:=(e_i-e_j)(e_i-e_j)^T.
\]
Its only nonzero entries are
\[
(E_{ij})_{ii}=(E_{ij})_{jj}=1,
\qquad
(E_{ij})_{ij}=(E_{ij})_{ji}=-1.
\]
Then we associate each $c=(c_{ij})_{i<j}\in\R^{m(m+1)/2}$ with a $(m+1)\times(m+1)$ matrix
\begin{equation}\label{eq:Lc-def}
L_c:=\sum_{1\le i<j\le m+1}c_{ij}E_{ij}.
\end{equation}
With the convention $c_{ji}=c_{ij}$ for $i<j$, we have
\[
(L_c)_{ij}=-c_{ij}\quad(i\ne j),
\qquad
(L_c)_{ii}=\sum_{j\ne i}c_{ij}.
\]
In particular, $L_c\mathbf{1}=0$. Here, since $c_{ij}>0$, we can
write the eigenvalues of $L_c$ as
\[
0=\mu_1(c)<\mu_2(c)\le\cdots\le\mu_{m+1}(c).
\]
Summing the energy and mass contributions of the rooms and tubes gives
\[
\begin{aligned}
\int_{\Omega_\tau}|\nabla J_\tau x|^2\,dz
&=\tau\sum_{i<j}c_{ij}(x_i-x_j)^2
=\tau x^TL_cx,\\
\int_{\Omega_\tau}|J_\tau x|^2\,dz
&=|x|^2+O(\tau)|x|^2\ge |x|^2.
\end{aligned}
\]
For details of the above identity, see \eqref{eq:exactgraphenergy}. The min--max principle therefore yields
\[
\lambda_j^{\mathrm N}(\Omega_\tau)\le\tau\mu_j(c),
\qquad 1\le j\le m+1.
\]

For the reverse bound, let $u_\tau$ satisfy
$\|u_\tau\|_{L^2(\Omega_\tau)}=1$ and
$\int_{\Omega_\tau}|\nabla u_\tau|^2\,dz\le C\tau$, and set
\[
x_{i,\tau}
=\frac1{\sqrt v}\int_{\widehat R_i}u_\tau\,dz,
\qquad
x_\tau=(x_{1,\tau},\dots,x_{m+1,\tau}).
\]
Applying the two-room estimates to each connection gives, uniformly
for such functions,
\[
\sum_i\int_{\widehat R_i}
\left|u_\tau-\frac{x_{i,\tau}}{\sqrt v}\right|^2\,dz
=O(\tau),
\qquad
\sum_e\int_{T_{e,\tau}}|u_\tau|^2\,dz=o(1).
\]
Consequently, $
|x_\tau|^2=1+o(1).
$
Similarly, by the two-room estimate \eqref{Sec 2 energy}, we obtain
\[
\begin{aligned}
\frac1\tau\int_{\Omega_\tau}|\nabla u_\tau|^2\,dz
&\ge
\sum_{i<j}c_{ij}(x_{i,\tau}-x_{j,\tau})^2-o(1)\\
&=x_\tau^TL_cx_\tau-o(1).
\end{aligned}
\]
The room-average map is consequently injective on the span of the
first $j$ Neumann eigenfunctions when $\tau$ is sufficiently small.
Applying the matrix min--max principle to its image yields
\[
\mu_j(c)\le
\frac{\lambda_j^{\mathrm N}(\Omega_\tau)}{\tau}+o(1),
\qquad 1\le j\le m+1.
\]
Combining the two bounds, we obtain
\begin{equation}
\label{eq:overview-neumann-limit}
\lambda_j^{\mathrm N}(\Omega_\tau)
=\tau\mu_j(c)+o(\tau),
\qquad 1\le j\le m+1.
\end{equation}
The room Poincar\'e inequalities and the tube estimates also give
a uniform positive lower bound for
$\lambda_{m+2}^{\mathrm N}(\Omega_\tau)$. Recall $\mu_1(c)=0$. 
Thus the matrix $L_c|_H$ governs the $m$ small positive Neumann
eigenvalues, while the direction $\mathbf1$ corresponds to the
constant eigenfunction. See \cite{CdV} for more details.

\textbf{The difference between the Dirichlet problem and the Neumann problem}:
\begin{itemize}
\item In the Neumann case, since the working space is $H^1(\Omega)$, Colin de Verdi\`ere can construct trial functions directly on the finite graph. For example, he can set $u\equiv $ const on each room $\widehat{R_i}$. However, such a construction does not work in the Dirichlet case.

\item On the same cylinder $T_{e,\tau}$, a function $u\in H^1_0(\Omega_{\tau})$ satisfies
\[
\int_{T_{e,\tau}}|\nabla_{z'}u|^2\,dz
\ge C_N\eta_{e,\tau}^{-2}
\int_{T_{e,\tau}}u^2\,dz,
\qquad z=(s,z').
\]
Thus,  the transverse Rayleigh quotient diverges as
$\eta_{e,\tau}\to0$. 
\end{itemize}
Therefore, Colin de Verdi\`ere's construction can not be applied easily to the Dirichlet problem, and new ideas are required.

\subsection{Our construction in three steps}
In this paper,
we modify Colin de Verdi\`ere's room-value mechanism to the Dirichlet problem through
the following three new steps. 

\medskip
\noindent
\textbf{Step 1: A weighted Dirichlet problem on the unit ball.}
Fix the configuration $D\Subset\mathbb B$ of rooms $R_i$ and tubes
$T_e$ constructed in Section~\ref{sec:graph}. The tubes $T_e$ have fixed geometry independent
of $\tau$, and their coupling is controlled by the densities.
To apply the Brouwer's fixed point theorem (see Section \ref{sec:finish}), here we need to introduce a new parameter $p\in \Sym(H)$ to vary the edge weights.
By Lemma~\ref{lem:graphweights}, we choose positive edge weights
$c_{ij}(p)$ depending on $p$ such that
\[
L_{c(p)}|_H=(m+1)I_H+p,
\qquad
p\in K_r:=\{p\in\Sym(H):\|p\|\le r\}.
\]
We encode these weights in the density
\[
\widehat\rho_{\tau,p}(z)=
\begin{cases}
1, & z\in R_i,\\
\tau\alpha_e(p), & z\in T_e,\\
\tau^2, & z\in\mathbb B\setminus\overline D,
\end{cases}
\qquad
\alpha_e(p)=\frac{v c_e(p)}{\kappa_e},
\]
where $\kappa_e$ is the fixed connector capacity defined in
Section~\ref{sec:graph}. We then consider the weighted Dirichlet
eigenvalue problem
\[
\begin{cases}
-\operatorname{div}(\widehat\rho_{\tau,p}\nabla f)
=\nu\,\widehat\rho_{\tau,p}f
& \text{in }\mathbb B,\\
f=0 & \text{on }\partial\mathbb B,
\end{cases}
\]
understood in the weak sense.
The room and tube boundaries are internal interfaces, so the trial functions may be constant on the rooms and interpolate through the tubes, and can then be extended to zero away from $D$, see Section \ref{sec:graph} and Section \ref{sec:potential}.

The estimates in Section~\ref{sec:potential} give a first eigenvalue
of order at most $\tau^2$ and an isolated cluster of $m$ eigenvalues
of order $\tau$. On this cluster, the operator matrix satisfies
\[
\sup_{p\in K_r}
\left\|
\widehat B_{\tau,p}
-\tau\bigl((m+1)I_H+p\bigr)
\right\|
\le C\tau^2.
\]

\medskip
\noindent
\textbf{Step 2: Transformation to a Schr\"odinger operator.}
For fixed $\tau$, smooth the density $\widehat\rho_{\tau,p}$ to
$\widehat\rho_{\tau,p,\delta}$ and set
\[
\phi_p=\sqrt{\widehat\rho_{\tau,p,\delta}},
\qquad
V_p=\frac{\Delta\phi_p}{\phi_p},
\qquad
L_p=-\Delta_{\mathbb B}^D+V_p.
\]
The substitution $u=\phi_pf$ transforms the smoothed weighted
problem into
\[
\begin{cases}
-\Delta u+V_pu=\nu u & \text{in }\mathbb B,\\
u=0 & \text{on }\partial\mathbb B.
\end{cases}
\]
Indeed, multiplication by $\phi_p$ is unitary from the weighted
$L^2$ space to $L^2(\mathbb B)$ and satisfies
\[
\int_{\mathbb B}\widehat\rho_{\tau,p,\delta}|\nabla f|^2\,dz
=
\int_{\mathbb B}\bigl(|\nabla u|^2+V_pu^2\bigr)\,dz.
\]
For sufficiently small $\tau$ and then sufficiently small $\delta$,
the operator matrix $\mathcal B_p$ of $L_p=-\Delta_{\mathbb B}^D+V_p$ on the spectral subspace for
eigenvalues $2,\dots,m+1$ depends continuously on $p$ and satisfies
\[
\sup_{p\in K_r}
\left\|\mathcal B_p-\tau\bigl((m+1)I_H+p\bigr)\right\|
<\frac{\tau r}{4}.
\]

\medskip
\noindent
\medskip
\noindent
\textbf{Step 3: Approximation by strictly star-shaped domains.}
Choose a common constant $C_*$ so that $W_p=V_p+C_*$ is positive
and $r^2W_p(r\theta)$ increases with $r$. Approximating Schr\"odinger operators with nonnegative potentials
by Dirichlet Laplacians on perforated domains is a classical
technique based on capacity estimates; see, for instance,
\cite{CM,MM95,Jing,KP,RT}.
In Section~\ref{sec:radialconstruction}, radial cuts with suitably
prescribed capacities, followed by boundary smoothing, produce
smooth strictly star-shaped domains $\Omega_{k,p}$ approximating the spectrum of
$L_p+C_*$ uniformly in $p$:
\[\sup_{p\in K_r}
 \left\|
 E_{\Omega_{k,p}}(-\Delta_{\Omega_{k,p}}^D+1)^{-1}E_{\Omega_{k,p}}^*
 -(L_p+C_*+1)^{-1}
 \right\|\to0\quad(k\to+\infty).\]
Together with Step~2, this gives continuous spectral matrices
$\mathcal C_p$ satisfying
\[
\sup_{p\in K_r}
\left\|
\mathcal C_p-(C_*+\tau(m+1))I_H-\tau p
\right\|<\frac{\tau r}{2}.
\]
Brouwer's fixed point theorem then yields, for every
$A\in\Sym(H)$ with $\|A\|\le r/4$, the existence of a parameter $p_*\in K_r$ such that
\begin{equation}
\label{eq:overview-exact-matrix}
\mathcal C_{p_*}
=(C_*+\tau(m+1))I_H+\tau A;
\end{equation}
see Section~\ref{sec:finish}.
Choosing $A$ with distinct scalar blocks of dimensions
$m_1,\dots,m_\ell$, where $\sum_jm_j=m$, gives the prescribed
multiplicities. 

\section{Weighted Graph Laplacians and Uniform Estimates}
\label{sec:graph}

Fix $N\ge3$ and $m\geq 1$. The rooms and tubes lie in
$\mathbb{B}\subset\R^N$, whereas $\R^{m+1}$ is the space of values
assigned to the $m+1$ vertices. Recall the $(m+1)\times(m+1)$ symmetric matrix $L_c$ defined in \eqref{eq:Lc-def} for $c=(c_{ij})_{i<j}\in\R^{m(m+1)/2}$, and recall that
\begin{align*}
     H:=\mathbf{1}^{\perp}
    =\left\{x\in\mathbb{R}^{m+1}:\sum_{i=1}^{m+1}x_i=0\right\}.
\end{align*}
Since $L_c^T=L_c$ and $L_c\mathbf{1}=0$, we have
$L_c(H)\subset H$ and $L_c|_H\in\Sym(H)$.

We first construct positive edge weights such that the graph
Laplacian on $H$ equals $(m+1)I_H+p$ for every sufficiently small
$p\in\Sym(H)$.

\begin{lemma}\label{lem:graphweights}
There exist $0<r<(m+1)/4$ and an affine family of positive weights
\[
c(p)=(c_{ij}(p))_{i<j},
\qquad
p\in K_r:=\{p\in\Sym(H):\|p\|\le r\},
\]
such that
\[
L_{c(p)}|_H=(m+1)I_H+p,
\qquad
\frac12\le c_{ij}(p)\le\frac32
\quad(1\le i<j\le m+1).
\]
\end{lemma}

\begin{proof}
Consider the linear map
\[
\mathcal L:\R^{m(m+1)/2}\longrightarrow\Sym(H),
\qquad
\mathcal L(a):=L_a|_H.
\]
If $\mathcal L(a)=0$, then $L_a|_H=0$ and
$L_a\mathbf{1}=0$. Since
$\R^{m+1}=H\oplus\R\mathbf{1}$, it follows that $L_a=0$,
and hence $a_{ij}=0$ for every $i<j$. Thus $\mathcal L$ is
injective. Moreover,
\[
\dim\R^{m(m+1)/2}
=\frac{m(m+1)}2
=\dim\Sym(H),
\]
so $\mathcal L$ is a linear isomorphism. Let $\mathbf c^{\,0}$ be given by $c^0_{ij}=1$. Then
\[
L_{\mathbf c^{\,0}}
=(m+1)I_{m+1}-\mathbf{1}\mathbf{1}^T,
\qquad
L_{\mathbf c^{\,0}}|_H=(m+1)I_H.
\]
Let $C_{\mathcal L}$ be the operator norm of $\mathcal L^{-1}$,
where $\Sym(H)$ carries the operator norm and
$\R^{m(m+1)/2}$ carries the maximum norm. Choose
\[
0<r<\min\left\{\frac{m+1}{4},\frac1{2C_{\mathcal L}}\right\},
\qquad
c(p):=\mathbf c^{\,0}+\mathcal L^{-1}(p).
\]
Then for $p\in K_r$,
\[
\max_{i<j}|c_{ij}(p)-1|
\le C_{\mathcal L}\|p\|<\frac12,
\qquad
L_{c(p)}|_H=(m+1)I_H+p.
\]
This proves the lemma.
\end{proof}

\subsection{Rooms, Tubes, and Parameter-Dependent Densities}
\label{sec:graphgeometry}

We construct a fixed connected domain from $m+1$ congruent balls
and pairwise connecting tubes.
Choose distinct real numbers $t_1,\ldots,t_{m+1}$ and set
\[
q_i=c_0(t_i,t_i^2,t_i^3,0,\ldots,0)
\in \mathbb B_{1/4}\subset\R^N,
\]
where $c_0>0$ is small. Any four of these points are
affinely independent (i.e., not coplanar), by the Vandermonde identity
\[
\det\bigl(t_{i_a}^{\,b-1}\bigr)_{a,b=1}^4
=\prod_{1\le a<b\le4}(t_{i_b}-t_{i_a})\ne0.
\]
Thus segments $[q_i,q_j]$ and $[q_k,q_l]$ with disjoint index pairs
cannot intersect. Choose congruent balls $R_i$ centered at $q_i$,
with pairwise disjoint closures and radii small enough that
$\overline{R_i}$ meets $[q_k,q_l]$ only when $i\in\{k,l\}$. \emph{Here, we emphasize that the assumption $N\geq 3$ is necessary. Indeed, when $N=2$, one can have at most four points whose pairwise connecting segments are mutually disjoint except at their endpoints, and hence $m\leq 3$. Thus, only when $N\geq 3$ can $m$ be chosen as an arbitrary positive integer.}

 To each edge $e=\{i,j\}$, $1\le i<j\le m+1$,
we take a sufficiently thin round cylinder along $[q_i,q_j]$ and
remove the two closed endpoint balls to obtain a tube $T_e$.
Let $\Gamma_{ei}\subset\partial R_i$ and
$\Gamma_{ej}\subset\partial R_j$ be its spherical interfaces.
The cylinder radii can be chosen so that
\[
\overline{T_e}\cap\overline{T_f}=\varnothing\quad(e\ne f),
\qquad
\overline{\Gamma_{ei}}\cap\overline{\Gamma_{ej}}=\varnothing,
\qquad
|\Gamma_{ei}|=|\Gamma_{ej}|>0.
\]
Their locations are illustrated in the following figure for
$m+1=4$.

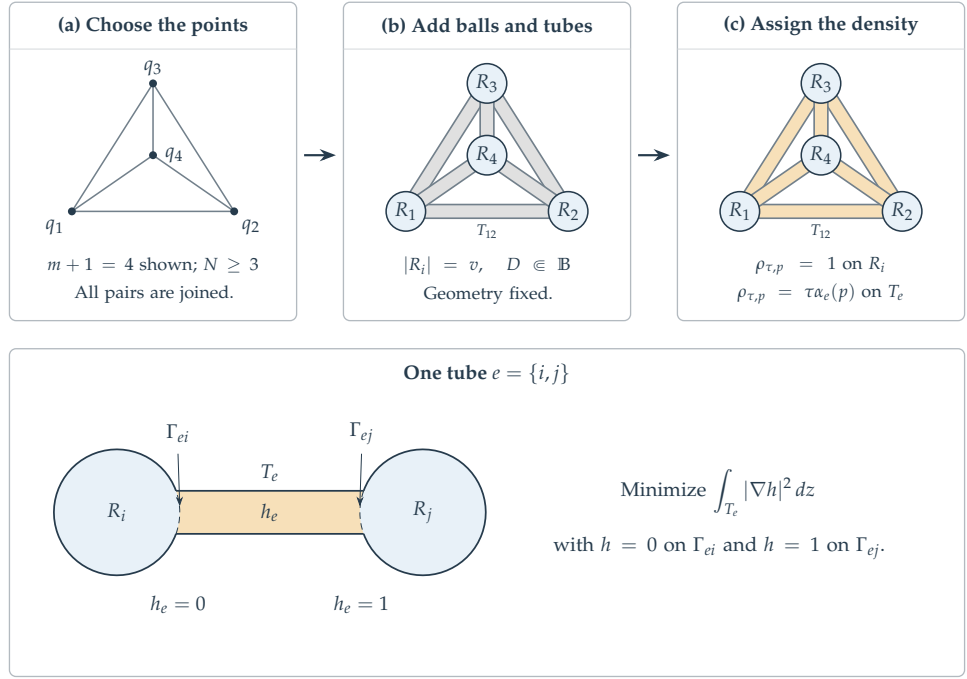
\begin{figure}[htbp]
\centering
\resizebox{\linewidth}{!}{%
\begin{tikzpicture}[
    font=\small,
    ink/.style={draw=sec31ink},
    panel/.style={
        draw=sec31ink!35,
        rounded corners=3pt,
        line width=0.6pt
    },
    heading/.style={
        font=\small\bfseries,
        text=sec31ink,
        align=center
    },
    room/.style={
        circle,
        draw=sec31ink,
        fill=sec31room,
        minimum size=18pt,
        inner sep=1pt,
        line width=0.7pt
    },
    flow/.style={
        -{Stealth[length=2.6mm]},
        draw=sec31ink,
        line width=1pt
    },
    every node/.style={text=sec31ink}
]


\foreach \x in {0,5.35,10.70}{
    \draw[panel] (\x,0) rectangle ++(4.60,5.10);
    \draw[sec31ink!20] (\x,4.35) -- ++(4.60,0);
}

\node[heading] at (2.30,4.72)
    {(a) Choose the points};

\node[heading] at (7.65,4.72)
    {(b) Add balls and tubes};

\node[heading] at (13.00,4.72)
    {(c) Assign the density};

\draw[flow] (4.72,2.65) -- (5.23,2.65);
\draw[flow] (10.07,2.65) -- (10.58,2.65);


\begin{scope}[shift={(2.30,2.65)}]

    \coordinate (P1) at (-1.30,-0.90);
    \coordinate (P2) at ( 1.30,-0.90);
    \coordinate (P3) at ( 0.00, 1.15);
    \coordinate (P4) at ( 0.00, 0.00);

    \foreach \i/\j in {1/2,1/3,1/4,2/3,2/4,3/4}{
        \draw[sec31ink!65,line width=0.7pt]
            (P\i) -- (P\j);
    }

    \foreach \i in {1,2,3,4}{
        \fill[sec31ink] (P\i) circle (1.8pt);
    }

    \node[below left]  at (P1) {$q_1$};
    \node[below right] at (P2) {$q_2$};
    \node[above]       at (P3) {$q_3$};
    \node[right=2pt]   at (P4) {$q_4$};

\end{scope}

\node[
    align=center,
    font=\footnotesize,
    text width=4.20cm
] at (2.30,0.65)
    {$m+1=4$ shown; $N\ge3$\\[3pt]
     All pairs are joined.};


\foreach \X/\tubeFill in {7.65/black!12,13.00/sec31tube}{

    \begin{scope}[shift={(\X,2.65)}]

        \coordinate (P1) at (-1.30,-0.90);
        \coordinate (P2) at ( 1.30,-0.90);
        \coordinate (P3) at ( 0.00, 1.15);
        \coordinate (P4) at ( 0.00, 0.00);

        \foreach \i/\j in {1/2,1/3,1/4,2/3,2/4,3/4}{
            \draw[sec31ink!65,line width=7pt]
                (P\i) -- (P\j);
            \draw[\tubeFill,line width=5.5pt]
                (P\i) -- (P\j);
        }

        \foreach \i in {1,2,3,4}{
            \node[room] at (P\i) {$R_{\i}$};
        }

        \node[font=\scriptsize] at (0,-1.22)
            {$T_{12}$};

    \end{scope}
}

\node[
    align=center,
    font=\footnotesize,
    text width=4.20cm
] at (7.65,0.65)
    {$|R_i|=v$, \quad $D\Subset \mathbb{B}$\\[3pt]
     Geometry fixed.};

\node[
    align=center,
    font=\footnotesize,
    text width=4.20cm
] at (13.00,0.65)
    {$\rho_{\tau,p}=1$ on $R_i$\\[3pt]
     $\rho_{\tau,p}=\tau\alpha_e(p)$ on $T_e$};


\draw[panel] (0,-0.45) rectangle (15.30,-5.70);

\node[heading] at (7.65,-0.85)
    {One tube $e=\{i,j\}$};

\begin{scope}[shift={(-0.10,0.50)},scale=1.40]

    \coordinate (Li) at (1.30,-2.55);
    \coordinate (Lj) at (4.80,-2.55);

    \coordinate (LU) at ($(Li)+(20:0.72)$);
    \coordinate (LD) at ($(Li)+(-20:0.72)$);
    \coordinate (RU) at ($(Lj)+(160:0.72)$);
    \coordinate (RD) at ($(Lj)+(200:0.72)$);

    \path[fill=sec31tube]
        (LU)--(RU)--(RD)--(LD)--cycle;

    \fill[sec31room] (Li) circle (0.72);
    \fill[sec31room] (Lj) circle (0.72);

    \draw[ink,line width=0.8pt]
        (LU) arc (20:340:0.72);

    \draw[ink,line width=0.8pt]
        (RU) arc (160:-160:0.72);

    \draw[ink,line width=0.8pt]
        (LU)--(RU)
        (LD)--(RD);

    \draw[ink,densely dashed]
        (LD) arc (-20:20:0.72);

    \draw[ink,densely dashed]
        (RU) arc (160:200:0.72);

    \node at (Li) {$R_i$};
    \node at (Lj) {$R_j$};

    \node at (3.05,-2.10) {$T_e$};
    \node at (3.05,-2.55) {$h_e$};

    \node at (2.00,-1.65) {$\Gamma_{ei}$};
    \node at (4.10,-1.65) {$\Gamma_{ej}$};

    \draw[ink,-{Stealth[length=1.5mm]}]
        (2.00,-1.88) -- (2.02,-2.48);

    \draw[ink,-{Stealth[length=1.5mm]}]
        (4.10,-1.88) -- (4.08,-2.48);

    \node at (2.00,-3.60) {$h_e=0$};
    \node at (4.10,-3.60) {$h_e=1$};

\end{scope}

\node[
    align=center,
    font=\small,
    text width=6.20cm
] at (11.35,-3.10)
    {Minimize
     $\displaystyle\int_{T_e}|\nabla h|^2\,dz$\\[8pt]
     with $h=0$ on $\Gamma_{ei}$
     and $h=1$ on $\Gamma_{ej}$.};

\end{tikzpicture}%
}
\label{fig:sec31-construction}
\end{figure}
The cylindrical and spherical boundary pieces meet transversely,
and no three boundary faces meet. Consequently,
\begin{equation}\label{eq:fixed-domain}
D:=\operatorname{int}\left(
\bigcup_{i=1}^{m+1}\overline{R_i}
\cup\bigcup_e\overline{T_e}
\right)\Subset\mathbb{B}
\end{equation}
is connected and Lipschitz, and
$O:=\mathbb{B}\setminus\overline D$ is also Lipschitz.
All these geometric choices are fixed; in particular, the tube
widths do not tend to zero (Recall Section \ref{sec:overview} that $\eta_{e,\tau}\to 0$ for the Neumann problem studied in \cite{CdV}). Write $v:=|R_i|>0$.

\textbf{The minimizing problem on each tube.}
For $e=\{i,j\}$, consider
\begin{equation}\label{Sec 3 def h}
\min\left\{
\int_{T_e}|\nabla h|^2\dd:
h\in H^1(T_e),\quad
h|_{\Gamma_{ei}}=0,\quad h|_{\Gamma_{ej}}=1
\right\}.
\end{equation}
No trace condition is imposed on the lateral boundary.
The admissible class is nonempty because the two closed interfaces
are disjoint. A standard contrary argument shows that differences of admissible functions satisfy
\[
	\|w\|_{H^1(T_e)}\le C\|\nabla w\|_{L^2(T_e)}
	\qquad\mathrm{for~any}\qquad
		w|_{\Gamma_{ei}\cup\Gamma_{ej}}=0.
\]
The direct method therefore yields a minimizer $h_e$, satisfying
\begin{equation}\label{Sec equ-q}
\int_{T_e}\nabla h_e\cdot\nabla w\dd=0
	\qquad\mathrm{for~any}\qquad
		w|_{\Gamma_{ei}\cup\Gamma_{ej}}=0.
\end{equation}
	If $h$ and $\widetilde h$ are two minimizers, then $h-\widetilde h$ also satisfies \eqref{Sec equ-q}. Since $(h-\widetilde h)|_{\Gamma_{ei}\cup\Gamma_{ej}}=0$, we may test with $h-\widetilde h$ to obtain
	$
	\int_{T_e}|\nabla(h-\widetilde h)|^2\,dz=0.
	$
	Their identical traces imply $h=\widetilde h$. Let $\tilde{h}_e=\min\{1,\max\{h_e,0\}\}\in H^1(T_e)$, then $\tilde{h}_e|_{\Gamma_{ei}}=0,\tilde{h}_e|_{\Gamma_{ej}}=1 $ and $\int_{T_e}|\nabla \tilde{h}_e|^2\,dz \leq \int_{T_e}|\nabla h_e|^2\,dz$, so uniqueness gives $0\le h_e\le1$. Define
\begin{equation}\label{eq:tube-energy}
\kappa_e:=\int_{T_e}|\nabla h_e|^2\dd>0.
\end{equation}
The strict positivity of $h_e$ in the interior of $T_e$ follows
from $\Delta h_e=0$ in $T_e$.

\textbf{The density and the local forms.} For $p\in K_r$ and $0<\tau\le1$, we assign the weight
$c_e(p):=c_{ij}(p)$ for each edge $e=\{i,j\}$ from Lemma~\ref{lem:graphweights}, and
set
\begin{equation}\label{eq:parameter-density}
\alpha_e(p):=\frac{v\,c_e(p)}{\kappa_e},
\qquad
\rho_{\tau,p}(z):=
\begin{cases}
1,&z\in R_i,\\
\tau\alpha_e(p),&z\in T_e.
\end{cases}
\end{equation}
This defines the density almost everywhere on $D$. Since
$1/2\le c_e(p)\le3/2$ and there are finitely many tubes,
\[
0<\min_e\frac{v}{2\kappa_e}
\le\alpha_e(p)
\le\max_e\frac{3v}{2\kappa_e}<\infty.
\]
For $f,g\in H^1(D)$, define
\[
a_{\tau,p}(f,g)
:=\int_D\rho_{\tau,p}\nabla f\cdot\nabla g\dd,
\qquad
b_{\tau,p}(f,g)
:=\int_D\rho_{\tau,p}fg\dd.
\]
These forms will be used to estimate the contributions from $D$
to the weighted Dirichlet problem on $\mathbb{B}$.

\subsection{Trial Functions and Uniform Estimates}
\label{sec:graphmatrix}

Define a linear map $J:\R^{m+1}\to H^1(D)$ by
\begin{equation}\label{eq:trialmap}
(Jx)|_{R_i}:=\frac{x_i}{\sqrt v},
\qquad
(Jx)|_{T_e}
:=\frac{(1-h_e)x_i+h_ex_j}{\sqrt v},
\quad e=\{i,j\}.
\end{equation}
The traces agree at both ends, so $Jx\in H^1(D)$.
Direct integration gives
\begin{align}
a_{\tau,p}(Jx,Jy)&=\sum_{e=\{i,j\}}\tau\alpha_e(p)\int_{T_e}\nabla Jx\cdot\nabla Jy\,dz\notag\\
&=\sum_{e=\{i,j\}}\frac{\tau\alpha_e(p)}v
(x_j-x_i)(y_j-y_i)\kappa_e\notag\\
&=\tau\sum_{i<j}c_{ij}(p)(x_i-x_j)(y_i-y_j)\notag\\
&=\tau\sum_{i<j}c_{ij}(p)x^TE_{ij}y
=\tau x^TL_{c(p)}y,
\label{eq:exactgraphenergy}
\end{align}
and
\begin{equation}\label{eq:graphmass}
b_{\tau,p}(Jx,Jy)=x\cdot y+\tau x^TM_py,
\qquad
\sup_{p\in K_r}\|M_p\|<\infty,
\end{equation}
where $M_p\ge0$ is defined by
\[
x^TM_py:=\sum_e\alpha_e(p)\int_{T_e}(Jx)(Jy)\dd.
\]
In particular, $J$ is injective and
\begin{equation}\label{eq:JxH}
  \|Jx\|_{H^1(D)}\le C|x|.  
\end{equation}

\begin{lemma}\label{lem:graphestimates}
There is $C>0$, independent of $\tau\in(0,1]$ and $p\in K_r$,
such that
\begin{equation}\label{eq:codimpoincare}
b_{\tau,p}(f,f)\le C a_{\tau,p}(f,f)
\qquad\mathrm{if}
\quad
f\in H^1(D),~\int_{R_i}f\dd=0\ \text{~for~ all~ }i.
\end{equation}
Moreover, for every $x\in\R^{m+1}$ and $g\in H^1(D)$,
\begin{equation}\label{eq:residual}
|a_{\tau,p}(Jx,g)|
\le C\tau|x|
\bigl(a_{\tau,p}(g,g)+b_{\tau,p}(g,g)\bigr)^{1/2}.
\end{equation}
\end{lemma}

\begin{proof}
	Let \(f\in H^1(D)\) satisfy \(\int_{R_i}f\dd=0\) for every
		\(i\). The room Poincaré and trace inequalities give
		\begin{equation}\label{eq:roomzero}
			\|f\|_{H^1(R_i)}^2+\|f\|_{L^2(\Gamma_{ei})}^2
			\le C\|\nabla f\|_{L^2(R_i)}^2.
		\end{equation}
On each fixed tube,
\begin{equation}\label{eq:tubepoincare}
\|f\|_{L^2(T_e)}^2
\le C\left(\|\nabla f\|_{L^2(T_e)}^2
+\|f\|_{L^2(\Gamma_{ei})}^2\right).
\end{equation}
Indeed, otherwise there would be a sequence with
$\|f_n\|_{L^2(T_e)}=1$ and
$\|\nabla f_n\|_{L^2(T_e)}+
\|f_n\|_{L^2(\Gamma_{ei})}\to0$.
Weak $H^1$ and strong $L^2$ convergence would give a limit $f$
with $\nabla f=0$, $f|_{\Gamma_{ei}}=0$, and
$\|f\|_{L^2(T_e)}=1$, a contradiction.

Multiplying \eqref{eq:tubepoincare} by $\tau\alpha_e(p)$ and
summing, we obtain from \eqref{eq:roomzero} and $\tau\le1$ that
\begin{align*}
b_{\tau,p}(f,f)
&=\sum_i\|f\|_{L^2(R_i)}^2
+\tau\sum_e\alpha_e(p)\|f\|_{L^2(T_e)}^2\\
&\le C\sum_i\|\nabla f\|_{L^2(R_i)}^2
+C\tau\sum_e\alpha_e(p)\|\nabla f\|_{L^2(T_e)}^2
\le C a_{\tau,p}(f,f).
\end{align*}
This proves \eqref{eq:codimpoincare}.

Let $\mathcal E_i:H^1(R_i)\to H^1(\R^N)$ be bounded extension
operators.  For \(e=\{i,j\}\), choose sufficiently small neighborhoods \(V_{ei}\) and \(V_{ej}\) of \(\Gamma_{ei}\) and \(\Gamma_{ej}\), respectively, and choose smooth cutoff functions
		\(\chi_{ei},\chi_{ej}\) such that \(\mathrm{supp}(\chi_{ei})\cap\mathrm{supp}(\chi_{ej})=\emptyset\) and
		\begin{align*}
			\chi_{ei}=1\quad\mathrm{in}\quad V_{ei},\quad\chi_{ei}=0\quad\mathrm{in}\quad V_{ej};\\
			\chi_{ej}=1\quad\mathrm{in}\quad V_{ej},\quad \chi_{ej}=0\quad\mathrm{in}\quad V_{ei}.
		\end{align*}
     For $g\in H^1(D)$, set
\[
E_eg:=\left.
\left(\chi_{ei}\mathcal E_i(g|_{R_i})
+\chi_{ej}\mathcal E_j(g|_{R_j})\right)\right|_{T_e}.
\]
The traces of $g$ agree across the interfaces (The trace and extension results used here hold on Lipschitz
domains; see e.g. \cite[Chapter~3]{McLean}). Hence
\begin{equation}\label{eq:endextension}
(g-E_eg)|_{\Gamma_{ei}\cup\Gamma_{ej}}=0,
\qquad
\|E_eg\|_{H^1(T_e)}
\le C\left(\|g\|_{H^1(R_i)}+\|g\|_{H^1(R_j)}\right).
\end{equation}
By \eqref{Sec equ-q} and \eqref{eq:trialmap},
\begin{align*}
a_{\tau,p}(Jx,g)
&=\sum_{e=\{i,j\}}
\frac{\tau\alpha_e(p)(x_j-x_i)}{\sqrt v}
\int_{T_e}\nabla h_e\cdot\nabla g\dd\\
&=\sum_{e=\{i,j\}}
\frac{\tau\alpha_e(p)(x_j-x_i)}{\sqrt v}
\int_{T_e}\nabla h_e\cdot\nabla E_eg\dd.
\end{align*}
Since $\rho_{\tau,p}=1$ on each room implies 
$$a_{\tau,p}(g,g)+b_{\tau,p}(g,g)\geq \sum_i\int_{R_i}|\nabla g|^2+g^2dz,$$
it follows that
\begin{align*}
|a_{\tau,p}(Jx,g)|&\le C\tau|x|\sum_{e=\{i,j\}}\|E_eg\|_{H^1(T_e)}
\le C\tau|x|\sum_i\|g\|_{H^1(R_i)}\\
&\le C\tau|x|
\bigl(a_{\tau,p}(g,g)+b_{\tau,p}(g,g)\bigr)^{1/2}.
\end{align*}
The proof is complete.
\end{proof}

\section{From Weighted Dirichlet Problems to Schr\"odinger operator}
\label{sec:potential}

Retain $D$, $H$, $K_r$, $\rho_{\tau,p}$, and $J$ from
Section~\ref{sec:graph}. We assign weight $\tau^2$ to
$\mathbb{B}\setminus\overline D$ and impose Dirichlet conditions
on $\partial\mathbb{B}$. Thus the room, tube, and exterior
contributions occur at the scales $1$, $\tau$, and $\tau^2$,
respectively. Using the local estimates from
Lemma~\ref{lem:graphestimates}, we prove the spectral approximation
directly on $\mathbb{B}$. We then smooth the density and convert
the weighted problem into a Schr\"odinger operator problem.

\begin{proposition}\label{prop:potential}
There exist small $\tau>0$ and a family
$V_p\in C^\infty(\overline{\mathbb{B}};\R)$, continuous in
$p\in K_r$ in every $C^k$ norm, such that
$L_p:=-\Delta_{\mathbb{B}}^D+V_p$ is positive and its eigenvalues satisfy
\begin{align}\label{eq:potentialgaps}
\begin{gathered}
0<\xi_1(p)<\frac{\tau(m+1)}4,\\
\frac{\tau(m+1)}2<\xi_2(p)\le\cdots\le\xi_{m+1}(p)
<\frac{3\tau(m+1)}2,\\
\xi_{m+2}(p)>2\tau(m+1).
\end{gathered}
\end{align}
The subspace for eigenvalues $2,\ldots,m+1$ has a continuous
orthonormal frame, whose operator matrix $\mathcal B_p$ on $H$ (i.e., $\mathcal B_p\in\Sym(H)$ and its eigenvalues are precisely \(\xi_2(p),\ldots,\xi_{m+1}(p)\)) satisfies
\begin{equation}\label{eq:potentialmatrix}
\sup_{p\in K_r}
\|\mathcal B_p-\tau((m+1)I_H+p)\|<\frac{\tau r}{4}.
\end{equation}
\end{proposition}

\begin{remark}\label{rmk3-2} Let us explain the operator matrix $\mathcal B_p$ and the estimate \eqref{eq:potentialmatrix} in details.
Let \(E_p\subset L^2(\mathbb B;\mathbb R)\) be the \(m\)-dimensional
spectral subspace of \(L_p\) corresponding to
\(\xi_2(p),\ldots,\xi_{m+1}(p)\), counted with multiplicity.
A continuous orthonormal frame of \(E_p\) means a choice of
orthonormal bases \(u_{1,p},\ldots,u_{m,p}\) such that each map
\(p\mapsto u_{j,p}\) is continuous in \(L^2(\mathbb B)\) for any $j=1,\cdots,m$.
These functions need not be eigenfunctions.
Using the fixed orthonormal basis \(h_1,\ldots,h_m\) of \(H\), we
define the linear isometric isomorphism
\[
J_p:H\longrightarrow E_p,
\qquad J_ph_j:=u_{j,p},\quad j=1,\ldots,m.
\] 
Let $J_p^*: E_p\to H$ be the adjoint of $J_p$, namely
$$\langle J_p x, f\rangle_{L^2(\mathbb B)}=\langle x, J_p^*f\rangle_{H},\quad x\in H, \quad f\in E_p.$$
Then $J_p^*J_p=I_H$, and
the operator \(\mathcal B_p\in\Sym(H)\) is defined by
\[
\mathcal B_p:=J_p^*(L_p|_{E_p})J_p.
\]
Its matrix entries in the basis \(h_1,\ldots,h_m\) are
\[
\begin{aligned}
\langle\mathcal B_ph_k,h_j\rangle_H
=\langle L_pu_{k,p},u_{j,p}\rangle_{L^2(\mathbb B)}
=\int_{\mathbb B}
\bigl(\nabla u_{k,p}\cdot\nabla u_{j,p}
      +V_pu_{k,p}u_{j,p}\bigr)\,dz.
\end{aligned}
\]
Thus \(\mathcal B_p\) is represented by a real symmetric matrix
whose eigenvalues are precisely
\(\xi_2(p),\ldots,\xi_{m+1}(p)\).
The estimate \eqref{eq:potentialmatrix} states that, under these
chosen identifications of \(E_p\) with the fixed space \(H\),
the restricted operators \(L_p|_{E_p}\) uniformly approximate
\(\tau((m+1)I_H+p)\in \Sym (H)\).
\end{remark}

The rest of this section is devoted to the proof of Proposition \ref{prop:potential}.

\subsection{Spectral Reduction for the Weighted Dirichlet Problem}

Denote $O=\mathbb{B}\setminus\overline D$ and define
\begin{equation}\label{eq:outsideweight}
\hat{\rho}_{\tau,p}=
\begin{cases}
\rho_{\tau,p},&\text{in }D,\\
\tau^2,&\text{in }O.
\end{cases}
\end{equation}
For $0<\tau\le1$, we have \begin{align}\label{hat rho}
    c\tau^2\le\hat{\rho}_{\tau,p}\le C,
\end{align}
uniformly in $p$. On $H_0^1(\mathbb{B})$, set
\[
\widehat a_{\tau,p}(f,g)
:=\int_{\mathbb{B}}\hat{\rho}_{\tau,p}\nabla f\cdot\nabla g\dd,
\qquad
\widehat b_{\tau,p}(f,g)
:=\int_{\mathbb{B}}\hat{\rho}_{\tau,p}fg\dd.
\]
Consider the weighted Dirichlet eigenvalue problem
\[
\widehat a_{\tau,p}(f,g)
=\nu\,\widehat b_{\tau,p}(f,g),
\qquad g\in H_0^1(\mathbb{B}),
\]
and denote its eigenvalues, counted with multiplicity, by
$\nu_j(\tau,p)$.

\begin{lemma}\label{lem:weighted-dirichlet}
There exist $C,\widehat\gamma>0$, independent of $\tau$ and $p$,
such that, for all sufficiently small $\tau>0$ and every $p\in K_r$,
\begin{equation}\label{eq:weighted-dirichlet-bounds}
0<\nu_1(\tau,p)\le C\tau^2,
\qquad
\nu_{m+1}(\tau,p)\le C\tau,
\qquad
\nu_{m+2}(\tau,p)\ge\widehat\gamma.
\end{equation}
The spectral subspace corresponding to eigenvalues $2,\ldots,m+1$
has a $\widehat b_{\tau,p}$-orthonormal frame
$u_1=u_1(\tau,p),\ldots,u_m=u_m(\tau,p)\in H_0^1(\mathbb{B})$, whose operator matrix
$\widehat B_{\tau,p}$ on $H$ satisfies
\begin{equation}\label{eq:weighted-dirichlet-matrix}
\sup_{p\in K_r}
\|\widehat B_{\tau,p}-\tau((m+1)I_H+p)\|
\le C\tau^2.
\end{equation}
\end{lemma}

\begin{proof}
Write $\widehat a=\widehat a_{\tau,p}$ and
$\widehat b=\widehat b_{\tau,p}$.
All error bounds below are uniform in $p\in K_r$.

First, we aim to prove the eigenvalue estimates in \eqref{eq:weighted-dirichlet-bounds}.
Choose a fixed bounded extension
$\mathcal E:H^1(D)\to H_0^1(\mathbb{B})$ with
$(\mathcal Ef)|_D=f$, and set $\widehat J=\mathcal EJ$, where $J: \mathbb R^{m+1}\to H^1(D)$ is the linear injective map defined in \eqref{eq:trialmap}.
Then \eqref{eq:JxH} implies \begin{equation}\label{eq:wJx}\|\widehat Jx\|_{H^1(\mathbb{B})}\le C|x|.\end{equation}
By \eqref{eq:exactgraphenergy}--\eqref{eq:graphmass} and \eqref{eq:outsideweight},
\begin{align}
\widehat a(\widehat Jx,\widehat Jy)
&=\tau x^TL_{c(p)}y+O(\tau^2)|x||y|,\notag\\
\widehat b(\widehat Jx,\widehat Jy)
&=x\cdot y+\tau x^TM_py+O(\tau^2)|x||y|.
\label{eq:weighted-dirichlet-trials}
\end{align}
Since $\widehat b(\widehat Jx,\widehat Jx)\ge b(Jx, Jx)\ge|x|^2$, the
min--max principle
\begin{align*}
    \nu_k
=
\min_{\substack{ \dim V=k}}
\max_{f\in V\setminus\{0\}}
\frac{\widehat a(f,f)}{\widehat b(f,f)}
=
\max_{\substack{ \dim V=k-1}}
\inf_{\substack{f\in H^1_0(\mathbb{B})\setminus\{0\}\\ f\perp_{\widehat b} V}}
\frac{\widehat a(f,f)}{\widehat b(f,f)},
\end{align*}
gives $$\nu_{m+1}\le \max_{x\in\mathbb R^{m+1}\setminus\{0\}}\frac{\widehat a(\widehat Jx,\widehat Jx)}{\widehat b(\widehat Jx,\widehat Jx)}\le C\tau.$$
Moreover, $J\mathbf{1}\equiv 1/\sqrt{v}$ is constant on $D$, so testing with
$\widehat J\mathbf{1}$ and $L_{c(p)}\mathbf{1}=0 $ gives $$\nu_1\le \frac{\widehat a(\widehat J\mathbf{1},\widehat J\mathbf{1})}{\widehat b(\widehat J\mathbf{1},\widehat J\mathbf{1})} \le C\tau^2.$$
The  Poincar\'e inequality and \eqref{hat rho} imply $\nu_1>0$.

If $f\in H_0^1(\mathbb{B})$ satisfies $\int_{R_i}f\dd=0$ for
every $i$, then \eqref{eq:codimpoincare} gives
\[
b_{\tau,p}(f|_D,f|_D)\le C a_{\tau,p}(f|_D,f|_D).
\]
Furthermore, again using Poincar\'e inequality and \eqref{hat rho}, we have
\[
\tau^2\int_O f^2\dd
\le C\tau^2\int_{\mathbb{B}}|\nabla f|^2\dd
\le C\widehat a(f,f).
\]
Hence$$\widehat b(f,f)= b_{\tau,p}(f|_D,f|_D)+\tau^2\int_O f^2\dd\le C\widehat a(f,f).$$
For any $(m+2)$-dimensional subspace
$V\subset H_0^1(\mathbb{B})$,
\[
\dim\ker\left(f\in V\longmapsto
\left(\int_{R_i}f\dd\right)_{i=1}^{m+1}\right)\ge1.
\]
The min--max principle therefore yields
\begin{align}\label{Sec 3 equ-a1}
    \nu_{m+2}
=\inf_{\dim V=m+2}\ \sup_{0\ne f\in V}
\frac{\widehat a(f,f)}{\widehat b(f,f)}
\ge C^{-1}=:\widehat\gamma.
\end{align}

Next, we estimate the error terms and define the corresponding operator matrix.
By \eqref{eq:residual} and \eqref{eq:wJx}, for $x\in\R^{m+1}$ and
$g\in H_0^1(\mathbb{B})$,
\begin{align}
|\widehat a(\widehat Jx,g)|
&\le |a_{\tau,p}(Jx,g|_D)|
+C\tau^2|x|\|\nabla g\|_{L^2(O)}\notag\\
&\le C\tau|x|
\bigl(\widehat a(g,g)+\widehat b(g,g)\bigr)^{1/2}.
\label{eq:weighted-dirichlet-improved}
\end{align}
Let $\varphi_1^{\tau,p},
	\ldots,
	\varphi_{m+1}^{\tau,p}$ be the $\widehat b $-orthonormal eigenfunctions corresponding to eigenvalues $\nu_1,\ldots,\nu_{m+1}$ respectively, namely, $\widehat b(\varphi_k^{\tau,p},\varphi_l^{\tau,p})=\delta_{kl} $.
Define
	\[
	E_{\tau,p}
	:=
	\operatorname{span}
	\left\{
	\varphi_1^{\tau,p},
	\ldots,
	\varphi_{m+1}^{\tau,p}
	\right\}
	\]
	and let
	\begin{equation}\label{eq:widehatP}
	\widehat P:
	L^2(\mathbb{B})
	\longrightarrow
	E_{\tau,p}
	\end{equation}
	be the $\hat{b}$-orthogonal projection. Explicitly,
	\[
	\widehat Pf
	:=
	\sum_{k=1}^{m+1}
	\left(
	\int_{\mathbb{B}}
	\hat{\rho}_{\tau,p}(z)
	f(z)\varphi_k^{\tau,p}(z)\,dz
	\right)
	\varphi_k^{\tau,p}.
	\]
Set $z_x:=(I-\widehat P)\widehat Jx$.
Spectral orthogonality, the gap and \eqref{Sec 3 equ-a1} give
\begin{align}\label{Sec 3 Lem 3 euq-0}
\widehat b(z_x,z_x)
\le\widehat\gamma^{-1}\widehat a(z_x,z_x), \qquad \widehat a(\widehat P\widehat Jx,z_x)=\widehat b(\widehat P\widehat Jx,z_x)=0
\end{align}
		Thus, it follows from \eqref{eq:weighted-dirichlet-improved} and \eqref{Sec 3 Lem 3 euq-0} that
		\begin{align*}
			\widehat a(z_x,z_x)
&=\widehat a(\widehat Jx,z_x)
\le C\tau|x|\widehat a(z_x,z_x)^{1/2}.
		\end{align*}
Consequently,
\[
\widehat a(z_x,z_x)+\widehat b(z_x,z_x)\le C\tau^2|x|^2.
\]
This, together with spectral orthogonality, yields
$$\widehat b(\widehat Je_i,
\widehat Je_j)=\widehat b(\widehat P\widehat Je_i,
\widehat P\widehat Je_j)+O(\tau^2),\quad \widehat a(\widehat Je_i,
\widehat Je_j)=\widehat a(\widehat P\widehat Je_i,
\widehat P\widehat Je_j)+O(\tau^2).$$
From here and
\eqref{eq:weighted-dirichlet-trials}, we introduce the Gram matrix $\widehat G$ and energy matrix $\widehat K$ by
\begin{align*}
\widehat G
&:=\bigl(\widehat b(\widehat P\widehat Je_i,
\widehat P\widehat Je_j)\bigr)_{i,j=1}^{m+1}
=I_{m+1}+\tau M_p+O(\tau^2),\\
\widehat K
&:=\bigl(\widehat a(\widehat P\widehat Je_i,
\widehat P\widehat Je_j)\bigr)_{i,j=1}^{m+1}
=\tau L_{c(p)}+O(\tau^2).
\end{align*}
For small $\tau$, $\widehat G$ is positive definite and
$\widehat G^{-1/2}=I_{m+1}+O(\tau)$.
Define
\[
\widehat U:=\widehat P\widehat J\widehat G^{-1/2}
:\R^{m+1}\longrightarrow E_{\tau,p}\subset H_0^1(\mathbb{B}).
\]
Since $\widehat b(\widehat P\widehat Jx,\widehat P\widehat Jy)=\sum_{i,j}x_iy_j\widehat b( \widehat P\widehat Je_i,\widehat P\widehat Je_j)=x^{T}\widehat Gy$, then 
\begin{align}\label{Sec 3 isometry b}
    \widehat b(\widehat Ux,\widehat Uy)=\widehat b(\widehat P\widehat J(\widehat G^{-1/2}x),\widehat P\widehat J(\widehat G^{-1/2}y))=x\cdot y.
\end{align}
Thus, $\{\widehat U e_1,\cdots,\widehat U e_{m+1}\}$ is an orthonormal frame of the spectral space $E_{\tau,p}$ for $\nu_1,\cdots,\nu_{m+1}$ in the weighted space. 
Its operator matrix $\widehat A:=\bigl(\widehat a(\widehat Ue_i,
\widehat U e_j)\bigr)_{i,j=1}^{m+1} $satisfies
\begin{equation}\label{eq:weighted-dirichlet-low-matrix}
\widehat A
:=\widehat G^{-1/2}\widehat K\widehat G^{-1/2}
=\tau L_{c(p)}+O(\tau^2).
\end{equation}
Using \eqref{Sec 3 isometry b}, the eigenvalues of $\widehat A$ are exactly
$\nu_1,\ldots,\nu_{m+1}$. To verify this claim, we write
\begin{equation}\label{eq:qk-defn}
\varphi_k^{\tau,p}
=\sum_{j=1}^{m+1}q_j^{(k)}\widehat Ue_j=\widehat Uq^{(k)},
\qquad
q^{(k)}:=(q_1^{(k)},\ldots,q_{m+1}^{(k)})^T.
\end{equation}
Testing the eigenvalue equation for \(\varphi_k^{\tau,p}\)
with \(\widehat Ue_i\), and using the symmetry of \(\widehat a\), gives
\[
\sum_{j=1}^{m+1}\widehat A_{ij}q_j^{(k)}
=\widehat a(\varphi_k^{\tau,p},\widehat Ue_i)
=\nu_k\widehat b(\varphi_k^{\tau,p},\widehat Ue_i)
=\nu_kq_i^{(k)}.
\]
Hence \(\widehat A q^{(k)}=\nu_kq^{(k)}\).
Moreover,
\[
q^{(k)}\cdot q^{(l)}
=\widehat b(\varphi_k^{\tau,p},\varphi_l^{\tau,p})
=\delta_{kl}.
\]
Thus the matrix \(Q=(q^{(1)},\ldots,q^{(m+1)})\) is orthogonal
and satisfies
\[
Q^{-1}\widehat A Q
=\operatorname{diag}(\nu_1,\ldots,\nu_{m+1}).
\]
Therefore the eigenvalues of \(\widehat A\) are exactly \(\nu_1,\ldots,\nu_{m+1}\).

Finally, we remove the first eigenvalue and complete the proof.
Since $\mathbb R^{m+1}=\mathbb{R}\mathbf{1}\oplus H$, $L_{c(p)}\mathbf{1}=0$ and
$L_{c(p)}|_H=(m+1)I_H+p$ with $\|p\|\leq r$, the min--max principle,
\eqref{eq:weighted-dirichlet-low-matrix} and $\nu_1\leq C\tau^2$ imply
\begin{equation}\label{eq:nu2nu1}
\nu_2\ge\tau(m+1-r)-C\tau^2>\nu_1
\end{equation}
for sufficiently small $\tau$.

Recalling \eqref{eq:qk-defn}, we write $v=q^{(1)}$ for simplicity. Then $|v|=1$ and $\widehat Av=\nu_1v$. Define orthogonal projection maps $\Pi_H:\mathbb R^{m+1}\to H=\mathbf{1}^{\perp}$ and $ \widehat\Pi:\mathbb{R}^{m+1}
\to v^{\perp}=\operatorname{Span}\{q^{(2)},\ldots,q^{(m+1)}\}$ by
\[
\Pi_H:=I_{m+1}-\frac{\mathbf1\mathbf1^T}{m+1},
\qquad
\widehat\Pi:=I_{m+1}-vv^T.
\]
Projecting the eigenvalue equation $\widehat Av=\nu_1v$ onto $H$ gives
\[\Pi_H(\tau L_{c(p)}v-\widehat Av)=(\tau(m+1)-\nu_1)\Pi_Hv+\tau p\Pi_Hv.\]
Using $\langle p\Pi_Hv,\Pi_Hv\rangle\geq -r|\Pi_Hv|^2$ and \eqref{eq:weighted-dirichlet-low-matrix}, we get
\[
\bigl(\tau(m+1-r)-\nu_1\bigr)|\Pi_Hv|^2\leq\langle\Pi_H(\tau L_{c(p)}v-\widehat Av),\Pi_Hv\rangle
\le C\tau^2|\Pi_Hv|.
\]
Together with $\tau(m+1-r)-\nu_1\ge c\tau$, we obtain $|\Pi_Hv|\leq C\tau$. Write $v=\alpha\frac{\mathbf{1}}{\sqrt{m+1}}+\Pi_Hv$, then $\alpha^2=1-|\Pi_Hv|^2$. By replacing $v$ with $-v$ if necessary, we may assume $\alpha=\sqrt{1-|\Pi_Hv|^2}>0$, so $1-\alpha\leq |\Pi_Hv|^2\leq C\tau^2$, which implies
\[\left|v-\frac{\mathbf{1}}{\sqrt{m+1}}\right|\leq 1-\alpha+|\Pi_Hv|\leq C\tau.\]
Consequently,
\[
\norm{\widehat\Pi-\Pi_H}=\left\|vv^T-\frac{\mathbf1\mathbf1^T}{m+1}\right\|\leq 2\left|v-\frac{\mathbf{1}}{\sqrt{m+1}}\right|
\le C\tau.
\]

Since the desired inequality \eqref{eq:weighted-dirichlet-matrix} is formulated on the space $H$, we seek an orthonormal transformation $\widehat F:H\to v^{\perp}$ satisfying
\begin{align}\label{Sec 3 F def}
    \widehat F^*\widehat F=I_H,
\qquad
\sup_{x\in H,\,|x|=1}|\widehat Fx-x|\le C\tau.
\end{align}
Indeed, $|\Pi_Hv|\leq C\tau$ implies that the linear map $\widehat M: H\to H$ defined by
\[
\widehat M:=\Pi_H\widehat\Pi|_H
=I_H-(\Pi_Hv)(\Pi_Hv)^T=I_H+O(\tau^2)\]
is positive definite and hence $\widehat F: H\to v^{\perp}$ by
\[\widehat Fx:=\widehat\Pi\widehat M^{-1/2}x
\]
is well-defined. Denote $T:=\widehat\Pi|_H: H\to v^{\perp}$. Since $\widehat\Pi$ is an orthogonal projection, we have $T^*T=\widehat M$. From here and $\widehat F=T\widehat M^{-1/2}$, we get $\widehat F^*\widehat F=I_H$. Hence, $\widehat F $ is orthonormal, and then 
$\widehat B_{\tau,p}:=\widehat F^*\widehat A\widehat F$
has eigenvalues $\nu_2,\ldots,\nu_{m+1}$. 
For any $x\in H$, 
\begin{align*}
    |\widehat Fx-x|&=|\widehat\Pi(\widehat M^{-1/2}-I_H)x+(\widehat\Pi-\Pi_H)x|\\
    &\leq (\|\widehat M^{-1/2}-I_H\|+\|\widehat\Pi-\Pi_H\|)|x|\leq C\tau|x|.
\end{align*}
This proves \eqref{Sec 3 F def}.

For any $x,y\in H$, we have
		\begin{align*}
			&\langle \widehat B_{\tau,p}x,y\rangle
			-\tau\langle((m+1)I_H+p)x,y\rangle\\
			&=
			\langle \widehat A(\widehat Fx-x),
			\widehat Fy\rangle
			+
			\langle \widehat Ax,\widehat Fy-y\rangle
			+
			\langle(\widehat A-\tau L_{c(p)})x,y\rangle.
		\end{align*}
        	Since $\|\widehat A\|\le C\tau$ and
		$\|\widehat A-\tau L_{c(p)}\|\le C\tau^2$,
		we conclude that
		\begin{equation}\label{eq:upperblockrate}
			\|\widehat B_{\tau,p}-\tau((m+1)I_H+p)\|\le C\tau^2.
		\end{equation}
which proves \eqref{eq:weighted-dirichlet-matrix}.

It remains to give the $\widehat b_{\tau,p}$-orthonormal frame
$u_1,\ldots,u_m\in H_0^1(\mathbb{B})$ such that its operator matrix is $\widehat B_{\tau,p}$.
Recall the fixed
orthonormal basis $h_1,\ldots,h_m$ of $H=\mathbf 1^\perp$ defined in \eqref{eq:hk-def}, we have $\langle h_k,h_l\rangle_H=\delta_{kl}$.
Define
\begin{align}\label{Def u_k}
u_k:=\widehat U\widehat Fh_k\in E_{\tau,p}\subset H_0^1(\mathbb B),
\qquad k=1,\ldots,m.
\end{align}
Then it follows from $\widehat F: H\to v^{\perp}$ and \eqref{eq:qk-defn} that $\{u_k\}_{k=1}^m$ spans the weighted spectral subspace $\operatorname{span}
	\{
	\varphi_2^{\tau,p},
	\ldots,
	\varphi_{m+1}^{\tau,p}
	\}$ corresponding to $\nu_2,\ldots,\nu_{m+1}$, and by \eqref{Sec 3 isometry b},
\[
\widehat b(u_k,u_l)
=\langle\widehat Fh_k,\widehat Fh_l\rangle_{\R^{m+1}}
=\langle h_k,h_l\rangle_H=\delta_{kl}.
\]
Thus $\{u_k\}_{k=1}^m$ is a $\widehat b_{\tau,p}$-orthonormal basis
of $\operatorname{span}
	\{
	\varphi_2^{\tau,p},
	\ldots,
	\varphi_{m+1}^{\tau,p}
	\}$.
By the definition of $\widehat A$, the matrix entries in this basis are
\[
\begin{aligned}
\langle\widehat B_{\tau,p}h_k,h_l\rangle_H=\langle\widehat A\widehat Fh_k, \widehat Fh_l\rangle=\widehat a(\widehat U\widehat Fh_k, \widehat U\widehat Fh_l)=\widehat a(u_k,u_l).
\end{aligned}
\]
Hence $\widehat B_{\tau,p}$ is the operator matrix of $u_1,\ldots, u_m$
with respect to the fixed basis
$h_1,\ldots,h_m$ of $H$. The proof is complete.
\end{proof}

\begin{lemma}\label{lem:framecontinuity}
For each sufficiently small fixed $\tau>0$, the matrices
$\widehat G$, $\widehat K$, $\widehat A$ and
$\widehat B_{\tau,p}$ depend continuously on $p\in K_r$.
Moreover, for each $k=1,\ldots,m$, the map $p\mapsto u_k=u_k(\tau,p)$
is continuous in $H_0^1(\mathbb{B})$.
\end{lemma}

\begin{proof}
Fix $\tau>0$ and let $p_n\to p$.
Write $\hat{\rho}_n=\hat{\rho}_{\tau,p_n}$ and
$\hat{\rho}=\hat{\rho}_{\tau,p}$. Since
\[
\|\hat{\rho}_n-\hat{\rho}\|_{L^{\infty}(\mathbb{B})}\to0,
\qquad 0<c_\tau\le\hat{\rho}_n,\hat{\rho}\le C,
\]
the corresponding forms satisfy
\[
\begin{aligned}
(1-\varepsilon_n)\widehat a_{\tau,p}
&\le\widehat a_{\tau,p_n}
\le(1+\varepsilon_n)\widehat a_{\tau,p},\\
(1-\varepsilon_n)\widehat b_{\tau,p}
&\le\widehat b_{\tau,p_n}
\le(1+\varepsilon_n)\widehat b_{\tau,p},
\end{aligned}
\qquad\varepsilon_n\to0.
\]
The min--max principle implies
\begin{equation}\label{eq:eigenvalue-continuity}
\nu_k(\tau,p_n)\to\nu_k(\tau,p),\qquad k\ge1.
\end{equation}

For $1\le k\le m+1$, choose
$\widehat b_{\tau,p_n}$-orthonormal eigenfunctions $f_{k,n}$ (Here, we note that $f_{k,n}$ need not be chosen as $u_k(\tau,p_n)$).
By \eqref{eq:weighted-dirichlet-bounds},
$\|f_{k,n}\|_{H^1(\mathbb{B})}\le C_\tau$.
After passing to a subsequence,
\[
f_{k,n}\rightharpoonup f_k\quad\text{in }H_0^1(\mathbb{B}),
\qquad
f_{k,n}\to f_k\quad\text{in }L^2(\mathbb{B}).
\]
Passing to the limit in the eigenvalue equations gives
\[
\widehat a_{\tau,p}(f_k,g)
=\nu_k(\tau,p)\widehat b_{\tau,p}(f_k,g),
\qquad
\widehat b_{\tau,p}(f_k,f_l)=\delta_{kl}.
\]
Hence $f_1,\ldots,f_{m+1}$ form an orthonormal basis of
$E_{\tau,p}$. Moreover,
\begin{align*}
\widehat a_{\tau,p}(f_{k,n},f_{k,n})
+\widehat b_{\tau,p}(f_{k,n},f_{k,n})
=\nu_k(\tau,p_n)+1+o(1)
\longrightarrow\nu_k(\tau,p)+1.
\end{align*}
Together with weak convergence, this yields
$f_{k,n}\to f_k$ strongly in $H_0^1(\mathbb{B})$.
Recalling the $\hat{b}$-orthogonal projection $\widehat P_{\tau, p_n}$ in \eqref{eq:widehatP}, it follows that
\begin{align*}
\widehat P_{\tau,p_n}\widehat Jx
=\sum_{k=1}^{m+1}
\widehat b_{\tau,p_n}(\widehat Jx,f_{k,n})f_{k,n}
\rightarrow
\sum_{k=1}^{m+1}
\widehat b_{\tau,p}(\widehat Jx,f_k)f_k
=\widehat P_{\tau,p}\widehat Jx
\end{align*}
in $H_0^1(\mathbb{B})$.
The limit is independent of the chosen eigenbasis. The subsequence
argument and finite dimensionality therefore give
\begin{equation}\label{eq:projection-continuity}
\widehat P_{\tau,p_n}\widehat J
\longrightarrow\widehat P_{\tau,p}\widehat J
\quad\text{in }\mathcal L(\R^{m+1},H_0^1(\mathbb{B})).
\end{equation}
By their definitions, $\widehat G$ and $\widehat K$ are continuous.
Since $\widehat G$ is positive definite, $\widehat G^{-1/2}$
is continuous in $p$. Hence the maps
\[
p\longmapsto\widehat U
=\widehat P\widehat J\widehat G^{-1/2},
\qquad
p\longmapsto\widehat A
=\widehat G^{-1/2}\widehat K\widehat G^{-1/2}
\]
are continuous, where continuity of the first map is in
$\mathcal L(\R^{m+1},H_0^1(\mathbb{B}))$.

From \eqref{eq:nu2nu1}, the first eigenvalue of $\widehat A$ is simple and separated from
the others. If $p_n\to p$ and $v_n$ are normalized first
eigenvectors, every convergent subsequence has a limit in the
first eigenspace of $\widehat A_{\tau,p}$.
Since $vv^T$ is independent of the sign of $v$, it follows that
$\widehat\Pi_{\tau,p_n}\to\widehat\Pi_{\tau,p}$.
The formulas for $\widehat M$, $\widehat F$ and
$\widehat B_{\tau,p}$ now give their continuity, and hence that of
$\widehat U\widehat F$ in $\mathcal L(H,H_0^1(\mathbb{B}))$.
Thus $u_k(\tau,p)=\widehat U\widehat Fh_k$ is continuous in
$H_0^1(\mathbb{B})$ for every $k=1,\ldots,m$.
\end{proof}

\subsection{Smoothing and  Transform to Schr\"odinger operator}

Fix $\tau>0$ sufficiently small for
Lemma~\ref{lem:weighted-dirichlet}. Extend $\widehat{\rho}_{\tau,p}$ by $\tau^2$
outside $\mathbb{B}$, and set
$\widehat{\rho}_{\tau,p,\delta}:=\psi_\delta*\widehat{\rho}_{\tau,p}\in C^\infty$, where $\psi\in C_c^\infty(\mathbb R^{N})$ is a
nonnegative mollifier of integral one and $\psi_\delta(z):=\delta^{-N}\psi(z/\delta)$ for $\delta>0$. For small $\delta$,
\[
0<c_\tau\le\widehat{\rho}_{\tau,p,\delta}\le C,
\qquad
\widehat{\rho}_{\tau,p,\delta}=\tau^2\quad\text{near }\partial\mathbb{B},
\]
uniformly in $p,\delta$. For fixed $\delta$, the density is continuous
in $p$ in every $C^k$ norm, since
\[
\|\widehat{\rho}_{\tau,p,\delta}-\widehat{\rho}_{\tau,p',\delta}\|_{C^k(\overline{\mathbb{B}})}
\le C_{\delta,k}
\|\widehat{\rho}_{\tau,p}-\widehat{\rho}_{\tau,p'}\|_{L^\infty(\R^N)}.
\]
Let $\widehat a_{\tau,p,\delta}$ and $\widehat b_{\tau,p,\delta}$ be the
energy and mass forms with this density, and let
$\nu_j(\tau,p,\delta)$ be the corresponding weighted Dirichlet
eigenvalues.

For the fixed orthonormal basis $h_1,\ldots,h_m$ of $H$, recall \eqref{Def u_k} that
$u_k(\tau,p)=\widehat U\widehat Fh_k$, $1\le k\le m$.
By Lemma~\ref{lem:framecontinuity}, these functions are continuous
in $H_0^1(\mathbb{B})$, and
\[
\widehat b_{\tau,p}(u_k,u_l)=\delta_{kl},
\qquad
\widehat B_{\tau,p}
=\bigl(\widehat a_{\tau,p}(u_k,u_l)\bigr)_{k,l=1}^m.
\]

\begin{lemma}\label{lem:smoothing}
For every fixed $j\ge1$,
\begin{equation}\label{eq:smoothingnorm}
\sup_{p\in K_r}
|\nu_j(\tau,p,\delta)-\nu_j(\tau,p)|\longrightarrow0
\qquad(\delta\downarrow0).
\end{equation}
For all sufficiently small $\delta$, the spectral subspace for
$\nu_2(\tau,p,\delta),\ldots,\nu_{m+1}(\tau,p,\delta)$ has a
$\widehat b_{\tau,p,\delta}$-orthonormal frame
$u_{1,\delta}(p),\ldots,u_{m,\delta}(p)$, continuous in $p$ in
$H_0^1(\mathbb{B})$.
Its energy matrix
\[
\widehat B_{\tau,p,\delta}
=\bigl(\widehat a_{\tau,p,\delta}
(u_{k,\delta},u_{l,\delta})\bigr)_{k,l=1}^m
\]
is continuous in $p$ and satisfies
\begin{equation}\label{eq:smoothingmatrix}
\sup_{p\in K_r}
\|\widehat B_{\tau,p,\delta}-\widehat B_{\tau,p}\|
\longrightarrow0\qquad(\delta\downarrow0).
\end{equation}
\end{lemma}
\begin{proof}
Fix $\tau>0$. Let $\delta_n\to0$ and $p_n\to p$, and write
\[
\begin{gathered}
\widehat\rho_n:=\widehat\rho_{\tau,p_n,\delta_n},
\qquad
\widehat\rho:=\widehat\rho_{\tau,p},\\
a_n:=\widehat a_{\tau,p_n,\delta_n},\qquad
b_n:=\widehat b_{\tau,p_n,\delta_n},\qquad
a:=\widehat a_{\tau,p},\qquad
b:=\widehat b_{\tau,p}.
\end{gathered}\]
Continuity in $p$ and approximation by mollifiers give
$\widehat\rho_n\to\widehat\rho$ in measure. Since
$0<c_\tau\le\widehat\rho_n,\widehat\rho\le C$, we have
\[
\|(\widehat\rho_n-\widehat\rho)w\|_{L^2(\mathbb B)}
\longrightarrow0
\qquad\text{for every fixed }w\in L^2(\mathbb B).
\]
Thus $a_n,b_n$ converge to $a,b$ uniformly on bounded subsets
of every fixed finite-dimensional subspace of $H_0^1(\mathbb B)$.
The min--max principle implies
\[
\limsup_{n\to\infty}\nu_j(\tau,p_n,\delta_n)
\le\nu_j(\tau,p).
\]

Fix $j\ge1$ and choose $b_n$-orthonormal eigenfunctions
$f_{i,n}$ corresponding to
$\nu_{i,n}:=\nu_i(\tau,p_n,\delta_n)$, $1\le i\le j$.
These functions are bounded in $H_0^1(\mathbb B)$.
After passing to a subsequence,
\[
\nu_{i,n}\to\ell_i\le\nu_j(\tau,p) ,\qquad
f_{i,n}\rightharpoonup f_i\ \text{in }H_0^1(\mathbb B),
\qquad
f_{i,n}\to f_i\ \text{in }L^2(\mathbb B).
\]
Passing to the limit gives
\[
a(f_i,g)=\ell_i b(f_i,g),\qquad
b(f_i,f_l)=\delta_{il}.
\]
Since $\ell_1\le\cdots\le\ell_j$, the min--max principle on
$\operatorname{span}\{f_1,\ldots,f_i\}$ yields
$\nu_i(\tau,p)\le\ell_i$, therefore $\ell_i=\nu_i(\tau,p) $.
Moreover,
\[
\begin{aligned}
a_n(f_{i,n}-f_i,f_{i,n}-f_i)
&=\nu_{i,n}\bigl(1-2b_n(f_{i,n},f_i)\bigr)
+a_n(f_i,f_i)
\longrightarrow0.
\end{aligned}
\]
Since $\widehat\rho_n $ has uniformly lower bound, then $f_{i,n}\to f_i$ strongly in $H_0^1(\mathbb B)$.
Applying the same argument to every subsequence gives
$\nu_j(\tau,p_n,\delta_n)\to\nu_j(\tau,p)$.
Together with \eqref{eq:eigenvalue-continuity} and compactness
of $K_r$, this proves \eqref{eq:smoothingnorm}.

The uniform spectral gaps persist for small $\delta$.
Choose $\widehat b_{\tau,p,\delta}$-orthonormal eigenfunctions
$f_{i,\delta}(p)$ corresponding to $\nu_i(\tau,p,\delta)$,
$2\le i\le m+1$, and set
\[
v_{k,\delta}(p)
:=
\sum_{i=2}^{m+1}
\widehat b_{\tau,p,\delta}
\bigl(u_k(\tau,p),f_{i,\delta}(p)\bigr)f_{i,\delta}(p),
\qquad 1\le k\le m.
\]
These projections are independent of the chosen eigenbasis.
For $\delta_n\to0$ and $p_n\to p$, the preceding compactness
argument and Lemma~\ref{lem:framecontinuity} give, along a subsequence,
\[
v_{k,\delta_n}(p_n)
\longrightarrow
\sum_{i=2}^{m+1}b(u_k(\tau,p),f_i)f_i
=u_k(\tau,p)
\quad\text{in }H_0^1(\mathbb B),
\]
where $f_2,\ldots,f_{m+1}$ form a $b$-orthonormal basis
of the limiting spectral subspace.
Since every subsequence admits such a further subsequence,
compactness of $K_r$ yields
\[
\sup_{p\in K_r}\max_{1\le k\le m}
\|v_{k,\delta}(p)-u_k(\tau,p)\|_{H^1(\mathbb B)}
\longrightarrow0.
\]
For fixed $\delta$, the same argument and continuity of the
density in $L^\infty$ prove continuity of $v_{k,\delta}$ in $p$.
Define the continuous Gram matrix
\[
G_\delta(p)
:=
\bigl(
\widehat b_{\tau,p,\delta}
(v_{k,\delta}(p),v_{l,\delta}(p))
\bigr)_{k,l=1}^m.
\]
The preceding convergence and the multiplier convergence of
the densities imply $G_\delta\to I_m$ uniformly in $p$.
For small $\delta$, set
\[
u_{k,\delta}(p)
:=
\sum_{l=1}^m
v_{l,\delta}(p)\bigl(G_\delta(p)^{-1/2}\bigr)_{lk}.
\]
These functions form a continuous
$\widehat b_{\tau,p,\delta}$-orthonormal frame and converge
to $u_k(\tau,p)$ uniformly in $H^1(\mathbb B)$.
Consequently, along every $\delta_n\to0$ and $p_n\to p$,
\[
\begin{aligned}
\bigl(\widehat B_{\tau,p_n,\delta_n}\bigr)_{kl}
=a_n\bigl(u_{k,\delta_n}(p_n),u_{l,\delta_n}(p_n)\bigr)
\rightarrow
a(u_k(\tau,p),u_l(\tau,p))
=\bigl(\widehat B_{\tau,p}\bigr)_{kl}.
\end{aligned}
\]
Continuity of $\widehat B_{\tau,p}$ and compactness of $K_r$
give \eqref{eq:smoothingmatrix}.
For fixed $\delta$, continuity of
$\widehat B_{\tau,p,\delta}$ follows from that of the
density and the frame.
\end{proof}

	\begin{proof}[\bf Proof of Proposition~\ref{prop:potential}] By using Lemma~\ref{lem:weighted-dirichlet}, \eqref{eq:weighted-dirichlet-bounds}, \eqref{eq:nu2nu1} and $r<(m+1)/4$, we can 
		take $\tau>0$ sufficiently small such that the following estimates hold uniformly for $p$:
		\[
		\begin{gathered}
			\nu_1(\tau,p)<\frac{\tau(m+1)}8,\quad
			\frac{5\tau(m+1)}8<\nu_2(\tau,p)\le\nu_{m+1}(\tau,p)<\frac{11\tau(m+1)}8,\\
		\nu_{m+2}(\tau,p)>3\tau(m+1),\quad	\|\widehat B_{\tau,p}-\tau((m+1)I_H+p)\|<\frac{\tau r}{8}.
		\end{gathered}
		\]
		By Lemma~\ref{lem:smoothing}, we choose one $\delta>0$ such that,
		for every $p$,
		\[
		\begin{gathered}
			\max_{1\le j\le m+2}
			|\nu_j(\tau,p,\delta)-\nu_j(\tau,p)|
			<\frac{\tau(m+1)}{16},\\
			\|\widehat B_{\tau,p,\delta}-\widehat B_{\tau,p}\|
			<\frac{\tau r}{8}.
		\end{gathered}
		\]
		Set
		\begin{equation}\label{eq:Vformula}
\phi_p:=\sqrt{\widehat\rho_{\tau,p,\delta}}>0,\qquad	V_p:=\frac{\Delta\phi_p}{\phi_p}.
		\end{equation}
		Multiplication by $\phi_p$ is a bijection of $H_0^1(\mathbb{B})$.
		For $u=\phi_pf$, integration by parts and density give
		\begin{equation}\label{eq:groundtransform}
			\begin{aligned}
				\int_{\mathbb{B}}\widehat\rho_{\tau,p,\delta}|\nabla f|^2\dd
				&=\int_{\mathbb{B}}|\nabla u-u\nabla\log\phi_p|^2\dd\\
				&=\int_{\mathbb{B}}\bigl(|\nabla u|^2+V_pu^2\bigr)\dd,
			\end{aligned}
		\end{equation}
		while $\int_{\mathbb{B}}\widehat\rho_{\tau,p,\delta}f^2\dd=\int_{\mathbb{B}}u^2\dd$. By polarization of \eqref{eq:groundtransform}, for any
		$f,g\in H_0^1(\mathbb{B})$, we have
		\begin{align*}
			\widehat a_{\tau,p,\delta}(f,g)
			&=\int_{\mathbb{B}}\left(
			\nabla(\phi_p f)\cdot\nabla(\phi_p g)
			+V_p\phi_p^2fg\right)\dd,\\
			\widehat b_{\tau,p,\delta}(f,g)
			&=\int_{\mathbb{B}}(\phi_p f)(\phi_p g)\dd.
		\end{align*}
		Consequently, multiplication by $\phi_p$ preserves the
		eigenvalues and their multiplicities, namely
		the eigenvalues of $L_p=-\Delta_{\mathbb B}^{D}+V_p$ satisfies
        \[
		\xi_j(p)=\nu_j(\tau,p,\delta),\qquad j\ge1.
		\]
		Moreover, $\{\phi_pu_{k,\delta}\}_{k=1}^m$ is a continuous
		orthonormal frame in $L^2(\mathbb{B})$, and its operator matrix satisfies
		\[
		\mathcal B_p
		=\bigl(\widehat a_{\tau,p,\delta}
		(u_{k,\delta},u_{l,\delta})\bigr)_{k,l=1}^m
		=\widehat B_{\tau,p,\delta}.
		\]
		
		Since $\phi_p$ is uniformly positive and continuous in $p$
		in every $C^k$ norm, so
		$V_p=\Delta\phi_p/\phi_p$ also depends continuously on \(p\) in every $C^k$ norm. Furthermore,
		\[
		\xi_1(p)
		=\inf_{0\ne f\in H_0^1(\mathbb{B})}
		\frac{\int_{\mathbb{B}}\widehat\rho_{\tau,p,\delta}|\nabla f|^2\dd}
		{\int_{\mathbb{B}}\widehat\rho_{\tau,p,\delta}f^2\dd}
		\ge\frac{c_\tau}{C}\lambda_1(\mathbb{B})>0.
		\]
		The preceding eigenvalue estimates imply
		\eqref{eq:potentialgaps}, while
		\begin{align*}
			&\|\mathcal B_p-\tau((m+1)I_H+p)\|\\
			&\quad\le
			\|\widehat B_{\tau,p,\delta}-\widehat B_{\tau,p}\|
			+\|\widehat B_{\tau,p}-\tau((m+1)I_H+p)\|
			<\frac{\tau r}{4},
		\end{align*}
		which proves \eqref{eq:potentialmatrix}.
	\end{proof}

\section{Approximation by radial cuts}\label{sec:radialconstruction}

In Section~\ref{sec:potential}, we constructed a continuous family
of smooth potentials $V_p$ on the fixed ball $\mathbb{B}$ such that the Schr\"odinger operator
$L_p=-\Delta_{\mathbb{B}}^D+V_p$ has an isolated $m$-dimensional spectral
subspace whose operator matrix satisfies
\[
\sup_{p\in K_r}
\|\mathcal B_p-\tau((m+1)I_H+p)\|
<\frac{\tau r}{4}.
\]
However, our main theorem concerns the ordinary Dirichlet Laplacian,
where only the domain is allowed to vary. The goal of this section
is therefore to realize the spectral effect of $V_p$, up to a
constant shift, through the geometry of smooth strictly
star-shaped domains. The traditional approach of removing isolated interior holes (see, e.g., \cite{CM,MM95,Jing,KP,RT}) can approximate Schr\"odinger operators, but it generally does not preserve the star-shapedness of the domain. The present construction is
different, as the effective potential is generated
by geometric radial cuts while preserving
star-shapedness.

The main idea is illustrated in Figure~\ref{fig:radial-idea-simple}.
We first add a constant $C_*$ so that $W_p=V_p+C_*>0$
and $r^2W_p(r\theta)$ is increasing in $r$. We then partition
$\mathbb S^{N-1}$ into small angular cells $Q_{j,\eps}$ and prescribe
the capacity of each radial cut according to
\[
\alpha_jr^2W_p(r\theta_j),
\qquad \alpha_j=|Q_{j,\eps}|.
\]
The resulting cuts widen outward, preserving the star-shapedness
of the remaining domain $G_{\eps,p}$. The key point is that the Dirichlet problem for $-\Delta^D_{G_{\eps,p}}$ in $G_{\eps,p}$ converges, as $\eps\to0$, to the Dirichlet problem for the Schr\"odinger operator $-\Delta_{\mathbb{B}}^D+W_p$ in $\mathbb{B}$; see Proposition \ref{prop:roughnorm}.
Finally, smoothing the radial boundary gives smooth strictly
star-shaped domains $\Omega_{k,p}$, whose Dirichlet Laplacians
approximate $-\Delta_{\mathbb{B}}^D+W_p$ in the operator norm sense, see Theorem \ref{thm:radial}.

\begin{figure}[!htbp]
\centering
\resizebox{0.97\linewidth}{!}{%
\begin{tikzpicture}[
    x=1cm,y=1cm,
    font=\footnotesize,
    line cap=round,line join=round,
    >=Stealth
]
\definecolor{figink}{RGB}{37,56,75}
\definecolor{figslit}{RGB}{225,133,61}
\path[use as bounding box] (0,0) rectangle (12.60,4.40);

\node[font=\small\bfseries] at (1.95,4.12)
    {(a) Smooth potential};
\node[font=\small\bfseries] at (6.30,4.12)
    {(b) Radial cuts};
\node[font=\small\bfseries] at (10.65,4.12)
    {(c) Smoothed domain};

\draw[-{Stealth[length=2mm]},line width=.85pt,figink]
   (3.43,2.51)--(4.53,2.51);
\node[font=\scriptsize] at (3.98,2.82)
   {capacity};

\draw[-{Stealth[length=2mm]},line width=.85pt,figink]
   (7.78,2.51)--(8.88,2.51);
\node[font=\scriptsize] at (8.33,2.82)
   {smoothing};

\begin{scope}[shift={(1.95,2.48)},scale=1.10]
  \fill[blue!7] (0,0) circle (1.11);
  \draw[figink,line width=.85pt] (0,0) circle (1.11);
\end{scope}

\begin{scope}[shift={(6.30,2.48)},scale=1.10]
  \fill[blue!7] (0,0) circle (1.11);
  \draw[figink,line width=.85pt] (0,0) circle (1.11);

  \foreach \a in {30,90,150,210,270,330}{
    \begin{scope}[rotate=\a]
      \path[fill=white]
        (1.28,-.22)--(.35,-.015)
        arc[start angle=-2.45,end angle=2.45,radius=.35]
        --(1.28,.22)--cycle;

      \draw[figslit,line width=.75pt]
        (1.095,-.188)--(.35,-.015)
        arc[start angle=-2.45,end angle=2.45,radius=.35]
        --(1.095,.188);
    \end{scope}
  }

  \draw[figink!65,densely dashed,line width=.55pt]
    (0,0) circle (.27);
  \fill[figink] (0,0) circle (1pt);

  \node[fill=white,inner sep=1pt,font=\scriptsize]
    at (0,-.04) {$0$};
  \draw[figink!70,thin] (-.19,-.19)--(-.40,-.42);
  \node[anchor=north east,inner sep=0pt,font=\scriptsize]
    at (-.4,-.42) {$\mathbb{B}_\varepsilon$};
\end{scope}

\begin{scope}[shift={(10.65,2.48)},scale=1.10]
  \path[fill=blue!7,draw=figink,line width=.85pt]
    plot[domain=0:360,samples=520,variable=\t]
    ({(1.11-.70*exp(-8*(1-cos(6*(\t-30)))))*cos(\t)},
     {(1.11-.70*exp(-8*(1-cos(6*(\t-30)))))*sin(\t)})
    --cycle;

  \draw[figink!65,densely dashed,line width=.55pt]
    (0,0) circle (.27);
  \fill[figink] (0,0) circle (1pt);
\end{scope}

\node at (1.95,.91)
  {$-\Delta_{\mathbb{B}}^D+W_p$};

\node at (6.30,.91)
  {$G_{\varepsilon,p}
    =\mathbb{B}\setminus\bigcup_jK_{j,\varepsilon,p}$};

\node at (10.65,.91)
  {$\Omega_{k,p}$ smooth, strictly star-shaped };
\end{tikzpicture}%
}\vspace{-8mm}
\caption{}
\label{fig:radial-idea-simple}
\end{figure}

For an open set $\Omega\subset \mathbb{B}$, let
$E_\Omega:L^2(\Omega)\to L^2(\mathbb{B})$ denote extension by zero.
 Given \(f\in L^2(\Omega)\), it is defined by
\begin{align*}
    E_{\Omega}f(z)=\begin{cases}
        \displaystyle f(z) \qquad&\mathrm{if}\qquad z\in \Omega,\\
        \displaystyle 0 \qquad&\mathrm{if}\qquad z\in \mathbb{B}\backslash \Omega.
    \end{cases}
\end{align*}
Let \(E^*_{\Omega}:L^2(\mathbb{B})\to L^2(\Omega)\) denote its adjoint operator. Then \(E^*_{\Omega}\) is simply the restriction operator, namely,
$
E^*_{\Omega}f=f|_{\Omega}.
$

\begin{theorem}\label{thm:radial}
Let \(N\ge3\), and let
\(p\mapsto V_p\) be continuous from \(K_r\) to
\(C^1(\overline{\mathbb{B}};\R)\). Set
\begin{align}
 C_*&=\sup_{p\in K_r}\norm{V_p}_{L^\infty(\mathbb{B})}+\tfrac12\sup_{p\in K_r}\norm{x\cdot\nabla V_p}_{L^\infty(\mathbb{B})}+1,& W_p&=V_p+C_*.
 \label{eq:shift}
\end{align}
Then there are domains \(\Omega_{k,p}\subset \mathbb{B}\) of the form
\eqref{eq:strictstar} such that, for each \(k\), their radial
functions depend continuously on \(p\) in \(C^\infty(\Sph^{N-1})\),
and
\begin{equation}\label{eq:radialnorm}
 \sup_{p\in K_r}
 \left\|
 E_{\Omega_{k,p}}(-\Delta_{\Omega_{k,p}}^D+1)^{-1}E_{\Omega_{k,p}}^*
 -(-\Delta_{\mathbb{B}}^D+W_p+1)^{-1}
 \right\|\to0\quad(k\to+\infty).
\end{equation}
Every fixed eigenvalue and every continuous, uniformly isolated
finite spectral projection converge uniformly in \(p\).
Whenever the limiting subspace has a continuous orthonormal
frame, the operator matrices in the transported frames also
converge uniformly.
\end{theorem}

We now prove Theorem~\ref{thm:radial}. Note that $N-1\ge2$ and
$W_p=V_p+C_*$. The choice of $C_*$ in \eqref{eq:shift} gives
\begin{equation}\label{eq:Wconditions}
W_p\ge1,\qquad 2W_p+x\cdot\nabla W_p\ge2.
\end{equation}
Define $A_p(r,\theta):=r^2W_p(r\theta)$. Uniformly in $p$, for $0<r\le1$,
\begin{equation}\label{eq:Abounds}
\begin{gathered}
r^2\le A_p\le Cr^2,\qquad
2r\le\partial_rA_p\le Cr,\\
|\nabla_{\Sph^{N-1}}A_p|\le Cr^3,\qquad
\frac{\partial_rA_p}{A_p}\le\frac Cr.
\end{gathered}
\end{equation}
All constants in this section are
uniform in $p$; the dependence on a fixed $\eps$ is indicated explicitly.

\subsection{Radial cuts calibrated by capacity}\label{sec:cutcalibration}

Choose $\eps\downarrow0$ and partitions $\{Q_{j,\eps}\}_j$ of
$\Sph^{N-1}$, up to measure zero, with centers $\theta_j$ such that
\begin{equation}\label{eq:cells}
\begin{gathered}
c\eps^{N-1}\le\alpha_j=\alpha_{j,\eps}:=|Q_{j,\eps}|\le C\eps^{N-1},\qquad
\operatorname{diam}Q_{j,\eps}\le C\eps,\\
B_{\Sph^{N-1}}(\theta_j,2b_\eps)\subset Q_{j,\eps},\qquad
b_\eps=c_0\eps,\qquad \sum_j\alpha_j=|\Sph^{N-1}|,
\end{gathered}
\end{equation}
where $B_{\Sph^{N-1}}(\theta_j,r) $ is the geodesic ball on $\mathbb{S}^{N-1}$ with center $\theta_j$,
$$B_{\Sph^{N-1}}(\theta_j,r):=\{\theta\in \mathbb{S}^{N-1}: \dist_{\Sph^{N-1}}(\theta,\theta_j)<r\}.$$
Such cells are obtained by subdividing the faces of
$[-1,1]^{N}$ into cubes of side $2/k$ and projecting radially onto
$\Sph^{N-1}$ (i.e., $x\mapsto x/|x|$), with $\eps=1/k$, see Figure \ref{F.2}. The face projections and their inverses
have uniformly bounded derivatives. Consequently, the rescaled
cells have uniform Lipschitz charts, and $c_0>0$ can be fixed small
enough to ensure \eqref{eq:cells}.
\begin{figure}[htbp]
\centering
\begingroup
\def\Cpt#1#2#3{({0.72*((#1)-(#2))},{0.415*(2*(#3)-(#1)-(#2))})}
\def\Spt#1#2#3{({1.55*((#1)-(#2))/sqrt(2*((#1)^2+(#2)^2+(#3)^2))},{1.55*(2*(#3)-(#1)-(#2))/sqrt(6*((#1)^2+(#2)^2+(#3)^2))})}

\resizebox{0.82\linewidth}{!}{%
\begin{tikzpicture}[
    font=\footnotesize,
    line join=round,
    line cap=round
]
\path[use as bounding box] (0,0) rectangle (10.8,5.13);
\definecolor{cubeink}{RGB}{46,72,93}
\definecolor{cellorange}{RGB}{222,124,51}

\node[font=\small\bfseries] at (2.05,4.73)
    {Subdivide the cube faces};
\node[font=\small\bfseries] at (8.20,4.73)
    {Radially project to $\mathbb S^{N-1}$};

\begin{scope}[shift={(2.05,2.60)}]

  \path[fill=blue!8]
    \Cpt{-1}{1}{1} -- \Cpt{1}{1}{1}
    -- \Cpt{1}{1}{-1} -- \Cpt{-1}{1}{-1} --cycle;

  \path[fill=blue!14]
    \Cpt{1}{1}{1} -- \Cpt{1}{-1}{1}
    -- \Cpt{1}{-1}{-1} -- \Cpt{1}{1}{-1} --cycle;

  \path[fill=blue!4]
    \Cpt{-1}{-1}{1} -- \Cpt{1}{-1}{1}
    -- \Cpt{1}{1}{1} -- \Cpt{-1}{1}{1} --cycle;

  \foreach \h in {-0.5,0,0.5} {
    \draw[cubeink!55,line width=.38pt]
      \Cpt{\h}{-1}{1} -- \Cpt{\h}{1}{1};
    \draw[cubeink!55,line width=.38pt]
      \Cpt{-1}{\h}{1} -- \Cpt{1}{\h}{1};

    \draw[cubeink!55,line width=.38pt]
      \Cpt{1}{\h}{-1} -- \Cpt{1}{\h}{1};
    \draw[cubeink!55,line width=.38pt]
      \Cpt{1}{-1}{\h} -- \Cpt{1}{1}{\h};

    \draw[cubeink!55,line width=.38pt]
      \Cpt{\h}{1}{-1} -- \Cpt{\h}{1}{1};
    \draw[cubeink!55,line width=.38pt]
      \Cpt{-1}{1}{\h} -- \Cpt{1}{1}{\h};
  }

  \path[fill=cellorange!55,draw=cellorange,line width=.85pt]
    \Cpt{0}{-.5}{1} -- \Cpt{.5}{-.5}{1}
    -- \Cpt{.5}{0}{1} -- \Cpt{0}{0}{1} -- cycle;

  \draw[cubeink,line width=.9pt]
    \Cpt{-1}{-1}{1} -- \Cpt{1}{-1}{1}
    -- \Cpt{1}{-1}{-1} -- \Cpt{1}{1}{-1}
    -- \Cpt{-1}{1}{-1} -- \Cpt{-1}{1}{1} -- cycle;

  \draw[cubeink,line width=.9pt]
    \Cpt{1}{1}{1} -- \Cpt{1}{-1}{1};
  \draw[cubeink,line width=.9pt]
    \Cpt{1}{1}{1} -- \Cpt{-1}{1}{1};
  \draw[cubeink,line width=.9pt]
    \Cpt{1}{1}{1} -- \Cpt{1}{1}{-1};

  \draw[cubeink!55,dashed,line width=.45pt]
    \Cpt{-1}{-1}{1} -- \Cpt{-1}{-1}{-1}
    -- \Cpt{1}{-1}{-1};

  \draw[cubeink!55,dashed,line width=.45pt]
    \Cpt{-1}{-1}{-1} -- \Cpt{-1}{1}{-1};

\end{scope}

\draw[->,line width=.85pt,cubeink]
    (4.35,2.63) -- (5.91,2.63);

\node[align=center,fill=white,inner sep=1pt]
    at (5.13,3.17) {$x\mapsto x/|x|$};

\node[align=center,inner sep=1pt]
    at (5.13,2.29) {radial projection};

\begin{scope}[shift={(8.20,2.60)}]

  \shade[inner color=white,outer color=blue!13]
    (0,0) circle (1.55);

  \foreach \h in {-1,-.5,0,.5,1} {
    \pgfmathsetmacro{\lo}{max(-1,-1-\h)}

    \draw[cubeink!40,line width=.36pt]
      plot[domain=\lo:1,samples=20,variable=\t]
      ({1.55*(\h-\t)/sqrt(2*(\h*\h+\t*\t+1))},
       {1.55*(2-\h-\t)/sqrt(6*(\h*\h+\t*\t+1))});

    \draw[cubeink!40,line width=.36pt]
      plot[domain=\lo:1,samples=20,variable=\t]
      ({1.55*(\t-\h)/sqrt(2*(\h*\h+\t*\t+1))},
       {1.55*(2-\t-\h)/sqrt(6*(\h*\h+\t*\t+1))});

    \draw[cubeink!40,line width=.36pt]
      plot[domain=\lo:1,samples=20,variable=\t]
      ({1.55*(1-\h)/sqrt(2*(1+\h*\h+\t*\t))},
       {1.55*(2*\t-1-\h)/sqrt(6*(1+\h*\h+\t*\t))});

    \draw[cubeink!40,line width=.36pt]
      plot[domain=\lo:1,samples=20,variable=\t]
      ({1.55*(1-\t)/sqrt(2*(1+\h*\h+\t*\t))},
       {1.55*(2*\h-1-\t)/sqrt(6*(1+\h*\h+\t*\t))});

    \draw[cubeink!40,line width=.36pt]
      plot[domain=\lo:1,samples=20,variable=\t]
      ({1.55*(\h-1)/sqrt(2*(1+\h*\h+\t*\t))},
       {1.55*(2*\t-1-\h)/sqrt(6*(1+\h*\h+\t*\t))});

    \draw[cubeink!40,line width=.36pt]
      plot[domain=\lo:1,samples=20,variable=\t]
      ({1.55*(\t-1)/sqrt(2*(1+\h*\h+\t*\t))},
       {1.55*(2*\h-1-\t)/sqrt(6*(1+\h*\h+\t*\t))});
  }

  \path[fill=cellorange!55,draw=cellorange,line width=.9pt]
    \Spt{0}{-.5}{1} -- \Spt{.5}{-.5}{1}
    -- \Spt{.5}{0}{1} -- \Spt{0}{0}{1} -- cycle;

  \draw[cubeink,line width=.9pt]
    (0,0) circle (1.55);

  \draw[->,cubeink!85,line width=.55pt]
    \Spt{.25}{-.25}{1} -- (1.42,1.64);

  \node[anchor=west,inner sep=1pt]
    at (1.43,1.64) {$Q_{j,\varepsilon}$};

\end{scope}

\node[align=center] at (2.05,.37)
    {face squares of side $2/k$};

\node[align=center] at (8.20,.37)
    {$|Q_{j,\varepsilon}|\simeq\varepsilon^{N-1},
      \quad\varepsilon=1/k$};

\end{tikzpicture}%
}
\endgroup\vspace{-8mm}
\caption{}
\label{F.2}
\end{figure}

For this cell decomposition, we have the following basic trace estimate, which will be used in Lemma \ref{lem:angularflux}.
\begin{lemma}\label{lem:celltrace}
For $g\in H^1(Q_{j,\eps})$, there holds
\begin{equation}\label{eq:celltrace}
\left|\avint_{\partial B_{\Sph^{N-1}}(\theta_j,b_\eps)}g\,dS
       -\avint_{Q_{j,\eps}}g\ds\right|
\le C\eps^{1-(N-1)/2}\|\nabla_{\Sph^{N-1}}g\|_{L^2(Q_{j,\eps})}.
\end{equation}
\end{lemma}

\begin{proof}
Let $g_Q=\avint_{Q_{j,\eps}}g\ds$. The rescaled Poincar\'e inequality gives
\[
\|g-g_Q\|^2_{L^2(Q_{j,\eps})}\le C\eps^2\|\nabla_{\Sph^{N-1}}g\|^2_{L^2(Q_{j,\eps})}.
\]
Applying the rescaled Sobolev trace inequality on $B_{\Sph^{N-1}}(\theta_j,b_\eps)\subset Q_{j,\eps}$ leads to 
\begin{align*}
    \|g-g_Q\|_{L^2(\partial B_{\Sph^{N-1}}(\theta_j,b_\eps))}^2
\le& C(\eps^{-1}\|g-g_Q\|_{L^2(B_{\Sph^{N-1}}(\theta_j,b_\eps))}^2+\eps\|\nabla_{\Sph^{N-1}}g\|_{L^2(B_{\Sph^{N-1}}(\theta_j,b_\eps))}^2 )\\
\leq&C\eps \|\nabla_{\Sph^{N-1}}g\|_{L^2(Q_{j,\eps})}^2.
\end{align*}
Therefore, we can obtain
\begin{align*}
    \left|\avint_{\partial B_{\Sph^{N-1}}(\theta_j,b_\eps)}(g-g_Q)\,dS \right|\le& |\partial B_{\Sph^{N-1}}(\theta_j,b_\eps)|^{-1/2}\|g-g_Q\|_{L^2(\partial B_{\Sph^{N-1}}(\theta_j,b_\eps))}\\
    \le& C\eps^{1-(N-1)/2}\|\nabla_{\Sph^{N-1}}g\|_{L^2(Q_{j,\eps})}.
\end{align*}
\end{proof}

\textbf{The capacity of cylinder on $\Sph^{N-1}$}: 
Write 
\begin{equation}\label{eq:Id}
I_{N-1}(\rho;b):=\int_\rho^b(\sin s)^{2-N}\,ds,
\qquad 0<\rho<b<\pi/2.
\end{equation}
On a spherical annulus $$C_{\rho,b}:=B_{\Sph^{N-1}}(\theta,b)\backslash \overline{B_{\Sph^{N-1}}(\theta,\rho)},$$ the function
$h(s):=1-I_{N-1}(s;b)/I_{N-1}(\rho;b)$ is harmonic, namely,
\begin{align*}
  \Delta_{\Sph^{N-1}}h=0 \quad\mathrm{in} \quad  C_{\rho,b},\qquad h=0 \quad\mathrm{on} \quad  s=\rho,\qquad h=1 \quad\mathrm{on} \quad  s=b.
\end{align*}   Direct integration gives its energy
\begin{align}\label{eq:capformula}
\int_{C_{\rho,b}}|\nabla_{\Sph^{N-1}}h|^2\ds&=\int_{\mathbb S^{N-2}}\int_\rho^b|h'(s)|^2(\sin s)^{N-2}dsd\omega
=\frac{|\Sph^{N-2}|}{I_{N-1}(\rho;b)}\notag\\
&=\inf_{\substack{u\in H^1(C_{\rho,b})\\ u|_{s=\rho}=0,u|_{s=b}=1}}\int_{C_{\rho,b} }|\nabla_{\Sph^{N-1}}u|^2\ds .
\end{align}
Thus \eqref{eq:capformula} is the relative capacity of the annulus.

Recall \eqref{eq:Abounds}-\eqref{eq:cells} that $A_p(r,\theta)=r^2W_p(r\theta)$ and  $b_\eps=c_0\eps$. 
For each $j,p,r$, we can  define $\rho_{j,p}(r) $ as follows:
\begin{equation}\label{eq:kappacut}
I_{N-1}(\rho_{j,p}(r);b_\eps)=\frac{|\Sph^{N-2}|}{\alpha_jA_p(r,\theta_j)}.
\end{equation}
Since for $N\ge3$, $I_{N-1}(\rho;b_\eps)$ decreases continuously in $(0,b_{\eps})$ and
\begin{align*}
    \lim_{\rho\to 0}I_{N-1}(\rho;b_\eps)=+\infty,\qquad \lim_{\rho\to b_{\eps}}I_{N-1}(\rho;b_\eps)=0.
\end{align*}
Thus $\rho_{j,p}(r)\in(0,b_\eps)$ exists uniquely. Differentiation gives
\begin{equation}\label{eq:rhomonotone}
\rho_{j,p}'(r)=|\Sph^{N-2}|
\frac{A'_p(r,\theta_j)}{\alpha_jA^2_p(r,\theta_j)}
(\sin\rho_{j,p}(r))^{N-2}>0.
\end{equation}
This follows from \eqref{eq:Abounds} and ensures that the resulting
domains are star-shaped.

{\bf The domains before smoothing}:
For $0<\eps<1/4$, we remove the cuts
\begin{equation}\label{eq:removed}
K_{j,\eps,p}:=
\{r\theta:\eps\le r<1,\
\dist_{\Sph^{N-1}}(\theta,\theta_j)\le\rho_{j,p}(r)\},
\end{equation}
and define
\begin{align}\label{Def G}
    G_{\eps,p}:=\mathbb{B}\setminus\bigcup_jK_{j,\eps,p}.
\end{align}
Their locations are listed by the following picture:
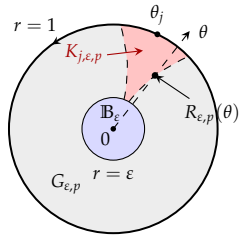
\begin{center}
\begin{tikzpicture}[
    scale=0.69,
    >=stealth,
    every node/.style={font=\scriptsize}
]

\fill[black!7] (0,0) circle (2);

\path[fill=red!23]
    (60:0.60)
    .. controls (58:0.85) and (54:1.10) .. (52:1.30)
    .. controls (49:1.55) and (46:1.82) .. (44:2)
    arc (44:86:2)
    .. controls (84:1.82) and (81:1.55) .. (78:1.30)
    .. controls (76:1.10) and (72:0.85) .. (70:0.60)
    arc (70:60:0.60) -- cycle;

\draw[dashed,thin]
    (60:0.60)
    .. controls (58:0.85) and (54:1.10) .. (52:1.30)
    .. controls (49:1.55) and (46:1.82) .. (44:2);

\draw[dashed,thin]
    (70:0.60)
    .. controls (72:0.85) and (76:1.10) .. (78:1.30)
    .. controls (81:1.55) and (84:1.82) .. (86:2);

\draw[line width=0.75pt] (0,0) circle (2);
\fill[blue!15] (0,0) circle (0.60);
\draw (0,0) circle (0.60);

\draw[->,dashed] (0,0) -- (52:2.38);
\fill (0,0) circle (1.5pt);
\fill (52:1.30) circle (1.8pt);
\fill (65:2) circle (1.6pt);

\node[below left,inner sep=1pt] at (0,0) {$0$};
\node at (-0.05,0.30) {$\mathbb{B}_{\varepsilon}$};
\node[below] at (0,-0.60) {$r=\varepsilon$};

\node[above] at (65:2) {$\theta_j$};
\node[right] at (52:2.38) {$\theta$};

\node at (-1.5,1.99) {$r=1$};
\draw[->,thin] (-1.14,1.64) -- (126:2);

\node[red!65!black] at (-0.56,1.43)
    {$K_{j,\varepsilon,p}$};
\draw[->,thin,red!65!black]
    (-0.12,1.42) -- (0.57,1.50);

\node at (-0.90,-1.05) {$G_{\varepsilon,p}$};

\node[anchor=west] at (1.18,0.29)
    {$R_{\varepsilon,p}(\theta)$};
\draw[->,thin] (1.46,0.58) -- (52:1.30);

\end{tikzpicture}
\end{center}

\textbf{Representation of star-shaped}: For $\theta\in \mathbb S^{N-1}$, we denote
$$s(\theta):=\min_{j}\dist_{\Sph^{N-1}}(\theta,\theta_j).$$
We call $B_{\Sph^{N-1}}(\theta_j,b_{\eps})$ the $j$-th cap. The cuts lie in disjoint caps and $\mathbb{B}_\eps\subset G_{\eps,p}$.
By \eqref{eq:rhomonotone}, every ray meets $G_{\eps,p}$ in an initial
interval. In the $j$-th cap, we have
$s(\theta)=\dist_{\Sph^{N-1}}(\theta,\theta_j)$, and its radial function is
\begin{equation}\label{eq:roughR}
R_{\eps,p}(\theta)=
\begin{cases}
\eps,&s(\theta)\le\rho_{j,p}(\eps),\\
\rho_{j,p}^{-1}(s(\theta)),&\rho_{j,p}(\eps)<s(\theta)<\rho_{j,p}(1),\\
1,&b_{\eps}>s(\theta)\ge\rho_{j,p}(1).
\end{cases}
\end{equation}
Outside the caps, it equals one; equivalently, $R_{\eps,p}(\theta)=1$ for $s(\theta)\ge b_{\eps}$.
For fixed $\eps$, compactness gives
\[
c_\eps:=\min_{j,p,\ \eps\le r\le1}\rho_{j,p}'(r)>0.
\]
The clipped inverse in \eqref{eq:roughR} is Lipschitz, and comparison
of two such inverses gives
\begin{equation}\label{eq:radialparameter}
\|R_{\eps,p}-R_{\eps,p'}\|_{L^{\infty}(\Sph^{N-1})}
\le c_\eps^{-1}\max_j
\sup_{\eps\le r\le1}|\rho_{j,p}(r)-\rho_{j,p'}(r)|.
\end{equation}
Thus $G_{\eps,p}$ is a Lipschitz star-shaped domain (see formula \eqref{eq:normal}),
\begin{align*}
    G_{\eps,p}=\left\{r\theta: 0\leq r< R_{\eps,p}(\theta)\right\},
\end{align*}
and $p\mapsto R_{\eps,p}$ is continuous
in the uniform norm. Remark that $c_\eps^{-1}$ may tend to infinity as $\eps\downarrow0$.

\subsection{An important auxiliary function and estimates}\label{sec:correctorpart}

 Recall that $s(\theta)=\min_j\dist_{\Sph^{N-1}}(\theta,\theta_j)$. For every $r>0$, we set $w_{\eps,p}(r,\theta):=1$ whenever $s(\theta)\geq b_{\eps}$. For $s(\theta)<b_{\eps}$, i.e., $\theta\in B_{\Sph^{N-1}}(\theta_j,b_{\eps})$ for some $j$, we have $s(\theta)=\dist_{\Sph^{N-1}}(\theta,\theta_j)$, and we define
\begin{equation}\label{eq:wdef}
w_{\eps,p}(r,\theta):=
\begin{cases}
0,&s(\theta)\le\rho_{j,p}(r),\\
1-\dfrac{I_{N-1}(s(\theta);b_\eps)}{I_{N-1}(\rho_{j,p}(r);b_\eps)},
 &\rho_{j,p}(r)<s(\theta)<b_\eps.
\end{cases}
\end{equation}
The precise values of $w_{\eps,p}$ are summarized in the following figure:

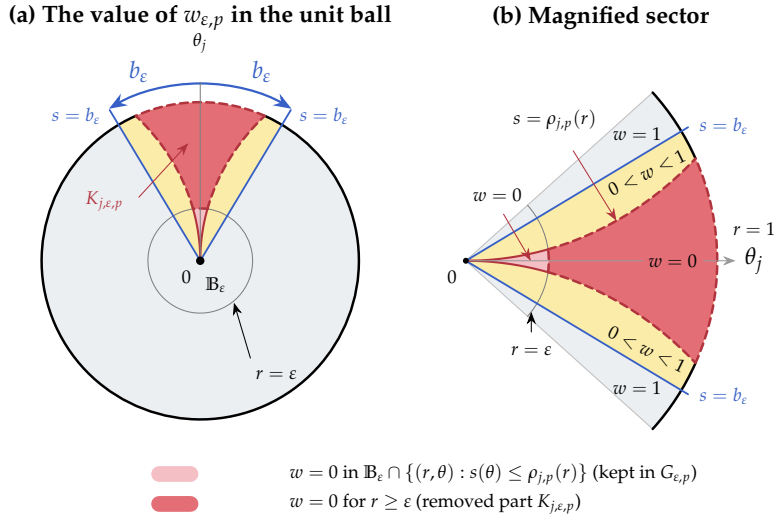
\begin{center}
\resizebox{0.82\linewidth}{!}{%
\begin{tikzpicture}[
    font=\small, >=Stealth,
    line cap=round, line join=round,
    removed/.style={borderred,line width=.95pt,
                    dash pattern=on 3pt off 2pt}
]
\definecolor{trans}{RGB}{252,236,169}
\definecolor{zero}{RGB}{227,110,118}
\definecolor{zerocore}{RGB}{244,192,197}
\definecolor{keep}{RGB}{236,239,242}
\definecolor{capblue}{RGB}{54,96,198}
\definecolor{borderred}{RGB}{177,51,61}

\def\EE{0.33}
\def\RMAX{24}
\def\CAP{31}

\begin{scope}[shift={(2.8,0)}]
\def\RR{2.12}

\node[font=\bfseries\small]
    at (0,3.25) {(a) The value of $w_{\eps,p}$ in the unit ball};

\fill[keep] (0,0) circle (\RR);
\fill[trans]
    (0,0)--({90-\CAP}:\RR)
    arc[start angle={90-\CAP},end angle={90+\CAP},radius=\RR]
    --cycle;

\path[fill=zero]
    (0,0)
    -- plot[domain=0:1,samples=120,variable=\t]
      ({\RR*\t*cos(90-\RMAX*\t)},
       {\RR*\t*sin(90-\RMAX*\t)})
    arc[start angle={90-\RMAX},end angle={90+\RMAX},radius=\RR]
    -- plot[domain=1:0,samples=120,variable=\t]
      ({\RR*\t*cos(90+\RMAX*\t)},
       {\RR*\t*sin(90+\RMAX*\t)})
    --cycle;

\path[fill=zerocore]
    (0,0)
    -- plot[domain=0:\EE,samples=80,variable=\t]
      ({\RR*\t*cos(90-\RMAX*\t)},
       {\RR*\t*sin(90-\RMAX*\t)})
    arc[start angle={90-\RMAX*\EE},
        end angle={90+\RMAX*\EE},radius={\RR*\EE}]
    -- plot[domain=\EE:0,samples=80,variable=\t]
      ({\RR*\t*cos(90+\RMAX*\t)},
       {\RR*\t*sin(90+\RMAX*\t)})
    --cycle;

\draw[line width=.95pt]
    ({90+\RMAX}:\RR)
    arc[start angle={90+\RMAX},end angle={450-\RMAX},radius=\RR];

\draw[black!50,thin]
    ({90+\RMAX*\EE}:{\RR*\EE})
    arc[start angle={90+\RMAX*\EE},
        end angle={450-\RMAX*\EE},radius={\RR*\EE}];

\foreach \sgn in {-1,1} {
    \draw[borderred,line width=.8pt]
      plot[domain=0:\EE,samples=80,variable=\t]
        ({\RR*\t*cos(90+\sgn*\RMAX*\t)},
         {\RR*\t*sin(90+\sgn*\RMAX*\t)});
    \draw[removed]
      plot[domain=\EE:1,samples=120,variable=\t]
        ({\RR*\t*cos(90+\sgn*\RMAX*\t)},
         {\RR*\t*sin(90+\sgn*\RMAX*\t)});
    \draw[capblue,line width=.7pt]
      (0,0)--({90+\sgn*\CAP}:2.35);
    \draw[capblue,->,line width=.85pt]
      (90:2.37)
      arc[start angle=90,end angle={90+\sgn*\CAP},radius=2.37];
}

\draw[removed]
    ({90-\RMAX}:\RR)
    arc[start angle={90-\RMAX},end angle={90+\RMAX},radius=\RR];
\draw[removed]
    ({90-\RMAX*\EE}:{\RR*\EE})
    arc[start angle={90-\RMAX*\EE},
        end angle={90+\RMAX*\EE},radius={\RR*\EE}];

\draw[black!50,thin] (0,0)--(90:2.38);
\node[capblue,fill=white,inner sep=.7pt]
    at (72:2.66) {$b_\varepsilon$};
\node[capblue,fill=white,inner sep=.7pt]
    at (108:2.66) {$b_\varepsilon$};
\node[font=\scriptsize,anchor=south] at (90:2.66) {$\theta_j$};
\node[capblue,font=\scriptsize,anchor=west]
    at (58:2.28) {$s=b_\varepsilon$};
\node[capblue,font=\scriptsize,anchor=east]
    at (122:2.28) {$s=b_\varepsilon$};

\node[borderred,font=\scriptsize,anchor=east]
    at (-.90,.80) {$K_{j,\varepsilon,p}$};
\draw[borderred,->,thin] (-.82,.85)--(-.16,1.57);
\node[font=\scriptsize] at (.21,-.32) {$\mathbb{B}_\varepsilon$};
\node[font=\scriptsize,anchor=west] at (.62,-1.52) {$r=\varepsilon$};
\draw[->,thin] (.67,-1.30)--(.45,-.56);
\fill (0,0) circle (1.6pt);
\node[font=\scriptsize,below left] at (0,0) {$0$};
\end{scope}

\begin{scope}[shift={(6.35,0)}]
\node[font=\bfseries\small]
    at (1.80,3.25) {(b) Magnified sector};
\begin{scope}[scale=.70]
\def\RR{4.80}
\def\CUT{42}

\fill[keep]
    (0,0)--(-\CUT:\RR)
    arc[start angle=-\CUT,end angle=\CUT,radius=\RR]--cycle;
\fill[trans]
    (0,0)--(-\CAP:\RR)
    arc[start angle=-\CAP,end angle=\CAP,radius=\RR]--cycle;

\path[fill=zero]
    (0,0)
    -- plot[domain=0:1,samples=120,variable=\t]
      ({\RR*\t*cos(-\RMAX*\t)},
       {\RR*\t*sin(-\RMAX*\t)})
    arc[start angle=-\RMAX,end angle=\RMAX,radius=\RR]
    -- plot[domain=1:0,samples=120,variable=\t]
      ({\RR*\t*cos(\RMAX*\t)},
       {\RR*\t*sin(\RMAX*\t)})
    --cycle;

\path[fill=zerocore]
    (0,0)
    -- plot[domain=0:\EE,samples=80,variable=\t]
      ({\RR*\t*cos(-\RMAX*\t)},
       {\RR*\t*sin(-\RMAX*\t)})
    arc[start angle={-\RMAX*\EE},
        end angle={\RMAX*\EE},radius={\RR*\EE}]
    -- plot[domain=\EE:0,samples=80,variable=\t]
      ({\RR*\t*cos(\RMAX*\t)},
       {\RR*\t*sin(\RMAX*\t)})
    --cycle;

\draw[black!25,thin] (-\CUT:\RR)--(0,0)--(\CUT:\RR);
\draw[line width=.95pt]
    (-\CUT:\RR) arc[start angle=-\CUT,end angle=-\RMAX,radius=\RR];
\draw[line width=.95pt]
    (\RMAX:\RR) arc[start angle=\RMAX,end angle=\CUT,radius=\RR];
\draw[black!50,thin]
    (-\CUT:{\RR*\EE})
    arc[start angle=-\CUT,end angle={-\RMAX*\EE},
        radius={\RR*\EE}];
\draw[black!50,thin]
    ({\RMAX*\EE}:{\RR*\EE})
    arc[start angle={\RMAX*\EE},end angle=\CUT,
        radius={\RR*\EE}];

\foreach \sgn in {-1,1} {
    \draw[borderred,line width=.8pt]
      plot[domain=0:\EE,samples=80,variable=\t]
        ({\RR*\t*cos(\sgn*\RMAX*\t)},
         {\RR*\t*sin(\sgn*\RMAX*\t)});
    \draw[removed]
      plot[domain=\EE:1,samples=120,variable=\t]
        ({\RR*\t*cos(\sgn*\RMAX*\t)},
         {\RR*\t*sin(\sgn*\RMAX*\t)});
    \draw[capblue,line width=.7pt]
      (0,0)--({\sgn*\CAP}:{\RR+.15});
    \node[capblue,font=\scriptsize,anchor=west]
      at ({\sgn*\CAP}:{\RR+.20}) {$s=b_\varepsilon$};
}

\draw[removed]
    (-\RMAX:\RR)
    arc[start angle=-\RMAX,end angle=\RMAX,radius=\RR];
\draw[removed]
    ({-\RMAX*\EE}:{\RR*\EE})
    arc[start angle={-\RMAX*\EE},
        end angle={\RMAX*\EE},radius={\RR*\EE}];

\draw[black!40,->,thin]
    (0,0)--({\RR+.35},0) node[right,text=black] {$\theta_j$};
\fill (0,0) circle (1.6pt);
\node[font=\scriptsize,below left] at (0,0) {$0$};
\node[font=\scriptsize,anchor=west] at (7:{\RR+.16}) {$r=1$};

\node[font=\scriptsize] at (36.5:4.07) {$w=1$};
\node[font=\scriptsize] at (-36.5:4.07) {$w=1$};
\node[font=\scriptsize,rotate=26] at (26:3.85) {$0<w<1$};
\node[font=\scriptsize,rotate=-26] at (-26:3.85) {$0<w<1$};
\node[font=\scriptsize,fill=zero,inner sep=2pt] at (3.95,0) {$w=0$};

\node[font=\scriptsize,anchor=south] at (1.70,2.30)
    {$s=\rho_{j,p}(r)$};
\draw[borderred,->,thin]
    (2.05,2.10)--({\RR*.62*cos(\RMAX*.62)},
                   {\RR*.62*sin(\RMAX*.62)});
\node[font=\scriptsize,anchor=south] at (.60,.95) {$w=0$};
\draw[borderred,->,thin] (.72,.89)--(1.24,0);
\node[font=\scriptsize,anchor=north] at (1.25,-1.45) {$r=\varepsilon$};
\draw[->,thin] (1.25,-1.36)--(-38:{\RR*\EE});

\end{scope}
\end{scope}

\draw[zerocore,line width=6pt] (2.25,-2.85)--(2.67,-2.85);
\node[font=\scriptsize,anchor=west] at (3.87,-2.85)
    {$w=0$ in $\mathbb{B}_\varepsilon \cap \{(r,\theta):s(\theta)\le\rho_{j,p}(r) \}$ (kept in $G_{\varepsilon,p}$)};
\draw[zero,line width=6pt] (2.25,-3.25)--(2.67,-3.25);
\node[font=\scriptsize,anchor=west] at (3.87,-3.25)
    {$w=0$ for $r\ge\varepsilon$ (removed part $K_{j,\varepsilon,p}$)};
\end{tikzpicture}%
}
\end{center}

For simplicity, write $w=w_{\eps,p}$. Then $0\le w\le1$, and it follows from
\eqref{eq:capformula}-\eqref{eq:kappacut} and \eqref{eq:Abounds}-\eqref{eq:cells} that
\[
\|\nabla_{\Sph^{N-1}}w(r,\cdot)\|_{L^2(\Sph^{N-1})}^2
=\sum_j \alpha_jA_p(r,\theta_j)\le Cr^2.
\]
Since $1-w$ has zero trace on each cap boundary $ \partial B_{\Sph^{N-1}}(\theta_j,b_{\eps})$, applying the rescaled Poincar\'e's
inequality in cap $B_{\Sph^{N-1}}(\theta_j,b_{\eps}) $ yields 
\begin{align*}
    \|1-w(r,\cdot)\|_{L^2(\Sph^{N-1})}^2&= \sum_j\|1-w(r,\cdot)\|_{L^2(B_{\Sph^{N-1}}(\theta_j,b_{\eps}))}^2\\
    &\leq C\varepsilon^2\|\nabla_{\Sph^{N-1}}w(r,\cdot)\|_{L^2(\Sph^{N-1})}^2\le C\eps^2r^2.
\end{align*}

On the  annuli $\mathscr{A}_{j,p}=\{\theta\in \Sph^{N-1}: \rho_{j,p}(r)<\dist_{\Sph^{N-1}}(\theta,\theta_j)<b_{\eps}\}$, we see from \eqref{eq:wdef} and \eqref{eq:kappacut}
that \[
\partial_rw=-\frac{A'_p(r,\theta_j)}{A_p(r,\theta_j)}(1-w),
\]
and elsewhere this derivative is zero almost everywhere.
 Hence, \eqref{eq:Abounds} implies $w\in H^1_{\mathrm{loc}}(\mathbb{B}\backslash\{0\})$ and
\begin{align}\label{eq:correctorslices}
&\|\nabla_{\Sph^{N-1}}w(r,\cdot)\|_{L^2(\Sph^{N-1})}^2\le Cr^2,\quad
\|1-w(r,\cdot)\|_{L^2(\Sph^{N-1})}^2\le C\eps^2r^2,\notag\\
&\|\partial_rw(r,\cdot)\|_{L^2(\Sph^{N-1})}^2\le C\eps^2.
\end{align}
Using $dz=r^{N-1}\,dr\,d\sigma$ and
$|\nabla w|^2=|\partial_rw|^2+r^{-2}|\nabla_{\Sph^{N-1}}w|^2$, we obtain
\begin{equation}\label{eq:correctorbounds}
\|1-w\|_{L^2(\mathbb{B})}+\|\partial_rw\|_{L^2(\mathbb{B})}\le C\eps,
\qquad \|\nabla w\|_{L^2(\mathbb{B})}\le C,
\end{equation}
and, for $0<a\le1$,
\begin{align}
\label{eq:coreenergy}
\int_{B_a}|\nabla w|^2\dd
\le&\int_0^ar^{N-1}\|\partial_rw(r,\cdot)\|_{L^2(\Sph^{N-1})}^2+r^{N-3}\|\nabla_{\Sph^{N-1}}w(r,\cdot)\|_{L^2(\Sph^{N-1})}^2\,dr\nonumber \\
\le&C\int_0^a(r^{N-1}\eps^2+r^{N-3}r^2)\,dr
\le Ca^{N}.  
\end{align}
These estimates imply $w\in H^1(\mathbb{B})$. 
Then for any $v\in C_c^\infty(\mathbb{B})$, it follows from $w=0$ on $\cup_jK_{j,\varepsilon,p}$ that
\begin{equation}\label{eq:testadmissible}
w_{\eps,p}v\in H_0^1(G_{\eps,p});
\end{equation}
see \cite[Chapter~3]{McLean}.

\medskip
We next identify the angular flux generated by the corrector and its
contribution to the limiting equation.

\begin{lemma}\label{lem:angularflux}
If $g\in H^1(\Sph^{N-1})$ and g vanishes in the inner boundary, i.e,
$g|_{s(\theta)=\rho_{j,p}(r)}=0$ for all $j$, then
\begin{equation}\label{eq:fluxidentity}
\int_{\Sph^{N-1}}\nabla_{\Sph^{N-1}}w_{\eps,p}\cdot\nabla_{\Sph^{N-1}}g\ds
=\sum_j \alpha_jA_p(r,\theta_j)
\avint_{\partial B_{\Sph^{N-1}}(\theta_j,b_\eps)}g\,dS.
\end{equation}
For arbitrary $g\in H^1(\Sph^{N-1})$,
\begin{equation}\label{eq:fluxerror}
\left|
\sum_j \alpha_jA_p(r,\theta_j)
\avint_{\partial B_{\Sph^{N-1}}(\theta_j,b_\eps)}g\,dS
-\int_{\Sph^{N-1}}A_p(r,\theta)g(\theta)\ds
\right|
\le C\eps r^2\|g\|_{H^1(\Sph^{N-1})}.
\end{equation}
\end{lemma}

\begin{proof}
Note that $w=w_{\eps,p}$ is harmonic on $\rho_{j,p}<s<b_{\eps}$ and $g|_{s(\theta)=\rho_{j,p}(r)}=0$. Integration by parts, together with \eqref{eq:kappacut} and the fact that $s$ is precisely the geodesic distance from $\theta_j$, gives
\begin{align*}
&\int_{\rho_{j,p}(r)<s<b_\eps}
\nabla_{\Sph^{N-1}}w\cdot\nabla_{\Sph^{N-1}}g\ds\\
=&\int_{s=b_{\eps}}
g\partial_sw\, dS=\frac{\alpha_jA_p(r,\theta_j)}{|\mathbb S^{N-2}|(\sin b_{\varepsilon})^{N-2}}\int_{s=b_{\eps}}
g\, dS\\
=&\alpha_jA_p(r,\theta_j)\avint_{\partial B_{\Sph^{N-1}}(\theta_j,b_\eps)}g\,dS.
\end{align*}
Summation proves \eqref{eq:fluxidentity}; the trace theorem extends
it to $H^1$ functions.

Since $\alpha_jA_p(r,\theta_j)\le C\eps^{N-1}r^2$ and there are $O(\eps^{1-N})$ cells,
Lemma~\ref{lem:celltrace} gives
\begin{align}\label{Sec 4 Lem-2 equ-a-0}
&\sum_j \alpha_jA_p(r,\theta_j)
\left|
\avint_{\partial B_{\Sph^{N-1}}(\theta_j,b_\eps)}g\,dS
-\avint_{Q_{j,\eps}}g\ds
\right|\nonumber\\
&\qquad\le C\eps^{1+(N-1)/2}r^2
\sum_j\|\nabla_{\Sph^{N-1}}g\|_{L^2(Q_{j,\eps})}\nonumber\\
&\qquad\le C\eps^{1+(N-1)/2}r^2
\left(\sum_j1\right)^{1/2}\left(\sum_j\|\nabla_{\Sph^{N-1}}g\|^2_{L^2(Q_{j,\eps})}\right)^{1/2}\nonumber\\
&\qquad\le C\eps r^2\|\nabla_{\Sph^{N-1}}g\|_{L^2(\Sph^{N-1})}.
\end{align}
Using \eqref{eq:Abounds}, we can get $|A_p(r,\theta_j)-A_p(r,\theta)|\leq C\eps r^3$ in $Q_{j,\eps}$, hence it follows that
\begin{align}\label{Sec 4 Lem-2 equ-a}
&\left|\sum_j \alpha_jA_p(r,\theta_j)\avint_{Q_{j,\eps}}g(\theta)\ds-\int_{\Sph^{N-1}}A_p(r,\theta)g(\theta)\ds\right|\notag\\
\leq&
\sum_j\int_{Q_{j,\eps}}
|A_p(r,\theta_j)-A_p(r,\theta)|\,|g(\theta)|\ds
\le C\eps r^3\|g\|_{L^2(\Sph^{N-1})}.
\end{align}
 Using \eqref{Sec 4 Lem-2 equ-a-0} and \eqref{Sec 4 Lem-2 equ-a}, we finish the proof of  \eqref{eq:fluxerror}.
\end{proof}

\begin{proposition}\label{prop:corrector}
For every $g\in H_0^1(G_{\eps,p})$, extended by zero to $\mathbb{B}$,
\begin{equation}\label{eq:corrector}
\left|
\int_{\mathbb{B}}\nabla w_{\eps,p}\cdot\nabla g\dd-\int_{\mathbb{B}}W_pg\dd
\right|
\le C\eps\|g\|_{H^1(\mathbb{B})}.
\end{equation}
\end{proposition}

\begin{proof}
First, take $g\in C_c^\infty(G_{\eps,p})$, extend it by zero, and write $w=w_{\eps,p}$. Using $$\nabla w\cdot\nabla g=\partial_rw\,\partial_rg+r^{-2}\nabla_{\Sph^{N-1}}w\cdot\nabla_{\Sph^{N-1}}g,$$ we decompose the integral into three parts:
\begin{align}\label{Lem 3 equ-b}
    &\int_{\mathbb{B}}\nabla w\cdot\nabla g\dd-W_pg\dd\nonumber\\
    =&\int_{\mathbb{B}\setminus \mathbb{B}_\eps}
 r^{-2}\nabla_{\Sph^{N-1}}w\cdot\nabla_{\Sph^{N-1}}g-W_pg\dd
 +\int_{\mathbb{B}\setminus \mathbb{B}_\eps}
 \partial_rw\partial_rg\dd\nonumber\\
 &+ \int_{B_{\eps}}\nabla w\cdot\nabla g\dd-W_pg\dd.
\end{align}
For $r>\eps$, since $g(r,\cdot)$ vanishes on the boundary of $G_{\eps,p}$, we have $g|_{s(\theta)=\rho_{j,p}(r)}=0$, and hence Lemma \ref{lem:angularflux} applies. By \eqref{eq:fluxidentity}--\eqref{eq:fluxerror} and
$A_p=r^2W_p$,
{\allowdisplaybreaks
\begin{align*}
&\left|
\int_{\mathbb{B}\setminus \mathbb{B}_\eps}
 r^{-2}\nabla_{\Sph^{N-1}}w\cdot\nabla_{\Sph^{N-1}}g\dd
-W_pg\dd
\right|\\
=&\left|\int_{\eps}^{1}r^{N-3}
\int_{\Sph^{N-1}}
 \nabla_{\Sph^{N-1}}w\cdot\nabla_{\Sph^{N-1}}g-A_pg\ds
dr
\right|\\
\le& C\eps\int_\eps^1r^{N-1}
\bigl(\|g(r,\cdot)\|_{L^2(\Sph^{N-1})}+\|\nabla_{\Sph^{N-1}}g(r,\cdot)\|_{L^2(\Sph^{N-1})}\bigr)\,dr\\
\le & C\eps\left(\int_\eps^1r^{N-1} dr\right)^{1/2}\left(\int_\eps^1r^{N-1}
\|g(r,\cdot)\|^2_{L^2(\Sph^{N-1})}dr\right)^{1/2}\\
&+C\eps\left(\int_\eps^1r^{N+1} dr\right)^{1/2}\left(\int_\eps^1r^{N-3}
\|\nabla_{\Sph^{N-1}}g(r,\cdot)\|^2_{L^2(\Sph^{N-1})}dr\right)^{1/2}\\
\le& C\eps\|g\|_{H^1(\mathbb{B})}.
\end{align*}
}%
 The second term  in \eqref{Lem 3 equ-b} is at most
$C\eps\|\nabla g\|_{L^2(\mathbb{B})}$ by \eqref{eq:correctorbounds}.
On $\mathbb{B}_\eps$, we see from $\|W_p\|_{L^{\infty}(\mathbb{B})}\leq C$ and \eqref{eq:coreenergy} that
\[
\left|\int_{\mathbb{B}_\eps}\nabla w\cdot\nabla g\dd\right|
+\left|\int_{\mathbb{B}_\eps}W_pg\dd\right|
\le C\eps^{N/2}\|g\|_{H^1(\mathbb{B})}
\le C\eps\|g\|_{H^1(\mathbb{B})}.
\]
This proves \eqref{eq:corrector} for smooth tests. Density in
$H_0^1(G_{\eps,p})$ gives the general case.
\end{proof}

\subsection{Identification of the limit equation}\label{sec:limitequation}
Recall that \(E_{G_{\eps,p}}\) is the extension operator from \(L^2(G_{\eps,p})\) to \(L^2(\mathbb{B})\) and \(E^*_{G_{\eps,p}}\) is simply the restriction operator, namely,
$
E^*_{G_{\eps,p}}f=f|_{G_{\eps,p}}.
$ The following proposition is the key point of this section, and the previous auxiliary function $w_{\eps,p}(r,\theta)$ plays a crucial role in its proof.
\begin{proposition}\label{prop:roughnorm}
The domains $G_{\eps,p}$ satisfy
\begin{equation}\label{eq:roughnorm}
\sup_{p\in K_r}
\left\|
E_{G_{\eps,p}}(-\Delta_{G_{\eps,p}}^D+1)^{-1}E_{G_{\eps,p}}^*
-(-\Delta_{\mathbb{B}}^D+W_p+1)^{-1}
\right\|\longrightarrow0
\quad(\eps\downarrow0).
\end{equation}
\end{proposition}

\begin{proof}
Suppose the conclusion fails. Then there are $\eps_n\to0$, $p_n\in K_r$,
and $f_n\in L^2(\mathbb{B})$ with $\|f_n\|_{L^2(\mathbb{B})}=1$ such that the weak solutions
\[
\begin{cases}
-\Delta u_n+u_n=f_n&\text{in }G_{\eps_n,p_n},\\
u_n\in H_0^1(G_{\eps_n,p_n}),
\end{cases}
\qquad
\begin{cases}
-\Delta v_n+(W_{p_n}+1)v_n=f_n&\text{in }\mathbb{B},\\
v_n\in H_0^1(\mathbb{B}),
\end{cases}
\]
satisfy $\|u_n-v_n\|_{L^2(\mathbb{B})}\ge c>0$, where $u_n$ is extended by zero to $\mathbb{B}$.
After taking a subsequence, $p_n\to p$ and $f_n\rightharpoonup f$ in
$L^2(\mathbb{B})$. Energy estimates give
\[
\|u_n\|_{H^1(\mathbb{B})}+\|v_n\|_{H^1(\mathbb{B})}\le C.
\]
After a further subsequence, $u_n\rightharpoonup u$ and
$v_n\rightharpoonup v$ in $H_0^1(\mathbb{B})$, with both convergences strongly
in $L^2(\mathbb{B})$. Passing to the limit in the equation for $v_n$ gives
$-\Delta v+(W_p+1)v=f$.

Write $w_n=w_{\eps_n,p_n}$ and $W_n=W_{p_n}$. The weak equation for $u_n$ is
\begin{equation}\label{eq:roughweak}
\int_{\mathbb{B}}(\nabla u_n\cdot\nabla\varphi+u_n\varphi)\dd
=\int_{\mathbb{B}}f_n\varphi\dd,
\qquad \varphi\in H_0^1(G_{\eps_n,p_n}).
\end{equation}
Fix $\eta\in C_c^\infty(\mathbb{B})$.  By Proposition \ref{prop:corrector} and using $u_n\eta$ as an
admissible test in \eqref{eq:corrector}, we get
\begin{equation}\label{eq:correctoru}
\begin{aligned}
\left|\int_{\mathbb{B}}\nabla (\eta u_n)\cdot\nabla w_n\dd
-\int_{\mathbb{B}}W_nu_n\eta\dd\right|
\le C_\eta\eps_n\|u_n\|_{H^1(\mathbb{B})}.
\end{aligned}
\end{equation}
Integration by parts gives
\begin{align*}
\int_{\mathbb{B}}u_n\nabla\eta\cdot\nabla w_n\dd
=-\int_{\mathbb{B}}(w_n-1)
(\nabla u_n\cdot\nabla\eta+u_n\Delta\eta)\dd.
\end{align*}
From here and \eqref{eq:correctorbounds}, we get 
\begin{align}\label{eq:crosssmall}
    \left|\int_{\mathbb{B}}u_n\nabla\eta\cdot\nabla w_n\dd\right|\le C_\eta\eps_n\|u_n\|_{H^1(\mathbb{B})}.
\end{align}
Applying \eqref{eq:correctoru} and \eqref{eq:crosssmall}, we derive that 
\begin{align}\label{Sec 4 Lem 3 equ-a}
    \int_{\mathbb{B}}\eta\nabla u_n\cdot\nabla w_n\dd=\int_{\mathbb{B}} W_nu_n\eta\dd+O(\eps_n).
\end{align}

On the other hand,
by \eqref{eq:testadmissible}, $w_n\eta$
is an admissible test in \eqref{eq:roughweak}, which gives
\begin{align}\label{Sec 4 Lem 3 equ-b}
    \int_{\mathbb{B}}\nabla u_n\cdot\nabla(w_n\eta)+u_nw_n\eta\dd
=\int_{\mathbb{B}}f_nw_n\eta\dd.
\end{align}
Plugging \eqref{Sec 4 Lem 3 equ-a} into \eqref{Sec 4 Lem 3 equ-b}, we derive
\begin{align*}
     \int_{\mathbb{B}}\nabla u_n\cdot\nabla\eta w_n+(W_n+w_n)u_n\eta\dd
=\int_{\mathbb{B}}f_nw_n\eta\dd+O(\eps_n).
\end{align*}
Using $\|w_n-1\|_{L^2(\mathbb{B})}\leq C\eps_n$ by \eqref{eq:correctorbounds}, we obtain
\begin{equation}\label{eq:limitresidual}
\begin{aligned}
\int_{\mathbb{B}}
\bigl(\nabla u_n\cdot\nabla\eta+(W_n+1)u_n\eta-f_n\eta\bigr)\dd=O(\eps_n).
\end{aligned}
\end{equation}
Since $W_n\to W_p$ uniformly, the limit satisfies
\[
\int_{\mathbb{B}}\bigl(\nabla u\cdot\nabla\eta+(W_p+1)u\eta\bigr)\dd
=\int_{\mathbb{B}}f\eta\dd,\quad\forall\eta\in C_c^\infty(\mathbb B).
\]
Density and uniqueness imply $u=v$, contradicting
$\|u_n-v_n\|_{L^2(\mathbb{B})}\ge c$. This proves \eqref{eq:roughnorm}.
\end{proof}

\subsection{Smooth strictly star-shaped approximations}\label{sec:boundarysmoothing} We begin with a perturbation estimate of domains.
For a positive continuous function $R$ on $\Sph^{N-1}$, write
$\Omega_R=\{r\theta:0\le r<R(\theta)\}$.
All Dirichlet solutions in this subsection are extended by zero
to $\R^N$.

\begin{lemma}\label{lem:radialcontinuity}
Suppose $a\le R,R'\le L$, with $a>0$, and
$\delta=\|R-R'\|_{L^{\infty}(\Sph^{N-1})}/a<1/2$.
For $h\in L^2(\R^N)$, let $u_R\in H_0^1(\Omega_R)$ and
$u_{R'}\in H_0^1(\Omega_{R'})$ solve
$-\Delta u+u=h$ on their respective domains. Define
\[
T_R:L^2(\mathbb R^N)\longrightarrow L^2(\mathbb R^N),
\qquad T_Rh:=u_R,
\]
\[
T_{R'}:L^2(\mathbb R^N)\longrightarrow L^2(\mathbb R^N),
\qquad T_{R'}h:=u_{R'},
\]
where the solutions are extended by zero.
Then
\begin{equation}\label{eq:radialresolventbound}
\|T_R-T_{R'}\|
\le \frac{C_{N,L}}{a}\|R-R'\|_{L^{\infty}(\Sph^{N-1})}.
\end{equation}
\end{lemma}

\begin{proof}
All functions below are extended by zero to $\R^N$. Set
\[
F(w):=2\int_{\R^N}hw\,dx-\|w\|_{H^1(\R^N)}^2,
\quad
\|w\|_{H^1(\R^N)}^2
=\int_{\R^N}\bigl(|\nabla w|^2+|w|^2\bigr)\,dx.
\]
For an open set $\Omega\subset \mathbb B_L$, let
$u_\Omega\in H_0^1(\Omega)$ solve $-\Delta u_\Omega+u_\Omega=h$.
Testing with $u_\Omega$ gives
\begin{align}\label{Sec 5 energy bound}
F(u_\Omega)
=\int_{\R^N}hu_\Omega\,dx
=\sup_{w\in H_0^1(\Omega)}F(w),
\qquad
\|u_\Omega\|_{H^1(\R^N)}
\le\|h\|_{L^2(\R^N)}.
\end{align}
For $w\in H_0^1(\Omega)$ and $1/2<t<2$, set
$
w_t(x):=t^{-N/2}w(x/t)\in H_0^1(t\Omega).
$
A change of variables yields
\[
\|w_t\|_{L^2(\R^N)}=\|w\|_{L^2(\R^N)},
\qquad
\|\nabla w_t\|_{L^2(\R^N)}
=t^{-1}\|\nabla w\|_{L^2(\R^N)}.
\]
For $w\in C_c^\infty(\Omega)$, differentiation gives
\[
\partial_t w_t(x)
=-\frac1t\left(\frac N2w_t(x)+x\cdot\nabla w_t(x)\right).
\]
Since $w_t$ is supported in $\mathbb B_{tL}$,
\[
\|\partial_t w_t\|_{L^2(\R^N)}
\le
\frac1t\left(
\frac N2\|w\|_{L^2(\R^N)}
+L\|\nabla w\|_{L^2(\R^N)}
\right)
\le C_{N,L}\|w\|_{H^1(\R^N)}.
\]
Integration in $t$ and density therefore imply
\begin{align}\label{Sec 5 equ-b}
\|w_t-w\|_{L^2(\R^N)}
\le C_{N,L}|t-1|\|w\|_{H^1(\R^N)}
\qquad\mathrm{for}\qquad w\in H_0^1(\Omega).
\end{align}
Using $(u_\Omega)_t\in H_0^1(t\Omega)$ as a trial function,
we obtain
\[
\begin{aligned}
\int_{\R^N}h(u_{t\Omega}-u_\Omega)\,dx
&=F(u_{t\Omega})-F(u_\Omega)
\ge F((u_\Omega)_t)-F(u_\Omega)\\
&=
2\int_{\R^N}h\bigl((u_\Omega)_t-u_\Omega\bigr)\,dx
-(t^{-2}-1)\|\nabla u_\Omega\|_{L^2(\R^N)}^2\\
&\ge -C_{N,L}|t-1|\|h\|_{L^2(\R^N)}^2,
\end{aligned}
\]
where the last step follows from \eqref{Sec 5 energy bound}, \eqref{Sec 5 equ-b} and $|t^{-2}-1|\le C|t-1|$.

For the reverse bound, we use $(u_{t\Omega})_{1/t}\in H_0^1(\Omega)$.
Since $t\Omega\subset \mathbb B_{2L}$ and
$|t^{-1}-1|\le2|t-1|$, the same dilation estimate gives
\begin{align}\label{Sec 5 equ-c}
\|(u_{t\Omega})_{1/t}-u_{t\Omega}\|_{L^2(\R^N)}
\le& C_{N,2L}|t^{-1}-1|\|u_{t\Omega}\|_{H^1(\R^N)}\nonumber\\
\le& C_{N,L}|t-1|\|h\|_{L^2(\R^N)}.
\end{align}
Using $(u_{t\Omega})_{1/t}$ in the variational formula on
$\Omega$, we consequently obtain
\[
\begin{aligned}
\int_{\R^N}h(u_{t\Omega}-u_\Omega)\,dx
&=F(u_{t\Omega})-F(u_\Omega)\le F(u_{t\Omega})-F((u_{t\Omega})_{1/t})\\
&=
2\int_{\R^N}h\bigl(u_{t\Omega}-(u_{t\Omega})_{1/t}\bigr)\,dx
+(t^2-1)\|\nabla u_{t\Omega}\|_{L^2(\R^N)}^2\\
&\le C_{N,L}|t-1|\|h\|_{L^2(\R^N)}^2.
\end{aligned}
\]
Combining the two bounds yields
\begin{equation}\label{eq:dilationresolvent}
\left|\int_{\R^N}h(u_{t\Omega}-u_\Omega)\,dx\right|
\le C_{N,L}|t-1|\|h\|_{L^2(\R^N)}^2.
\end{equation}

Now, since $|R-R'|\le a\delta\le\delta R$, we have
\[
(1-\delta)\Omega_R
\subset\Omega_{R'}
\subset(1+\delta)\Omega_R.
\]
If $\Omega_1\subset\Omega_2$, then 
$H_0^1(\Omega_1)\subset H_0^1(\Omega_2)$ after zero extension, we can get $F(u_{\Omega_1})\leq F(u_{\Omega_2})$.
Consequently,
\[
\begin{aligned}
\left|\int_{\R^N}h(u_R-u_{R'})\,dx\right|
&\le
\int_{\R^N}h
\bigl(u_{(1+\delta)\Omega_R}
-u_{(1-\delta)\Omega_R}\bigr)\,dx\\
&\le C_{N,L}\delta\|h\|_{L^2(\R^N)}^2,
\end{aligned}
\]
where the last inequality follows by applying
\eqref{eq:dilationresolvent} separately with $t=1+\delta$
and $t=1-\delta$. 
Consequently,
\[
|\langle (T_R-T_{R'})h,h\rangle_{L^2}|
\leq C_{N,L}\delta\|h\|_{L^2(\mathbb R^N)}^2,\quad\forall~ h\in L^2(\mathbb R^N).
\]
Since both \(T_R\) and \(T_{R'}\) are bounded and self-adjoint,
the norm formula for bounded self-adjoint operators therefore yields
\[
\|T_R-T_{R'}\|
=\sup_{\|h\|_{L^2}=1}|\langle (T_R-T_{R'})h,h\rangle_{L^2}|
\leq C_{N,L}\delta.
\]
The proof is complete.
\end{proof}

Recall the radial function $R_{\eps,p}\in C(\mathbb S^{N-1})$ defined in \eqref{eq:roughR}.
For fixed $\eps$, smooth the function $R_{\eps,p}$ by
\[
R_{\eps,p}^{\,t}:=e^{t\Delta_{\Sph^{N-1}}}R_{\eps,p}\in C^{\infty}(\mathbb S^{N-1}),\qquad t>0.
\]
More precisely, $R_{\eps,p}^{\,t} $ is the solution of the heat equation
\begin{align*}
    \begin{cases}
        \displaystyle \partial_tR_{\eps,p}^{\,t}=\Delta_{\Sph^{N-1}} R_{\eps,p}^{\,t}\qquad&\mathrm{in}\qquad (0,+\infty)\times \Sph^{N-1}\\
          \displaystyle R_{\eps,p}^{0}=R_{\eps,p}\qquad&\mathrm{on}\qquad t=0.
    \end{cases}
\end{align*}
Since $\eps\le R_{\eps,p}\le1 $, and $R_{\eps,p}^{\,t} $ can be written as $R_{\eps,p}^{\,t}=K^{\,t}\ast R_{\eps,p}$, where $K^{\,t} $ is the heat kernel in $(0,+\infty)\times \Sph^{N-1}$ with $\int_{\Sph^{N-1} }K^t=1$,
then
$$\eps\le R_{\eps,p}^{\,t}\le1.$$ 

For fixed $t>0$, the smoothed function $R_{\eps,p}^{\,t}$
depends continuously on $p$ in every $C^k$ norm. Since
$\{R_{\eps,p}:p\in K_r\}$ is compact in $C^0$, strong continuity and
the contraction property give
\[
\sup_{p\in K_r}
\|R_{\eps,p}^{\,t}-R_{\eps,p}\|_{L^{\infty}(\Sph^{N-1})}\longrightarrow0
\qquad\mathrm{as}\qquad t\to 0.
\]
Choose $t_\eps>0$, uniformly in $p$, so that $\sup_{p\in K_r} \|R_{\eps,p}^{\,t_\eps}-R_{\eps,p}\|_{L^{\infty}(\Sph^{N-1})}\le \eps^2$, and let    
\begin{align*}
    \Omega_{\eps,p}:=\left\{r\theta: 0\leq r< R_{\eps,p}^{\,t_\eps}(\theta)\right\}.
\end{align*}
By
\eqref{eq:normal} and $R_{\eps,p}^{\,t_{\eps}}\geq \eps $, we get
\begin{align*}
    &x\cdot\nu_{\eps,p}^{\,t_{\eps}}=
 \frac{R_{\eps,p}^{\,t_{\eps}}(\theta)^2}
 {\sqrt{R_{\eps,p}^{\,t_{\eps}}(\theta)^2+|\nabla_{\Sph^{N-1}}R_{\eps,p}^{\,t_{\eps}}(\theta)|^2}}>\tilde c_{\eps}>0,
\end{align*}
hence $\Omega_{\eps,p}$ is smooth and strictly star-shaped.

\begin{proof}[\bf Proof of Theorem~\ref{thm:radial}]
Write
\[
\mathcal R_\Omega
:=E_\Omega(-\Delta_\Omega^D+1)^{-1}E_\Omega^*,
\qquad
T_p:=(-\Delta_{\mathbb{B}}^D+W_p+1)^{-1},
\]
which are compact nonnegative self-adjoint operators on $L^2(\mathbb{B})$.

Since $\eps\le R_{\eps,p},R_{\eps,p}^{\,t_\eps}\le1$ and
$\sup_{p\in K_r}\|R_{\eps,p}^{\,t_\eps}-R_{\eps,p}\|_{L^{\infty}(\Sph^{N-1})}\leq\eps^2$,
Lemma \ref{lem:radialcontinuity} gives
\begin{equation}\label{eq:smoothnorm}
\sup_{p\in K_r}
\|\mathcal R_{\Omega_{\eps,p}}-\mathcal R_{G_{\eps,p}}\|
\le C\eps.
\end{equation}
Choose $\eps_k\downarrow0$ along the cell construction and set
$\Omega_{k,p}:=\Omega_{\eps_k,p}$.
Together with Proposition~\ref{prop:roughnorm}, this yields
\[
\sup_{p\in K_r}\|\mathcal R_{\Omega_{k,p}}-T_p\|
\le C\eps_k+
\sup_{p\in K_r}\|\mathcal R_{G_{\eps_k,p}}-T_p\|
\longrightarrow0,
\]
proving \eqref{eq:radialnorm}.

We next verify the continuity in $p$. Fix $\eps>0$ and let $p_n\to p$.
By \eqref{eq:radialparameter},
$\|R_{\eps,p_n}-R_{\eps,p}\|_{L^{\infty}(\Sph^{N-1})}\to0$.
Since $t_\eps$ is independent of $p$, the heat semigroup estimate gives,
for every integer $\ell\ge0$,
\[
\|R_{\eps,p_n}^{\,t_\eps}-R_{\eps,p}^{\,t_\eps}\|_{C^\ell(\Sph^{N-1})}
\le C_{t_\eps,\ell}
\|R_{\eps,p_n}-R_{\eps,p}\|_{L^{\infty}(\Sph^{N-1})}
\longrightarrow0.
\]
Thus the radial function $R_{\eps_k,p}^{\,t_{\eps_k}}$ of $\Omega_{k,p}$ depends continuously
on $p$ in every $C^\ell$ norm.
Moreover, Lemma \ref{lem:radialcontinuity} and the  identity $T_{p_n}-T_p=T_{p_n}(W_p-W_{p_n})T_p$,
using $W_p\ge1$, give
\[
\begin{aligned}
\|\mathcal R_{\Omega_{\eps,p_n}}
  -\mathcal R_{\Omega_{\eps,p}}\|
&\le \frac{C}{\eps}
\|R_{\eps,p_n}^{\,t_\eps}-R_{\eps,p}^{\,t_\eps}\|_{L^{\infty}(\Sph^{N-1})}
\longrightarrow0,\\
\|T_{p_n}-T_p\|
&\le \|W_{p_n}-W_p\|_{L^{\infty}(\mathbb{B})}
\longrightarrow0.
\end{aligned}
\]
Here the first estimate holds for all sufficiently large $n$.
Hence both operator families are norm-continuous in $p$.

The remaining spectral assertions follow from
Remark~\ref{rem:counts} and Lemma~\ref{lem:transfer},
with $g(z)=z^{-1}-1$ for the operator matrices.
This completes the proof.\end{proof}

\section{Exact spectral matrices and prescribed multiplicities}
\label{sec:finish}

This section is devoted to the proof of the main result Theorem~\ref{thm:finite}.

\subsection{An exact matrix on one spectral subspace}
\label{sec:exactmatrix}

Recall the fixed orthonormal basis of $H=\mathbf 1^\perp$:
\[
h_k=\frac{e_1+\cdots+e_k-ke_{k+1}}{\sqrt{k(k+1)}},
\qquad 1\le k\le m,
\]
namely $\langle h_k,h_l\rangle_H=\delta_{kl}$.
Recall the functions $u_{k,\delta}(p)$ in Lemma~\ref{lem:smoothing}. By Proposition~\ref{prop:potential} and Remark \ref{rmk3-2}, it follows from the proof of Proposition~\ref{prop:potential} that $\sqrt{\hat{\rho}_{\tau,p,\delta}}u_{1,\delta}(p),\ldots, \sqrt{\hat{\rho}_{\tau,p,\delta}}u_{m,\delta}(p)$ is the continuous orthonormal frame for the subspace for eigenvalues $2,\ldots,m+1$.
Therefore, as in Remark \ref{rmk3-2}, we define
\[
J_p:H\longrightarrow L^2(\mathbb{B}),\qquad
J_px:=\sum_{k=1}^m
\langle x,h_k\rangle_H\,\sqrt{\hat{\rho}_{\tau,p,\delta}}u_{k,\delta}(p).
\]
In particular, $J_ph_k=\sqrt{\hat{\rho}_{\tau,p,\delta}}u_{k,\delta}(p)$.
Let $J_p^*:L^2(\mathbb{B})\to H$ be the adjoint of $J_p$, characterized by
\[
\langle J_px,f\rangle_{L^2(\mathbb{B})}
=\langle x,J_p^*f\rangle_H,
\qquad x\in H,\quad f\in L^2(\mathbb{B}).
\]
Explicitly,
\[
J_p^*f
=\sum_{k=1}^m
\left(\int_{\mathbb{B}} f\,\sqrt{\hat{\rho}_{\tau,p,\delta}}u_{k,\delta}(p)\dd\right)h_k.
\]
Then $J_p$ is an orthonormal frame for the spectral subspace
of $L_p=-\Delta_{\mathbb{B}}^D+V_p$ corresponding to
$\xi_2(p),\ldots,\xi_{m+1}(p)$, with
\[
J_p^*J_p=I_H,
\qquad
J_p^*L_pJ_p=\mathcal B_p,
\]
where $\mathcal B_p$ is the operator matrix on $H$ in Proposition~\ref{prop:potential}.
Moreover, $p\mapsto J_p$ is continuous in
$\mathcal L(H,L^2(\mathbb{B}))$.

\begin{proposition}\label{prop:exactmatrix}
Fix $N\ge3$ and $m\ge1$. There exist constants
$r\in(0,(m+1)/4)$, $\tau>0$, $C_*$, and a family of
smooth strictly star-shaped domains
$\Omega_p\subset\mathbb{B}$, $p\in K_r$, whose radial functions depend
continuously on $p$ in every $C^k$ norm, such that
\begin{equation}\label{eq:domaincluster}
\lambda_1(\Omega_p)<b_-<\lambda_2(\Omega_p)
\le\cdots\le\lambda_{m+1}(\Omega_p)
<b_+<\lambda_{m+2}(\Omega_p),
\end{equation}
where
\begin{equation}\label{eq:clusterthresholds}
b_-=C_*+\frac{3\tau(m+1)}8,
\qquad
b_+=C_*+\frac{7\tau(m+1)}4.
\end{equation}
The spectral subspace corresponding to eigenvalues $2,\ldots,m+1$
has an orthonormal frame
$\psi_{1,p},\ldots,\psi_{m,p}$, continuous in $L^2(\mathbb{B})$ after
extension by zero.
Relative to the fixed orthonormal basis $h_1,\ldots,h_m$ of $H$,
we define $\mathcal C_p\in\Sym(H)$ by
\[
\langle\mathcal C_ph_l,h_k\rangle_H
:=\int_{\Omega_p}
\nabla\psi_{k,p}\cdot\nabla\psi_{l,p}\dd,
\qquad 1\le k,l\le m,
\]
so that the eigenvalues of $\mathcal C_{p}$ are precisely $\lambda_2(\Omega_{p}),\ldots,\lambda_{m+1}(\Omega_{p})$.
Then for every $A\in\Sym(H)$ with $\|A\|\le r/4$, there exists $p_*\in K_r$ such that
\begin{equation}\label{eq:exacttarget}
\mathcal C_{p_*}=(C_*+\tau(m+1))I_H+\tau A.
\end{equation}
\end{proposition}

\begin{proof}
Let $V_p$ and $L_p=-\Delta_{\mathbb{B}}^D+V_p$ be given by
Proposition~\ref{prop:potential}, and let $J_p$ be 
defined above. Thus
\[
J_ph_k=\sqrt{\hat{\rho}_{\tau,p,\delta}}u_{k,\delta}(p),
\qquad
J_p^*J_p=I_H,
\qquad
J_p^*L_pJ_p=\mathcal B_p.
\]
Recall $C_*$ and $\Omega_{n,p}$ in Theorem~\ref{thm:radial}. 
Since $C_*$ is independent of $p$,
\begin{equation}\label{eq:BpLp}
J_p^*(L_p+C_*)J_p=\mathcal B_p+C_*I_H.
\end{equation}
Set
\[
T_p:=(L_p+C_*+1)^{-1},\quad
T_{n,p}:=E_{\Omega_{n,p}}
(-\Delta_{\Omega_{n,p}}^D+1)^{-1}E_{\Omega_{n,p}}^*,
\quad
\Pi_p:=J_pJ_p^*.
\]
By \eqref{eq:potentialgaps}, the spectral subspace
$$\operatorname{Ran}\Pi_p=\operatorname{Span}\left\{\sqrt{\hat{\rho}_{\tau,p,\delta}}u_{1,\delta}(p),\ldots,\sqrt{\hat{\rho}_{\tau,p,\delta}}u_{m,\delta}(p)\right\}$$ corresponds to a uniformly isolated
positive spectral cluster of $T_p$, and
$p\mapsto\Pi_p$ is norm-continuous.
Moreover, Theorem~\ref{thm:radial} gives
\[
\sup_{p\in K_r}\|T_{n,p}-T_p\|\longrightarrow0.
\]

For sufficiently large $n$, let $\Pi_{n,p}$ be the orthogonal
projection in $L^2(\mathbb{B})$ onto the zero-extended spectral subspace
corresponding to
$\lambda_2(\Omega_{n,p}),\ldots,\lambda_{m+1}(\Omega_{n,p})$.
By Lemma~\ref{lem:transfer}, the maps
\[
J_{n,p}
:=\Pi_{n,p}J_p
\bigl(J_p^*\Pi_{n,p}J_p\bigr)^{-1/2}
\]
are well defined, depend continuously on $p$, and satisfy
$J_{n,p}^*J_{n,p}=I_H$, $J_{n,p}J_{n,p}^*=\Pi_{n,p}$.
Set
\[
\psi_{k,n,p}:=J_{n,p}h_k,
\qquad 1\le k\le m,
\]
where these functions are viewed on $\Omega_{n,p}$ and extended
by zero to $\mathbb{B}$. They form an orthonormal basis of the selected
spectral subspace:
\[
\int_{\Omega_{n,p}}\psi_{k,n,p}\psi_{l,n,p}\dd
=\delta_{kl}.
\]
Define $\mathcal C_{n,p}\in\Sym(H)$ by
\[
\langle\mathcal C_{n,p}h_l,h_k\rangle_H
:=\int_{\Omega_{n,p}}
\nabla\psi_{k,n,p}\cdot\nabla\psi_{l,n,p}\dd,
\qquad 1\le k,l\le m.
\]
Thus $\mathcal C_{n,p}$ represents the Dirichlet Laplacian
on this spectral subspace in the basis
$\{\psi_{k,n,p}\}_{k=1}^m$, and the eigenvalues of $\mathcal C_{n,p}$ are precisely $\lambda_2(\Omega_{n,p}),\ldots,\lambda_{m+1}(\Omega_{n,p})$.

Now, in Lemma \ref{lem:transfer},  taking $g(z)=z^{-1}-1$ on a neighborhood of the
selected positive spectral clusters and using Remark \ref{rem:counts} and \eqref{eq:BpLp}, we obtain
\[
\mathcal C_{n,p}
=J_{n,p}^*[g(T_{n,p})\Pi_{n,p}]J_{n,p},\quad
\mathcal B_p+C_*I_H
=J_p^*[g(T_p)\Pi_p]J_p.
\]
Here the functional calculus is restricted to the selected
finite-dimensional spectral subspaces.

Theorem~\ref{thm:radial}, Remark~\ref{rem:counts} and Lemma~\ref{lem:transfer}
therefore imply that $p\mapsto\mathcal C_{n,p}$ is continuous
and allow one sufficiently large $n$, independently of $p$,
to be chosen such that
\begin{equation}\label{eq:lastmatrixerror}
\begin{aligned}
\sup_{p\in K_r}\max_{1\le j\le m+2}
|\lambda_j(\Omega_{n,p})-(C_*+\xi_j(p))|
&<\frac{\tau(m+1)}{16},\\
\sup_{p\in K_r}
\|\mathcal C_{n,p}-(\mathcal B_p+C_*I_H)\|
&<\frac{\tau r}{4}.
\end{aligned}
\end{equation}
Fix this $n$ and set
\[
\Omega_p:=\Omega_{n,p},
\qquad
\psi_{k,p}:=\psi_{k,n,p},
\qquad
\mathcal C_p:=\mathcal C_{n,p}.
\]
The radial functions of $\Omega_p$ depend continuously on $p$
in every $C^\ell$ norm, and $p\mapsto\mathcal C_p$ is continuous. By \eqref{eq:potentialgaps}, \eqref{eq:clusterthresholds},
and \eqref{eq:lastmatrixerror},
\[
\begin{aligned}
\lambda_1(\Omega_p)
&<C_*+\frac{5\tau(m+1)}{16}<b_-,\\
\lambda_2(\Omega_p)
&>C_*+\frac{7\tau(m+1)}{16}>b_-,\\
\lambda_{m+1}(\Omega_p)
&<C_*+\frac{25\tau(m+1)}{16}<b_+,\\
\lambda_{m+2}(\Omega_p)
&>C_*+\frac{31\tau(m+1)}{16}>b_+.
\end{aligned}
\]
Thus \eqref{eq:domaincluster} holds, and the eigenvalues of
$\mathcal C_p$ are precisely
$\lambda_2(\Omega_p),\ldots,\lambda_{m+1}(\Omega_p)$.
Define the continuous map $F:K_r\to\Sym(H)$ by
\begin{equation}\label{eq:Ffinal}
F(p):=\tau^{-1}
\bigl(\mathcal C_p-(C_*+\tau(m+1))I_H\bigr).
\end{equation}
Using \eqref{eq:potentialmatrix} and \eqref{eq:lastmatrixerror},
with $\|\cdot\|$ denoting the operator norm, we obtain
\[
\begin{aligned}
\tau\sup_{p\in K_r}\|F(p)-p\|
&\le\sup_{p\in K_r}
\|\mathcal C_p-(\mathcal B_p+C_*I_H)\|
+\sup_{p\in K_r}
\|\mathcal B_p-\tau((m+1)I_H+p)\|\\
&<\frac{\tau r}{4}+\frac{\tau r}{4}
=\frac{\tau r}{2}.
\end{aligned}
\]

Given any $A\in\Sym(H)$ with $\|A\|\le r/4$, we define $G_A:K_r\to\Sym(H)$ by 
       \[
     G_A(p):=p-F(p)+A.
    \]
    Then for any \(p\in K_r\),
    \[
       \|G_A(p)\|\le\|p-F(p)\|+\|A\|\leq \frac{3}{4}r <r.
    \]
    Thus $G_A(K_r)\subset K_r$. Since the closed
ball \(K_r=\{p\in\Sym(H):\|p\|\le r\}\) is compact and convex in the finite-dimensional
space \(\Sym(H)\), by using
Brouwer's fixed point theorem \cite[Chapter~1]{Deimling},
there exists \(p_*\in K_r\) such that
\[
p_*=G_A(p_*)=p_*-F(p_*)+A,
\]
namely $F(p_*)=A.$  
Consequently,
\[
\mathcal C_{p_*}
=(C_*+\tau(m+1))I_H+\tau A,
\]
which proves \eqref{eq:exacttarget}. The proof is complete.
\end{proof}

\subsection{Proof of the main theorem}\label{sec:finiteproof}

\begin{proof}[\bf Proof of Theorem~\ref{thm:finite}]
Denote $m=\sum_{j=1}^\ell m_j=s_\ell-1\geq 1$.
Proposition~\ref{prop:exactmatrix} provides constants
$r,\tau,C_*$ and a family $\{\Omega_p\}_{p\in K_r}$.
Using the fixed orthonormal basis $\{h_k\}_{k=1}^m$ of $H$, set
\[
H_j:=\operatorname{span}\{h_k:s_{j-1}\le k<s_j\}.
\]
Then $H=\bigoplus_{j=1}^\ell H_j$ orthogonally and
$\dim H_j=m_j$. Define
\begin{equation}\label{eq:blocktarget}
a_j:=\frac{r}{8\ell\sqrt m}j,
\qquad
A:=\bigoplus_{j=1}^\ell a_jI_{H_j}\in\Sym(H).
\end{equation}
Since
\[
0<a_1<\cdots<a_\ell=\frac{r}{8\sqrt m},
\qquad
\|A\|
=a_\ell=\frac{r}{8\sqrt m}\le\frac r8<\frac r4.
\]
Proposition~\ref{prop:exactmatrix} therefore yields $p_*\in K_r$
such that
\[
\mathcal C_{p_*}
=(C_*+\tau(m+1))I_H+\tau A
=\bigoplus_{j=1}^\ell\mu_j^0I_{H_j},
\quad
\mu_j^0:=C_*+\tau(m+1+a_j).
\]
Set $\Omega_0:=\Omega_{p_*}$.
Since $\{\psi_{k,p_*}\}_{k=1}^m$ is an orthonormal basis of
the entire spectral subspace associated with $(b_-,b_+)$,
and $\mathcal C_{p_*}$ is block diagonal, we have
\[
-\Delta_{\Omega_0}^D\psi_{k,p_*}
=\mu_j^0\psi_{k,p_*},
\qquad s_{j-1}\le k<s_j.
\]
Since $\tau>0$ and $a_1<\cdots<a_\ell$,
\eqref{eq:domaincluster} gives
\[
\lambda_1(\Omega_0)<b_-<\mu_1^0<\cdots<\mu_\ell^0
<b_+<\lambda_{m+2}(\Omega_0).
\]
Consequently, for every $1\le j\le\ell$,
\[
\ker(-\Delta_{\Omega_0}^D-\mu_j^0)
=\operatorname{span}
\{\psi_{k,p_*}:s_{j-1}\le k<s_j\}.
\]
Thus
\[
\dim\ker(-\Delta_{\Omega_0}^D-\mu_j^0)=m_j.
\]
Since $\lambda_1(\Omega_0)$ is the only eigenvalue below $b_-$,
\[
\lambda_k(\Omega_0)=\mu_j^0
\qquad\mathrm{for}\qquad
s_{j-1}+1\le k\le s_j,\ 1\le j\le\ell.
\]
Since $\Omega_0$ is a smooth strictly star-shaped domain of the form \eqref{eq:strictstar},
it follows from the
isotopy extension theorem \cite[Chapter~8, Theorem~1.3]{Hirsch} that
$\overline{\Omega_0}$ is $C^\infty$-diffeomorphic to the closed ball. This completes the proof.
\end{proof}

\section{Multiplicity of Dirichlet Eigenvalues on Convex Domains}
\label{sec:convexproof}

In this section, we consider convex domains and prove
Theorem~\ref{thm:convex-all-multiplicities}.
The following argument combines the ground-state log-concavity
estimate of Andrews--Clutterbuck \cite{AC2011} with the
averaged-gradient argument of Barthe--Klartag
\cite[Lemma~14]{BK2020}.
We recall the log-concavity comparison theorem of
Andrews--Clutterbuck in the smooth-potential setting.

\begin{theorem}[Andrews--Clutterbuck
{\cite[Theorem~1.5]{AC2011}}]
Let $\Omega\subset\mathbb{R}^N$ be a bounded convex domain
of diameter $D$. Let $W\in C^\infty(\overline{\Omega})$ and
$\widetilde W\in C^\infty([-D/2,D/2])$, with $\widetilde W$
even. Denote by $\phi_1$ and $\widetilde\phi_1$ the positive
first Dirichlet eigenfunctions of
\[
-\Delta+W \quad\text{in }\Omega,
\qquad
-\frac{d^2}{ds^2}+\widetilde W
\quad\text{in }(-D/2,D/2),
\]
respectively. Suppose that, for all distinct $x,y\in\Omega$,
\[
\left\langle
\nabla W(y)-\nabla W(x),
\frac{y-x}{|y-x|}
\right\rangle
\ge
2\widetilde W'\!\left(\frac{|y-x|}{2}\right).
\]
Then, for all distinct $x,y\in\Omega$,
\begin{equation}\label{eq:AC-log-concavity-comparison}
\left\langle
\nabla\log\phi_1(y)-\nabla\log\phi_1(x),
\frac{y-x}{|y-x|}
\right\rangle
\le
2(\log\widetilde\phi_1)'
\!\left(\frac{|y-x|}{2}\right).
\end{equation}
\end{theorem}

\begin{remark}\label{Sec 6 rem}
Fix $x\in\Omega$ and a unit vector $e\in\mathbb{R}^N$,
and take
\[
W\equiv0,\qquad
\widetilde W\equiv0,\qquad
\widetilde\phi_1(s)=\cos\frac{\pi s}{D}.
\]
For sufficiently small $t>0$, we have $x+te\in\Omega$.
Taking $y=x+te$ in \eqref{eq:AC-log-concavity-comparison}
and dividing by $t$ gives
\[
\left\langle
\frac{\nabla\log\phi_1(x+te)-\nabla\log\phi_1(x)}{t},
e
\right\rangle
\le
-\frac{2\pi}{D}\,
\frac{\tan(\pi t/(2D))}{t}.
\]
Letting $t\downarrow0$, we obtain
\[
\left\langle
\nabla^2\log\phi_1(x)e,e
\right\rangle
\le
-\frac{2\pi}{D}\cdot\frac{\pi}{2D}
=
-\frac{\pi^2}{D^2}.
\]
Since $x$ and $e$ are arbitrary,
\[
\nabla^2\log\phi_1
\le -\frac{\pi^2}{D^2}I_N
\qquad\text{in }\Omega.
\]
\end{remark}

\begin{proof}[\bf Proof of
Theorem~\ref{thm:convex-all-multiplicities}]
We first assume that $\Omega$ has smooth boundary and is uniformly
convex. Let $\{\phi_j\}_{j\ge1}$ be a real orthonormal basis of
Dirichlet eigenfunctions in $L^2(\Omega)$, with $\phi_1>0$ in
$\Omega$. Set $d\mu=\phi_1^2\,dx$ and $V=-2\log\phi_1$,
and introduce
\[
L=\Delta-\nabla V\cdot\nabla,\qquad
\gamma_j=\lambda_j-\lambda_1,\qquad
\kappa=\frac{2\pi^2}{D^2}.
\]
In particular, $\mu(\Omega)=1$ and
$0=\gamma_1<\gamma_2\le\gamma_3\le\cdots$.
Using Remark~\ref{Sec 6 rem}, 
\begin{equation}\label{eq:convex-groundstate-hessian}
\nabla^2V
=-2\nabla^2\log\phi_1
\ge \frac{2\pi^2}{D^2}I_N
=\kappa I_N
\qquad\text{in }\Omega.
\end{equation}

For $f_j=\phi_j/\phi_1$, the boundary regularity and Hopf lemma
imply that $f_j$ extends smoothly to $\overline{\Omega}$.
Moreover,
\[
-Lf_j=\gamma_jf_j\quad\text{in }\Omega,
\qquad
\partial_\nu f_j=0\quad\text{on }\partial\Omega,
\qquad
\int_\Omega f_if_j\,d\mu=\delta_{ij};
\]
see also \cite[Proposition~3.1]{AC2011}.
In particular, $f_1=1$.
Multiplying $-\Delta\phi_1=\lambda_1\phi_1$ by
$\phi_1g^2$ and integrating by parts, we obtain
\begin{equation}\label{Sec 6 equ-a}
\int_\Omega |\nabla(\phi_1g)|^2\,dx
-\lambda_1\int_\Omega\phi_1^2g^2\,dx
=
\int_\Omega|\nabla g|^2\,d\mu,
\quad
g\in C^\infty(\overline{\Omega}).
\end{equation}

Fix an integer $k\ge2$. For $g\in C^\infty(\overline{\Omega})$,
set
\[
c_j=\int_\Omega gf_j\,d\mu,\qquad 1\le j\le k-1,
\qquad
u=\phi_1g-\sum_{j=1}^{k-1}c_j\phi_j.
\]
Then $u\in H^1_0(\Omega)$ and
\[
\int_\Omega u\phi_i\,dx=0,
\qquad 1\le i\le k-1.
\]
By orthonormality,
\[
\int_\Omega u^2\,dx
=
\int_\Omega g^2\,d\mu-\sum_{j=1}^{k-1}c_j^2.
\]
The eigenvalue equation $-\Delta\phi_j=\lambda_j\phi_j$  gives
\[
\int_\Omega
\nabla(\phi_1g)\cdot\nabla\phi_j\,dx
=
\lambda_j\int_\Omega\phi_1g\phi_j\,dx
=
\lambda_jc_j.
\]
Consequently,
\[
\int_\Omega|\nabla u|^2\,dx=\int_\Omega\bigg|\nabla\bigg(\phi_1g-\sum_{j=1}^{k-1}c_j\phi_j\bigg)\bigg|^2dx
=
\int_\Omega|\nabla(\phi_1g)|^2\,dx
-\sum_{j=1}^{k-1}\lambda_jc_j^2.
\]
Then the variational characterization of $\lambda_k$ yields
\[
\lambda_k
\left(
\int_\Omega g^2\,d\mu-\sum_{j=1}^{k-1}c_j^2
\right)
\le
\int_\Omega|\nabla(\phi_1g)|^2\,dx
-\sum_{j=1}^{k-1}\lambda_jc_j^2.
\]
Using \eqref{Sec 6 equ-a} and subtracting
$\lambda_1\int_\Omega u^2\,dx$, we obtain
\[
\gamma_k
\left(
\int_\Omega g^2\,d\mu-\sum_{j=1}^{k-1}c_j^2
\right)
\le
\int_\Omega|\nabla g|^2\,d\mu
-\sum_{j=1}^{k-1}\gamma_jc_j^2
\le
\int_\Omega|\nabla g|^2\,d\mu.
\]
Thus
\begin{equation}\label{eq:convex-weighted-poincare}
\gamma_k
\left[
\int_\Omega g^2\,d\mu
-\sum_{j=1}^{k-1}
\left(\int_\Omega gf_j\,d\mu\right)^2
\right]
\le
\int_\Omega|\nabla g|^2\,d\mu,\quad\forall g\in C^\infty(\overline{\Omega}).
\end{equation}

Let $f$ be a finite linear combination of the functions $f_j$,
and define
\[
\Theta_k(f)
:=
\left(
\int_\Omega f_j\partial_i f\,d\mu
\right)_{\substack{1\le i\le N\\1\le j\le k-1}}
\in\mathbb{R}^{N(k-1)}.
\]
Then
\[
|\Theta_k(f)|^2
=
\sum_{i=1}^N\sum_{j=1}^{k-1}
\left(\int_\Omega f_j\partial_i f\,d\mu\right)^2.
\]

The corresponding weighted Bochner identity is
\[
\frac12L(|\nabla f|^2)
=
|\nabla^2f|^2
+\langle\nabla f,\nabla Lf\rangle
+\nabla^2V(\nabla f,\nabla f).
\]
From here, we claim that
\begin{equation}\label{Sec 6 equ-b}
\int_\Omega L(|\nabla f|^2)\,d\mu=0,
\qquad
\int_\Omega\langle\nabla f,\nabla Lf\rangle\,d\mu
=
-\int_\Omega(Lf)^2\,d\mu,
\end{equation}
and hence
\begin{equation}\label{eq:convex-weighted-bochner}
\int_\Omega(Lf)^2\,d\mu
=
\int_\Omega|\nabla^2f|^2\,d\mu
+
\int_\Omega\nabla^2V(\nabla f,\nabla f)\,d\mu.
\end{equation}

To deal with the singularity of $V$ at the boundary,
we justify \eqref{Sec 6 equ-b} and
\eqref{eq:convex-weighted-bochner} by integration over
\[
\Omega_\varepsilon
:=
\{x\in\Omega:
\operatorname{dist}(x,\partial\Omega)>\varepsilon\}.
\]
Let $\nu_\varepsilon$ be its outward unit normal.
Using
$\phi_1^2Lq=\operatorname{div}(\phi_1^2\nabla q)$,
we obtain
\[
\begin{aligned}
\int_{\Omega_\varepsilon}L(|\nabla f|^2)\,d\mu
&=
\int_{\partial\Omega_\varepsilon}
\phi_1^2\partial_{\nu_\varepsilon}(|\nabla f|^2)\,dS,
\\
\int_{\Omega_\varepsilon}
\langle\nabla f,\nabla Lf\rangle\,d\mu
&=
-\int_{\Omega_\varepsilon}(Lf)^2\,d\mu
+
\int_{\partial\Omega_\varepsilon}
\phi_1^2(Lf)\partial_{\nu_\varepsilon}f\,dS.
\end{aligned}
\]
Since $f,Lf\in C^\infty(\overline{\Omega})$,
$\phi_1^2=O(\varepsilon^2)$ on
$\partial\Omega_\varepsilon$, and
$|\partial\Omega_\varepsilon|\le C$,
both boundary integrals are $O(\varepsilon^2)$ and all
 volume integrals converge as $\varepsilon\downarrow0$. This proves \eqref{Sec 6 equ-b}.

Next, we integrate the weighted Bochner identity over
$\Omega_\varepsilon$. The potentially singular term is
integrable because
\[
\phi_1^2\nabla^2V(\nabla f,\nabla f)
=
2(\nabla\phi_1\cdot\nabla f)^2
-
2\phi_1\nabla^2\phi_1(\nabla f,\nabla f)
\]
extends continuously to $\overline{\Omega}$.
We may therefore let $\varepsilon\downarrow0$ and use
\eqref{Sec 6 equ-b} to obtain
\eqref{eq:convex-weighted-bochner}.

Applying \eqref{eq:convex-weighted-poincare} to
$g=\partial_i f$ and summing over $i$ gives
\[
\int_\Omega|\nabla^2f|^2\,d\mu
\ge
\gamma_k
\left[
\int_\Omega|\nabla f|^2\,d\mu
-
|\Theta_k(f)|^2
\right].
\]
Together with \eqref{eq:convex-groundstate-hessian}
and \eqref{eq:convex-weighted-bochner}, this implies
\begin{equation}\label{eq:convex-averaged-gradient-estimate}
\int_\Omega(Lf)^2\,d\mu
\ge
(\gamma_k+\kappa)
\int_\Omega|\nabla f|^2\,d\mu
-
\gamma_k|\Theta_k(f)|^2.
\end{equation}

We now prove \eqref{eq:convex-quantitative-separation}.
Set
\[
d=N(k-1)+2.
\]
The space
$\operatorname{span}\{f_2,\ldots,f_d\}$
has dimension
$
d-1=N(k-1)+1>N(k-1).
$
Since $\Theta_k$ takes values in $\mathbb{R}^{N(k-1)}$,
there exists a nonzero function
\[
f=\sum_{j=2}^da_jf_j
\qquad\text{such that}\qquad
\Theta_k(f)=0.
\]
Then $-Lf=\sum_{j=2}^d \gamma_ja_jf_j$ and integration by part give
\[
\int_\Omega|\nabla f|^2\,d\mu
=
-\int_\Omega fLf\,d\mu
=
\sum_{j=2}^d\gamma_ja_j^2
>0.
\]
By \eqref{eq:convex-averaged-gradient-estimate},
\begin{align*}
(\gamma_k+\kappa)
\int_\Omega|\nabla f|^2\,d\mu
&\le
\int_\Omega(Lf)^2\,d\mu
=
\sum_{j=2}^d\gamma_j^2a_j^2
\\
&\le
\gamma_d\sum_{j=2}^d\gamma_ja_j^2
=
\gamma_d\int_\Omega|\nabla f|^2\,d\mu.
\end{align*}
Thus
\[
\lambda_d-\lambda_k=\gamma_d-\gamma_k\geq\kappa.
\]
This proves \eqref{eq:convex-quantitative-separation}
for smooth uniformly convex domains.

Now let $\Omega$ be an arbitrary bounded convex domain.
After a translation, assume that $0\in\Omega$.
Choose smooth uniformly convex domains $\Omega_\ell$ and
$\varepsilon_\ell\downarrow0$ such that
\[
(1-\varepsilon_\ell)\Omega
\subset\Omega_\ell
\subset(1+\varepsilon_\ell)\Omega.
\]
Domain monotonicity and scaling give, for every fixed $j$,
\[
(1+\varepsilon_\ell)^{-2}\lambda_j(\Omega)
\le
\lambda_j(\Omega_\ell)
\le
(1-\varepsilon_\ell)^{-2}\lambda_j(\Omega).
\]
Consequently,
\[
\lambda_j(\Omega_\ell)\longrightarrow\lambda_j(\Omega),
\qquad
\operatorname{diam}(\Omega_\ell)
\longrightarrow\operatorname{diam}(\Omega).
\]
For every fixed $k\geq2$, passing to the limit in
\eqref{eq:convex-quantitative-separation} for $\Omega_\ell$
proves the same estimate for $\Omega$.

It remains to prove
\eqref{eq:convex-all-multiplicity-bound}.
Denote
\[
r=\min\{j\ge1:\lambda_j(\Omega)=\lambda_k(\Omega)\}.
\]
Then $r\le k$.
Applying \eqref{eq:convex-quantitative-separation}
with $k$ replaced by $r$, we have
\[
\lambda_{N(r-1)+2}(\Omega)
\ge
\lambda_r(\Omega)+\frac{2\pi^2}{D^2}
>
\lambda_r(\Omega).
\]
It follows that
\[
r+m_k(\Omega)-1\le N(r-1)+1,
\]
and therefore
\[
m_k(\Omega)
\le
(N-1)(r-1)+1
\le
(N-1)(k-1)+1.
\]
This proves \eqref{eq:convex-all-multiplicity-bound}.

Finally, taking $k=2$, we obtain $m_2(\Omega)\le N$, and the multiplicity bound $N$ is sharp,
as shown by $\Omega=\mathbb B$.
Another example is the cube $\Omega=(0,\pi)^N$.
Its Dirichlet eigenvalues are
$k_1^2+\cdots+k_N^2$, where
$k_i\in\{1,2,\ldots\}$.
Thus $\lambda_1=N$ and $\lambda_2=N+3$, the latter having
exactly $N$ independent eigenfunctions, obtained by taking
one $k_i=2$ and all remaining $k_j=1$.
The proof is complete.
\end{proof}

\appendix

\section{Operator Continuity and Spectral Perturbation}
\label{sec:spectral}

Let $K_r$ be defined by Lemma \ref{lem:graphweights}, and let $T_{n,p}$ and $T_p$
be norm-continuous families of compact nonnegative self-adjoint
operators on a real Hilbert space $\mathcal H$.
For each $p\in K_r$, let
$\Sigma_p\subset\sigma(T_p)\cap(0,\infty)$ be a finite spectral
cluster of total multiplicity $\ell\ge1$, and let
$\Pi_p:\mathcal H\to\mathcal H$ be the orthogonal projection onto
\[
\operatorname{Ran}\Pi_p
=
\bigoplus_{\lambda\in\Sigma_p}\ker(T_p-\lambda I).
\]
Assume that $p\mapsto\Pi_p$ is norm-continuous and that
\begin{equation}\label{spec dist}
\delta:=
\inf_{p\in K_r}
\operatorname{dist}
\bigl(\Sigma_p,\{0\}\cup(\sigma(T_p)\setminus\Sigma_p)\bigr)>0.
\end{equation}
Similarly, let $\Pi_{n,p}$ denote the orthogonal spectral
projection of $T_{n,p}$ associated with
\[
\Sigma_{n,p}
:=
\left\{\lambda\in\sigma(T_{n,p}):
\operatorname{dist}(\lambda,\Sigma_p)<\delta/3\right\}.
\]

Given an orthonormal basis
$u_1(p),\ldots,u_\ell(p)$ of $\operatorname{Ran}\Pi_p$,
define its frame map by
\[
J_p:\mathbb R^\ell\longrightarrow\mathcal H,
\qquad
J_px:=\sum_{k=1}^{\ell}x_k u_k(p),
\]
where $x_k$ are the standard coordinates of $x\in\mathbb R^\ell$.
Its adjoint and the projection satisfy
\[
J_p^*f
=
\bigl(\langle f,u_k(p)\rangle_{\mathcal H}\bigr)_{k=1}^{\ell},
\qquad
\Pi_pf
=
\sum_{k=1}^{\ell}
\langle f,u_k(p)\rangle_{\mathcal H}u_k(p).
\]
In particular, $J_p^*J_p=I_\ell$ and $J_pJ_p^*=\Pi_p$.
The frame is called continuous if
$p\mapsto J_p$ is continuous in
$\mathcal L(\mathbb R^\ell,\mathcal H)$.

For a function $g$ holomorphic in a neighborhood of the selected spectrum,
$g(T_p)\Pi_p$ denotes the functional calculus of
$T_p|_{\operatorname{Ran}\Pi_p}$, extended by zero on
$(\operatorname{Ran}\Pi_p)^\perp$; the notation
$g(T_{n,p})\Pi_{n,p}$ is interpreted similarly.
We use complexifications when taking contour integrals or
applying complex-valued functions.
The norms $\|\cdot\|$  denote
the operator  norm.

The following lemma records a uniform version of standard
spectral perturbation results; see
\cite[Sections~I.4.6, I.5.6 and I.6.8, and
Theorem~IV.3.16]{Kato}.

\begin{lemma}\label{lem:transfer}
Under the above assumptions, suppose that
\[
\varepsilon_n:=
\sup_{p\in K_r}\|T_{n,p}-T_p\|\longrightarrow0.
\]
Then, for all sufficiently large $n$,
$p\mapsto\Pi_{n,p}$ is norm-continuous and
\begin{equation}\label{eq:transfer-projection}
\operatorname{rank}\Pi_{n,p}=\ell,
\qquad
\sup_{p\in K_r}\|\Pi_{n,p}-\Pi_p\|\longrightarrow0.
\end{equation}

If $J_p$ is a given continuous orthonormal frame, then
\begin{equation}\label{eq:transported-frame}
J_{n,p}
:=
\Pi_{n,p}J_p
\bigl(J_p^*\Pi_{n,p}J_p\bigr)^{-1/2}
\end{equation}
is well defined for large $n$ and is a continuous orthonormal
frame of $\operatorname{Ran}\Pi_{n,p}$, with
\[
\sup_{p\in K_r}\|J_{n,p}-J_p\|\longrightarrow0.
\]
In this case, for every open $V\subset\mathbb C$ containing
$\bigcup_{p\in K_r}\Sigma_p$ and every holomorphic
$g:V\to\mathbb C$, the matrices below are well defined and
continuous in $p$ for all sufficiently large $n$, and
\begin{equation}\label{eq:spectral-matrix-transfer}
\sup_{p\in K_r}
\left\|
J_{n,p}^*\bigl[g(T_{n,p})\Pi_{n,p}\bigr]J_{n,p}
-
J_p^*\bigl[g(T_p)\Pi_p\bigr]J_p
\right\|
\longrightarrow0.
\end{equation}
\end{lemma}

	\begin{proof}
By the stability of separated spectral subspaces
\cite[Theorem~IV.3.16]{Kato}, near each $p_0\in K_r$ we may use
a fixed positively oriented contour $\Gamma$ enclosing the
selected clusters of both $T_p$ and $T_{n,p}$, and no other
spectrum, for all sufficiently large $n$.
Here we use the continuity of $T_p$ and $\Pi_p$ and the gap
in \eqref{spec dist} to follow the same cluster.
The theorem also gives continuity of $\Pi_{n,p}$ and preservation
of its rank. Compactness of $K_r$ reduces the argument to finitely
many such contours, on which the resolvents are uniformly bounded.

The Riesz projection formula
\cite[Section~III.6.4-(6.19)]{Kato} and the resolvent identity give
\[
\Pi_{n,p}-\Pi_p
=\frac{1}{2\pi i}\int_\Gamma
(zI-T_{n,p})^{-1}(T_{n,p}-T_p)(zI-T_p)^{-1}\,dz.
\]
Consequently,
\begin{equation}\label{eq:Pi estimate}
\sup_{p\in K_r}\|\Pi_{n,p}-\Pi_p\|
\le C\varepsilon_n\longrightarrow0,
\end{equation}
where $C$ is independent of $p$ and $n$.
This proves \eqref{eq:transfer-projection}.

To construct the frame, set $M_{n,p}:=J_p^*\Pi_{n,p}J_p$.
This is the Gram matrix of the projected vectors
$\Pi_{n,p}u_1(p),\ldots,\Pi_{n,p}u_\ell(p)$, and $J_p^*\Pi_pJ_p=I_{\ell}$ gives
\[
\|M_{n,p}-I_\ell\|
\le\|\Pi_{n,p}-\Pi_p\|\le C\varepsilon_n.
\]
Thus $M_{n,p}$ is positive definite for large $n$, and
$M_{n,p}^{-1/2}\to I_\ell$ uniformly in $p$.
Multiplication by $M_{n,p}^{-1/2}$ orthonormalizes the projected
vectors, since
\[
J_{n,p}^*J_{n,p}
=M_{n,p}^{-1/2}M_{n,p}M_{n,p}^{-1/2}=I_\ell.
\]
Since $\operatorname{rank}\Pi_{n,p}=\ell$, these vectors span
$\operatorname{Ran}\Pi_{n,p}$. Using $\Pi_pJ_p=J_p$ and \eqref{eq:transported-frame}, we obtain
\begin{align*}
\|J_{n,p}-J_p\|&=\|\Pi_{n,p}J_pM_{n,p}^{-1/2}-\Pi_{p}J_pI_{\ell}\|\\
&\le C_1\|\Pi_{n,p}-\Pi_p\|
+C_2\|M_{n,p}^{-1/2}-I_\ell\|\longrightarrow0
\end{align*}
uniformly in $p$.

Finally, for the given $V$ and $g$, choose the local contours
and their enclosed regions inside $V$. The holomorphic functional
calculus \cite[Section~I.5.6-(5.47)]{Kato} gives
\[
g(T_p)\Pi_p
=\frac{1}{2\pi i}\int_\Gamma g(z)(zI-T_p)^{-1}\,dz,
\]
and the analogous formula for $g(T_{n,p})\Pi_{n,p}$. The same resolvent
estimate, again using a finite cover of $K_r$, yields
\begin{equation}\label{eq:g projection estimate}
\sup_{p\in K_r}
\|g(T_{n,p})\Pi_{n,p}-g(T_p)\Pi_p\|
\le C_g\varepsilon_n.
\end{equation}
These operators are norm-continuous in $p$ and uniformly bounded.
Combining this estimate with the uniform convergence of the frames
gives
\[
\begin{aligned}
&\left\|J_{n,p}^*\bigl[g(T_{n,p})\Pi_{n,p}\bigr]J_{n,p}
-J_p^*\bigl[g(T_p)\Pi_p\bigr]J_p\right\|\\
&\hspace{2em}\le C_g\bigl(\varepsilon_n+\|J_{n,p}-J_p\|\bigr)
\longrightarrow0
\end{aligned}
\]
uniformly in $p$, this proves
\eqref{eq:spectral-matrix-transfer}.
\end{proof}

\begin{remark}\label{rem:counts}\
\begin{itemize}
\item[(1)] In Lemma \ref{lem:transfer}, if we choose 
		\(T_{n,p}=(A+1)^{-1}\) and \(g(z)=z^{-1}-1\), then by $\Pi_{n,p}J_{n,p}=J_{n,p}$, we have
		\[g(T_{n,p})\Pi_{n,p}=A\Pi_{n,p},\quad J_{n,p}^*[g(T_{n,p})\Pi_{n,p}]J_{n,p}=J_{n,p}^*AJ_{n,p}.\]
\item[(2)]
The min--max principle gives
\[
\sup_{p\in K_r}
|\lambda_j(T_{n,p})-\lambda_j(T_p)|
\le\varepsilon_n.
\]
If $\lambda_j(T_p)>0$ for every $p\in K_r$, the continuity and compactness
imply
$c_j:=\min_{p\in K_r}\lambda_j(T_p)>0$.
Then for all sufficiently large $n$,
$\lambda_j(T_{n,p})\ge c_j/2$ holds uniformly in $p$, hence
\[
\sup_{p\in K_r}
\left|
\bigl(\lambda_j(T_{n,p})^{-1}-1\bigr)
-
\bigl(\lambda_j(T_p)^{-1}-1\bigr)
\right|
\le\frac{2\varepsilon_n}{c_j^2}
\longrightarrow0.
\]
\end{itemize}
\end{remark}

\bigskip

\subsection*{Acknowledgments}  Z. Chen is supported by National Key R\&D Program of China (No. 2023YFA1010002) and NSFC (No. 12222109). 
S. Zhang is supported by the Postdoctoral Fellowship Program and China Postdoctoral Science Foundation under Grant Numbers BX20250062 and 2026M793381.
The AI tools were used during the preparation of this manuscript. The authors assume full responsibility for the
mathematical validity, accuracy, and integrity of the proofs presented in the manuscript.

\end{document}